\documentclass[a4paper,12pt]{article}
\usepackage[latin1]{inputenc}
\usepackage{amsmath,amsrefs,amssymb,amsthm}
\usepackage{graphicx}
\usepackage[bookmarks=true,pdfstartview={FitH}]{hyperref}
\usepackage{enumerate,enumitem}
\usepackage{mathrsfs}
\usepackage{epstopdf}
\usepackage{fullpage}
\usepackage[small]{caption}

\newtheorem{defi}{Definition}[section]
\newtheorem{lem}[defi]{Lemma}
\newtheorem{thm}[defi]{Theorem}
\newtheorem{cor}[defi]{Corollary}
\newtheorem{prop}[defi]{Proposition}
\theoremstyle{remark}
\newtheorem{rem}[defi]{Remark}
\numberwithin{equation}{section}
\numberwithin{figure}{section}

\def\cB{\ensuremath{\mathcal{B}}}
\def\oB{\ensuremath{\mathcal{U}}}
\def\C{\ensuremath{\mathcal{C}}}
\def\Cc{\ensuremath{\mathscr{C}}}
\def\diam{\ensuremath{{\rm diam}\;}}

\def\Dd{\ensuremath{\mathfrak{D}}}
\def\dH{\ensuremath{d_{\rm H}}}
\def\dP{\ensuremath{d_{\rm P}}}
\def\dGP{\ensuremath{d_{\rm GP}}}
\def\dMGP{\ensuremath{d_{\rm mGP}}}
\def\dGHP{\ensuremath{d_{\rm GHP}}}
\def\dust{\ensuremath{\mathcal{M}_{\rm{dust}}}}
\def\E{\ensuremath{\mathbb{E}}}
\def\ep{\varepsilon}
\def\e{\mathrm{e}}
\def\F{\ensuremath{\mathcal{F}}}

\def\Ind{\ensuremath{\mathbf{1}}}
\def\p{\ensuremath{\mathcal{P}}}
\def\P{\ensuremath{\mathbb{P}}}
\def\Q{\ensuremath{\mathbb{Q}}}
\def\M{\ensuremath{\mathbb{M}}}
\def\Mm{\ensuremath{\mathcal{M}}}

\def\Mc{\ensuremath{\mathbb{M}_{\rm c}}}
\def\bM{\ensuremath{\mathbf{M}}}
\def\MCDI{\ensuremath{\mathcal{M}_{\rm{CDI}}}}
\def\N{\ensuremath{\mathbb{N}}}
\def\nd{\ensuremath{\mathcal{M}_{\rm nd}}}
\def\R{\ensuremath{\mathbb{R}}}
\def\S{\ensuremath{\mathcal{S}}}

\def\sa{\ensuremath{[a]}}
\def\sb{\ensuremath{[b]}}
\def\sk{\ensuremath{[k]}}

\def\sn{\ensuremath{[n]}}
\def\sell{\ensuremath{[\ell]}}

\def\supp{\ensuremath{{\rm supp}\;}}
\def\Thetaext{\ensuremath{\Theta^{\rm ext}}}
\def\U{\ensuremath{\mathbb{U}}}
\def\UUerg{\ensuremath{\mathcal{U}^{\rm erg}}}
\def\X{\ensuremath{\mathcal{X}}}

\newcommand{\abs}[1]{\mathop {\left|{ #1}\right|}}
\newcommand{\reff}[1]{(\ref{#1})}
\newcommand{\expp}[1]{\mathop {\mathrm{e}^{ #1}}}
\newcommand{\ew}[1]{\mathop {[\![{ #1}]\!]}}
\newcommand{\es}[1]{\mathop {[{ #1}]}}
\newcommand{\I}[1]{\mathop {\mathbf{1}{\left\{ #1\right\}}}}
\newcommand{\hparagraph}[1]{\paragraph{#1}\mbox{}\\}

\begin{document}
\title{Pathwise construction of tree-valued Fleming-Viot processes}
\author{
Stephan Gufler\thanks{Technion, Faculty of Industrial Engineering and Management, Haifa 3200003, Israel, \texttt{stephan.gufler@gmx.net}}
}
\maketitle
\label{Pathw}
\begin{abstract}
In a random complete and separable metric space that we call the lookdown space, we encode the genealogical distances between all individuals ever alive in a
lookdown model with simultaneous multiple reproduction events.
We construct families of probability measures on the lookdown space and on an extension of it that allows to include the case with dust.
From this construction, we read off the tree-valued $\Xi$-Fleming-Viot processes and deduce path properties.
For instance, these processes usually have a.\,s.\ càdlàg paths with jumps at the times of large reproduction events. In the case of coming down from infinity, the construction on the lookdown space also allows to read off a process with values in the space of measure-preserving isometry classes of compact metric measure spaces, endowed with the Gromov-Hausdorff-Prohorov metric. This process has a.\,s.\ càdlàg paths with additional jumps at the extinction times of parts of the population.
\end{abstract}
{\small
\emph{Keywords:} Lookdown model, tree-valued Fleming-Viot process, evolving coalescent, $\Xi$-coales\-cent, (marked) metric measure space, (marked) Gromov-weak topology, Gromov-Hausdorff-Prohorov topology.\\
\emph{AMS MSC 2010:} Primary 60K35, Secondary 60J25, 60G09, 92D10.}
\setcounter{tocdepth}{2}
{\scriptsize \tableofcontents }

\section{Introduction}
Similarly to the measure-valued Fleming-Viot process that is a model for the evolution of the type distribution in a large neutral haploid population, a tree-valued Fleming-Viot process models the evolution of the distribution of the genealogical distances between randomly sampled individuals.
The (neutral) tree-valued Fleming-Viot process is introduced in Greven, Pfaffelhuber, and Winter \cite{GPW13} and generalized in \cite{Sampl} to the setting with simultaneous multiple reproduction events.
The lookdown model of Donnelly and Kurtz \cites{DK96,DK99} provides a pathwise construction of the measure-valued Fleming-Viot process and more general measure-valued processes. In this article, we give a pathwise construction of the tree-valued Fleming-Viot process from the lookdown model.

Let us sketch the lookdown model that we state in more detail
in Section~\ref{Pathw:sec:ld-space}.
The time axis is $\R_+$. In the population model, there are countably infinitely many levels which are labeled by $\N$. Each level is occupied by one particle at each time. As time evolves, the particles undergo reproduction events in which particles can increase their levels.
We call a particle at a fixed instant in time an individual. We identify each element $(t,i)$ of $\R_+\times\N$ with the individual on level $i$ at time $t$.
From the genealogy that is determined by the reproduction events and from given genealogical distances between the individuals at time zero, we define the semi-metric $\rho$ on $\R_+\times\N$ of
the genealogical distances between all individuals.
We speak of the case with dust if each particle reproduces at finite rate.
In the general case, only the rate at which a particle reproduces and has offspring on a given level is finite.
In the case without dust, we introduce the lookdown space $(Z,\rho)$ as the metric completion of $(\R_+\times\N,\rho)$.
We allow for simultaneous multiple reproduction events so that we can obtain any $\Xi$-coalescent as the genealogy at a fixed time \cites{Schw00,MS01,Pitman99,Sagitov99,DK99,BBMST09}.

In Section~\ref{Pathw:sec:ld-sampl}, we state the central results in this article. Theorem~\ref{Pathw:thm:ld-nd} asserts that a.\,s.\ in the case without dust, the uniform measures
$\mu^n_t=n^{-1}\sum_{i=1}^n\delta_{(t,i)}$
on the individuals on the first $n$ levels at time $t$ converge uniformly in compact time intervals to some probability measures $(\mu_t,t\in\R_+)$ in the Prohorov metric $\dP^Z$ over the lookdown space,
\begin{equation}
\label{Pathw:eq:conv-unif}
\lim_{n\to\infty}\sup_{t\in[0,T]}\dP^Z(\mu^n_t,\mu_t)\quad\text{a.\,s.\ for all }T\in\R_+.
\end{equation}

We recall that a metric measure space $(X,r,\mu)$ is a triple that consists of a complete and separable metric space $(X,r)$ and a probability measure $\mu$ on the Borel sigma algebra on $(X,r)$. The Gromov-Prohorov distance between two metric measure spaces $(X,r,\mu)$ and $(X',r',\mu')$ is defined as
\begin{equation*}
\dGP((X,r,\mu),(X',r',\mu'))=\inf_{Y,\varphi,\varphi'}\dP^Y(\varphi(\mu),\varphi'(\mu'))
\end{equation*}
where the infimum is over all isometric embeddings $\varphi:X\to Y$, $\varphi':X'\to Y$ into complete and separable metric spaces $Y$, the Prohorov metric over $Y$ is denoted by $\dP^Y$, and pushforward measures are written as $\varphi(\mu)=\mu\circ\varphi^{-1}$. Two metric measure spaces are called isomorphic if their Gromov-Prohorov distance is zero, or equivalently, if there is a measure-preserving isometry between the closed supports of the measures. The Gromov-Prohorov distance is a complete and separable metric on the space $\M$ of isomorphy classes of metric measure spaces, it induces the Gromov-weak topology in which metric measure spaces converge if and only if the distributions of the matrices of the distances between iid samples (the so-called distance matrix distributions) converge weakly. For the theory of metric measure spaces, we refer to Greven, Pfaffelhuber, and Winter \cite{GPW09} and Gromov \cite{Gromov}.

In Section~\ref{Pathw:sec:proc}, we read off the tree-valued $\Xi$-Fleming-Viot processes.
While the lookdown model is used in \cite{Sampl} to characterize only versions of the tree-valued $\Xi$-Fleming-Viot processes (see Remark 4.4 in \cite{Sampl}), we obtain the whole paths in the present article.
Other than in \cite{Sampl}, we do not use ultrametricity of the initial state for the techniques in the present article. Therefore, we speak for instance of an $\M$-valued $\Xi$-Fleming-Viot process when the initial state not necessarily corresponds to an ultrametric tree. Such an $\M$-valued $\Xi$-Fleming-Viot process is given in the case without dust by $(\chi_t,t\in\R_+)$, where $\chi_t$ is the isomorphy class of the metric measure space $(Z,\rho,\mu_t)$.
We stress that to construct a tree-valued Fleming-Viot process, we first show a.\,s.\ convergence of probability measures in the Prohorov metric on the lookdown space, as in \reff{Pathw:eq:conv-unif}. Thereafter, we take isomorphy classes to obtain pathwise a tree-valued process.
By contrast, in \cite{GPW13}, first finite population models are considered and isomorphy classes are taken to obtain tree-valued processes whose convergence in distribution in the Gromov-Prohorov metric is then shown. The limit process is then characterized by a well-posed martingale problem.
From our approach, we can also deduce convergence in distribution of tree-valued processes that describe finite population models, see Remark \ref{Pathw:rem:tvC}.

From our approach on the lookdown space, it follows readily that tree-valued $\Xi$-Fleming-Viot processes have a.\,s.\ càdlàg paths with jumps at the times of large reproduction events (except for some settings in which there is no right-continuity at initial time). 
In particular, we retrieve the result from \cite{GPW13}
that paths are a.\,s.\ continuous in the Gromov-weak topology in the case with only binary reproduction events (which is the case associated with the Kingman coalescent). The Gromov-weak topology emphasizes the typical genealogical distances in a sample from the population.

As initially suggested to the author by G.\ Kersting and A.\ Wakolbinger,
we also consider a process whose state space is endowed with a stronger topology, the Gromov-Hausdorff-Prohorov topology, which highlights also the overall structure of the population.
This process has jumps already in the Kingman case, namely at the times when the shape of the whole genealogical tree changes as all descendants of an ancestor die out (Theorem~\ref{Pathw:thm:ld-CDI} and Proposition~\ref{Pathw:thm:proc-nd-Mvec}). We call this process a tree-valued evolving $\Xi$-coalescent, it can be defined in the case of coming down from infinity which is a subcase of the case without dust.

An ($\R_+$-)marked metric measure space is a triple $(X,r,m)$ that consists of a complete and separable metric space $(X,r)$ and a probability measure $m$ on the Borel sigma algebra on the product space $X\times\R_+$. The space of isomorphy classes of marked metric measure spaces is introduced in Depperschmidt, Greven, and Pfaffelhuber \cite{DGP11}, we recall basic facts in the beginning of Section~\ref{Pathw:sec:proc-gen}.
Tree-valued $\Xi$-Fleming-Viot processes can be defined as processes with values in the space of isomorphy classes of marked metric measure spaces to include the case with dust, as shown in \cite{Sampl}. To give a pathwise construction, we define in Section~\ref{Pathw:sec:ld-dec} for each individual $(t,i)$ a parent $z(t,i)$, and we introduce an extended lookdown space $(\hat Z,\rho)$ which also includes the parents of the individuals at time zero. We denote by $v_t(i)$ the genealogical distance between the individual $(t,i)$ and its parent, and we consider in Section~\ref{Pathw:sec:ld-sampl-dust} the probability measures
\[m^n_t=\frac{1}{n}\sum_{i=1}^n\delta_{(z(t,i),v_t(i))}\]
on $\hat Z\times\R_+$. By Theorem~\ref{Pathw:thm:ld-dust}, these measures $m^n_t$ converge in the Prohorov metric uniformly for $t$ in compact time intervals to probability measures $(m_t,t\in\R_+)$.
In Section~\ref{Pathw:sec:proc-gen}, we obtain the tree-valued $\Xi$-Fleming-Viot process as the process $(\hat\chi_t,t\in\R_+)$, where $\hat\chi_t$ is the isomorphy class of the marked metric measure space $(\hat Z,\rho,m_t)$. In the case without dust, $z(t,i)=(t,i)$ for all individuals $(t,i)$, which yields consistency with the construction in the case without dust.

Besides the space of isomorphy classes of marked metric measure spaces, another possible state space for tree-valued Fleming-Viot processes (in the case with or without dust) is a space of matrix distributions. Here the 
state $\xi_t$ at time $t$ is obtained from the marked metric measure space $(\hat Z,\rho,m_t)$ as follows: We sample an $m_t$-iid sequence $(x(i),v(i))_{i\in\N}$ from $\hat Z\times\R_+$ and consider the infinite matrix $(\rho'(i,j))_{i,j\in\N}$ given by
\begin{equation}
\label{eq:dmd-intr}
\rho'(i,j)=(v(i)+\rho(x(i),x(j))+v(j))\I{i\neq j}.
\end{equation}
Then we define $\xi_t$ as the conditional distribution of $\rho'$ given $(\eta,r_0,v_0)$. Intuitively, this means we condition on the marked metric measure space $(\hat Z,\rho,m_t)$. In equation \reff{eq:dmd-intr}, we obtain $\rho'(i,j)$ by sampling parents $x(i)$, $x(j)$ to whose mutual distance we add their respective distances $v(i)$, $v(j)$ to their descendants at time $t$. We consider the process $(\xi_t,t\in\R_+)$ in the end of Section \ref{Pathw:sec:proc-gen}.

In general, we work in one-sided time. In this way, we obtain the path regularity of the processes under consideration for arbitrary initial states, which is applied in \cite{Conv}. Complementing
the results on convergence to equilibrium in \cite{Sampl}*{Section 9} and \cite{GPW13}*{Theorem 3}, we also show in Section~\ref{Pathw:sec:proc-nd} that the tree-valued evolving $\Xi$-coalescent, started from any initial state, converges to a unique equilibrium, and we define a stationary tree-valued evolving $\Xi$-coalescent in two-sided time. We remark that in the Kingman case, the restriction of the lookdown space in two-sided time to the closure of the set of individuals at a fixed time $t$, endowed with the probability measure $\mu_t$, equals the a.\,s.\ compact metric measure space associated with the Kingman coalescent that is studied by Evans \cite{E00}.

We defer the proofs of the central theorems from Section~\ref{Pathw:sec:ld-sampl} to the second part of the article whose organization is outlined in Section~\ref{Pathw:sec:outline-proofs}.

Now we discuss more relations to the literature.
The lookdown graph of Pfaffelhuber and Wakolbinger \cite{PW06} can be viewed
as a semi-metric space whose completion is a lookdown space in two-sided time.
A lookdown construction of the measure-valued $\Xi$-Fleming-Viot process is given by Birkner et al.\ \cite{BBMST09}.
Véber and Wakolbinger \cite{VW15} give a lookdown construction of measure-valued spatial $\Lambda$-Fleming-Viot processes with dust using a skeleton structure.
To construct the probability measures $(\mu_t,t\in\R_+)$ on the lookdown space in the case without dust, we use the flow of partitions for which we refer to Foucart \cite{Foucart12} and Labbé \cite{Lab12}. A related description of evolving genealogies is the flow of bridges of Bertoin and Le Gall \cite{BLG03}.

Pfaffelhuber, Wakolbinger, and Weisshaupt \cite{PWW11} and Dahmer, Knobloch, and Wakolbinger \cite{DKW14} study the compensated total tree length of the evolving Kingman coalescent as a stochastic process with jumps, using also the lookdown model.
The times of these jumps correspond to the extinction times of parts of the population.
Functionals of evolving coalescents such as the external length have been studied in several works, see for example
\cites{KSW14,DK14}.

For the coming down from infinity property in the setting with simultaneous multiple reproduction events, see e.\,g.\ \cites{Schw00,HM12,BBMST09,Foucart12}.
By methods which differ from those in the present article, it is also shown in \cite{DGP13} that a.\,s., the states of the tree-valued Fleming-Viot process are non-atomic in the Kingman case.
We also mention the work of Athreya, Löhr, and Winter \cite{ALW16} where in particular the Gromov-weak topology and the Gromov-Hausdorff-Prohorov topology are compared.
Marked metric measure spaces are applied by Depperschmidt, Greven, and Pfaffelhuber \cites{DGP12,DGP13} to construct the tree-valued Fleming-Viot process with mutation and selection.

\section{The lookdown space}
\label{Pathw:sec:ld-space}
We write $\R_+=[0,\infty)$, $\N=\{1,2,\ldots\}$, and we denote the set of partitions of $\N$ by $\p$\label{Pathw:not:p}.
For $n\in\N$, we write $\sn=\{1,\ldots,n\}$ and we denote the set of partitions of $\sn$ by $\p_n$. We define the restriction $\gamma_n$ from $\p$ to $\p_n$,\label{Pathw:not:gamma-p}
$\gamma_n(\pi)=\{B\cap[n]:B\in\pi\}\setminus\{\emptyset\}$.
We endow $\p_n$ with the discrete topology, and $\p$ with the topology induced by the restriction maps.

Let us first repeat the lookdown model from \cite{Sampl}*{Section 5.1} which is determined by the genealogy at time $0$ and a point measure that encodes the reproduction events.

In the population model, there are countably infinitely many levels which are labeled by $\N$. The time axis is $\R_+$, and each level is occupied by one particle at each time. To encode the reproduction events that the particles undergo, we use a simple point measure $\eta$\label{Pathw:not:eta-det} on $(0,\infty)\times\p$ with
\begin{equation}
\label{Pathw:eq:ass-eta}
\eta((0,T]\times\p^n)<\infty\quad\text{for all }n\in\N\text{ and }T\in(0,\infty),
\end{equation}
where $\p^n$\label{Pathw:not:pn} denotes the subset of those partitions of $\N$ in which not all of the first $n$ integers are in different blocks, that is,
\begin{equation}
\label{Pathw:eq:def-pn}
\p^n=\{\pi\in\p: \gamma_n(\pi)\neq\{\{1\},\ldots,\{n\}\}\}.
\end{equation}

For a partition $\pi\in\p$ and $i\in\N$, we denote by $B_i(\pi)$\label{Pathw:not:Bi} the $i$-th block of $\pi$ when the blocks are ordered increasingly according to their smallest elements.
Each point $(t,\pi)$ of $\eta$ is interpreted as a reproduction event as follows.
At time $t-$, the particles on levels $i\in\N$ with $i>\#\pi$ are removed. Then, for each $i\in[\#\pi]$, the particle that was on level $i$ at time $t-$ is on level $\min B_i(\pi)$ at time $t$ and has offspring on all other levels in $B_i(\pi)$.
In this way, the level of a particle is non-decreasing as time evolves.
For each $n\in\N$, only finitely many particles in bounded time intervals are pushed away from one of the first $n$ levels by condition~\reff{Pathw:eq:ass-eta}.

We consider not only the process that describes the genealogical distances between the individuals at each fixed time as in \cite{Sampl}, but we are interested in the genealogical distances between all individuals which we describe by a complete and separable metric space, the lookdown space.
We define an individual as a particle at a fixed instant in time. We identify each element $(t,i)$\label{Pathw:not:ti} of $\R_+\times\N$ with the individual on level $i$ at time $t$. For $s\in[0,t]$, we denote by $A_s(t,i)$\label{Pathw:not:ancestor} the level of the ancestor of the individual $(t,i)$ such that the maps $s\mapsto A_s(t,i)$ and $t\mapsto A_s(t,i)$ are càdlàg. Let $\rho_0$ be a semi-metric on $\N$.
We define the genealogical distance between the individuals on levels $i$ and $j$ at time $0$ by
\[\rho((0,i),(0,j))=\rho_0(i,j).\]
More generally, we define the genealogical distance between individuals $(t,i), (u,j)\in\R_+\times\N$ by
\[\rho((t,i),(u,j))=
\left\{
\begin{aligned}
&t+u-2\sup\{s\leq t\wedge u: A_s(t,i)=A_s(u,j)\}\quad\text{if }A_0(t,i)=A_0(u,j)\\
&t+u+\rho_0(A_0(t,i),A_0(u,j))\quad\text{else.}
\end{aligned}
\right.\]
The genealogical distance $\rho((t,i),(u,j))$ can be seen as the sum of the distances to the most recent common ancestor of $(t,i)$ and $(t,j)$ if these individuals have a common ancestor after time zero.
Else it is the genealogical distance of their ancestors at time zero, augmented by the times at which the individuals live. The distinction between these two cases is needed as we work in one-sided time.

The distance $\rho$ is a semi-metric on $\R_+\times\N$ (offspring individuals from the same parent have genealogical distance zero at the time of the reproduction event). We identify individuals with genealogical distance zero, and we take the metric completion. We call the resulting metric space $(Z,\rho)$\label{Pathw:not:Z} the lookdown space associated with $\eta$ and $\rho_0$. In slight abuse of notation, we refer by $(t,i)\in\R_+\times\N$ also to the element of the metric space after the identification of elements with $\rho$-distance zero, in this sense we also assume $\R_+\times\N\subset Z$.

For $t\in\R_+$, we define a semi-metric $\rho_t$ on $\N$ by
\begin{equation}
\label{Pathw:eq:def-rhot}
\rho_t(i,j)=\rho((t,i),(t,j)),\quad i,j\in\N.
\end{equation}
Then $\rho_t$ describes the genealogical distances between the particles at fixed times, and the process $(\rho_t,t\in\R_+)$ is the process that is denoted in the same way in Section 5 of \cite{Sampl}.

The remainder of this section is organized as follows. In Subsection~\ref{Pathw:sec:ld-dec}, we replace $\rho_0$ with a decomposed semi-metric and we enlarge the lookdown space by parents of the individuals at time zero.
In Subsection~\ref{Pathw:sec:ld-ext}, we consider the two ways in which particles can die and we define extinction times for parts of the population.
The construction is randomized in Subsection~\ref{Pathw:sec:ld-Xi} where $\eta$ becomes a Poisson random measure.

\subsection{Parents and decomposed genealogical distances}
\label{Pathw:sec:ld-dec}
We will use the contents of this section to include the case with dust.

Let $r_0$ be a semi-metric on $\N$ and $v_0=(v_0(i))_{i\in\N}\in\R_+^\N$ such that $(r_0,v_0)$ satisfies
\[\rho_0(i,j)=(r_0(i,j)+v_0(i)+v_0(j))\I{i\neq j}\]
for all $i,j\in\N$. Then $(r_0,v_0)$ is a decomposition of the semi-metric $\rho_0$ in the sense of \cite{Sampl}*{Section 2}. The trivial decomposition $(r_0,v_0)=(\rho_0,0)$ always exists.

For each $(t,i)\in\R_+\times\N$, we define the quantity $v_t(i)$\label{Pathw:not:vti} as in \cite{Sampl}*{Section 6.1}: For $j\in\N$, let
$\p(j)=\{\pi\in\p:\{j\}\notin\pi\}$
be the set of partitions of $\N$ in which $j$ does not form a singleton block. If $\eta(\{s\}\times\p(A_s(t,i)))>0$ for some $s\in(0,t]$, then we set
\[v_t(i)=t-\sup\{s\in(0,t]: \eta(\{s\}\times\p(A_s(t,i)))>0\},\]
else we set
\[v_t(i)=t+v_0(A_0(t,i)).\]
The quantity $v_t(i)$ is the time back from the individual $(t,i)$ until the ancestral lineage is involved in a reproduction event in which it belongs to a non-singleton block, if there is such an event, else $v_t(i)$ is defined from $v_0$.
\begin{rem}
For $t\in\R_+$, if $\rho_t$ is a semi-ultrametric (that is, $\max\{\rho_t(j,k),\rho_t(k,\ell)\}\geq\rho_t(j,\ell)$ for all $j,k,\ell\in\N$) and the condition
\begin{equation}
\label{eq:ext-br}
v_t(i)=\tfrac{1}{2}\inf_{j\in\N\setminus\{i\}}\rho_t(i,j)
\end{equation}
is satisfied for some $i\in\N$, then $v_t(i)$ is the length of the external branch that ends in the individual $(t,i)$ in the genealogical tree at time $t$. See Remarks \ref{Pathw:rem:parent-sv} and \ref{rem:ass-v} for more details.
\end{rem}

Now we enlarge the set of individuals to
the disjoint union $(\R_+\times\N)\sqcup\N$. We call each element $i$ of $\N\subset(\R_+\times\N)\sqcup\N$ the parent of the individual $(0,i)$. We extend the semi-metric $\rho$ to $(\R_+\times\N)\sqcup\N$ by
\begin{align*}
\rho(i,j)&=r_0(i,j)\quad\text{for }i,j\in\N\\
\text{and}\quad\rho((t,i),j)&=t+v_0(A_0(t,i))+r_0(A_0(t,i),j)\quad\text{for }(t,i)\in\R_+\times\N,j\in\N.
\end{align*}
That is, the distance between the parents of the individuals $(0,i)$ and $(0,j)$ is given by $r_0(i,j)$.
Furthermore, we define for each individual $(t,i)\in\R_+\times\N$ the parent $z(t,i)$\label{Pathw:not:zti} as the individual $(t-v_t(i),A_{t-v_t(i)}(t,i))$ if $v_t(i)<t$, else we set $z(t,i)=A_0(t,i)$. Then $v_t(i)$ equals the genealogical distance between the individual $(t,i)$ and its parent.

We identify the elements of $(\R_+\times\N)\sqcup\N$ with distance zero and take the metric completion. We call the resulting metric space the extended lookdown space associated with $\eta$ and $(r_0,v_0)$, and we denote it by $(\hat Z,\rho)$\label{Pathw:not:hatZ}. Here we write again $\hat Z\supset(\R_+\times\N)\sqcup\N$ in slight abuse of notation. Note that the lookdown space $(Z,\rho)$ associated with $\eta$ and $\rho_0$ is contained in $(\hat Z,\rho)$ as a subspace. Figure \ref{Pathw:fig:ld-space-dust} below shows an extended lookdown space.
\begin{rem}
\label{Pathw:rem:Z-consistence}
If $v_0=0$, then $\rho_0=r_0$ and the individuals at time zero in the extended lookdown space $(\hat Z,\rho)$ are identified with their parents as $\rho((0,i),z(0,i))=0$ for all $i\in\N$. In this case, $(\hat Z,\rho)$ is equal to the lookdown space $(Z,\rho)$ associated with $\eta$ and $r_0$ from the beginning of Section~\ref{Pathw:sec:ld-space}.
\end{rem}
\begin{rem}[Relation to the decomposed genealogical distances in \cite{Sampl}]
\label{Pathw:rem:rho-r-v}
We denote the genealogical distances between the parents of individuals $(t,i),(t,j)\in\R_+\times\N$ by
\begin{equation*}
\label{Pathw:eq:rt}
r_t(i,j)=\rho(z(t,i),z(t,j)).
\end{equation*}
For $t=0$, this is consistent with the definition of $r_0$ above as $z(0,i)=i$, $z(0,j)=j$, and $\rho(i,j)=r_0(i,j)$. For all $t\in\R_+$ and $i,j\in\N$,
\begin{equation}
\label{Pathw:eq:rho-r-v}
\rho_t(i,j)=(v_t(i)+r_t(i,j)+v_t(j))\I{i\neq j}.
\end{equation}
That is, the process $((r_t,v_t),t\in\R_+)$ of the decomposed genealogical distances between the individuals at fixed times coincides with the process defined from $\eta$ and $(r_0,v_0)$ in Section 6.1 of \cite{Sampl}. For $t=0$, equation~\reff{Pathw:eq:rho-r-v} holds by definition of $\rho_0$, $r_0$, and $v_0$. That equation~\reff{Pathw:eq:rho-r-v} holds for all $t\in\R_+$ can be seen from Figure~\ref{Pathw:fig:ld-space-dust}. For a formal proof, we distinguish four cases. We always assume $i\neq j$ in the following.

\emph{Case 1: $v_t(i),v_t(j)<t$, $A_0(t,i)=A_0(t,j)$.}\quad
In this case, the definition of $v_t(i)$ and $v_t(j)$ implies $A_s(t,i)\neq A_s(t,j)$ for all $s\in(t-v_t(i)\vee v_t(j),t]$. By definition of $\rho$, it follows
\begin{align*}
&\rho(z(t,i),z(t,j))
=t-v_t(i)-2\sup\{s\leq t: A_0(t,i)=A_0(t,j)\}+t-v_j(t)\\
&=\rho((t,i),(t,j))-v_t(i)-v_t(j),
\end{align*}
which is equation~\reff{Pathw:eq:rho-r-v}.

\emph{Case 2: $v_t(i),v_t(j)<t$, $A_0(t,i)\neq A_0(t,j)$.}\quad
In this case, the definition of $z(t,i)$ and $z(t,j)$ yields $A_0(t,i)=A_0(z(t,i))$ and $A_0(t,j)=A_0(z(t,j))$. In particular, it follows $A_0(z(t,i))\neq A_0(z(t,j))$. With the definition of $\rho$, it follows
\begin{align*}
&\rho(z(t,i),z(t,j))=t-v_t(i)-\rho(A_0(t,i),A_0(t,j))+t-v_j(t)\\
&=\rho((t,i),(t,j))-v_t(i)-v_t(j),
\end{align*}
which is equation~\reff{Pathw:eq:rho-r-v}.

\emph{Case 3: $v_t(i)<t,v_t(j)\geq t$.}\quad
In this case, it follows that $A_0(t,i)\neq A_0(t,j)$. From the definitions, it follows that
\[\rho(z(t,i),z(t,j))=t-v_t(i)+v_0(A_0(t,i))+r_0(A_0(t,i),A_0(t,j)).\]
Using $v_t(j)=t+v_0(A_0(t,j))$ and equation~\reff{Pathw:eq:rho-r-v} for $t=0$, we deduce that
\[v_t(i)+r_t(i,j)+v_t(j)=
t+\rho_0(A_0(t,i),A_0(t,j))+t\]
which is equation~\reff{Pathw:eq:rho-r-v}.

\emph{Case 4: $v_t(i),v_t(j)\geq t$.}\quad
Again by the definitions and by equation~\reff{Pathw:eq:rho-r-v} for $t=0$, we have
\begin{align*}
&v_t(i)+\rho(z(t,i),z(t,j))+v_t(j)\\
&=t+v_0(A_0(t,i))+r_0(A_0(t,i),A_0(t,j))+t+v_0(A_0(t,j))
=2t+\rho_0(A_0(i,j),A_0(t,j)),
\end{align*}
which is equation~\reff{Pathw:eq:rho-r-v}.
\end{rem}
\begin{rem}[Parents and starting vertices of external branches]
\label{Pathw:rem:parent-sv}
In this remark, we assume $v_0(i)=\tfrac{1}{2}\inf_{j\in\N\setminus\{i\}}\rho_0(i,j)$ for all $i\in\N$, and that $\rho_0$ is a semi-ultrametric.
Then, as in Remark 1.1 of \cite{Sampl} and the references therein, we associate with $\rho_0$ the real tree $(T_0,d_0)$ that is obtained by identifying the points with distance zero in the semi-metric space $((-\infty,0]\times\N,d_0)$, where $d_0((s,i),(t,j))=\max\{\rho_0(i,j)+s+t,|s-t|\}$.
Now we briefly sketch how the extended lookdown space $(\hat Z,\rho)$ can be isometrically embedded into a real tree $(T,d)$ that contains the genealogical trees of the individuals at all times, and we interpret the parents as starting vertices of external branches.

We define a semi-metric $d$ on $\R\times\N$ that coincides on $\R_+\times\N$ with the semi-metric $\rho$ from the beginning of Section~\ref{Pathw:sec:ld-space}, that coincides with $d_0$ on $(-\infty,0]\times\N$, and for $(s,i)\in\R_+\times\N$, $(t,j)\in(-\infty,0]\times\N$, we set $d((s,i),(t,j))=s+d_0(A_0(s,i),(t,j))$. Then we identify points with $d$-distance zero and define $(T,d)$ as the metric completion. By construction, $\rho_0(i,j)=d((0,i),(0,j))$ and $r_0(i,j)=(\rho_0(i,j)-v_0(i)-v_0(j))\I{i\neq j}=d((-v_0(i),i),(-v_0(j),j))$ for all $i,j\in\N$. Hence, $(\hat Z,\rho)$ is embedded into $(T,d)$ by the isometry that maps $(t,i)\in\R_+\times\N\subset\hat Z$ to $(t,i)\in T$, and $i\in\N\subset\hat Z$ to $(-v_0(i),i)\in T$.

For each $t\in\R_+$, Remark 5.2 in \cite{Sampl} says that $\rho_t$ is a semi-ultrametric. The associated real tree is given by the subspace $T_t=(-\infty,t]\times\N$ of $(T,d)$.
If condition \reff{eq:ext-br} is satisfied, then the quantity $v_t(i)$ and the parent $z(t,i)$ can be interpreted as the length and the starting vertex, respectively, of the external branch that ends in the leaf $(t,i)$ of $T_t$, see also \cite{Sampl}*{Remark 2.2}.
\end{rem}
\begin{rem}
\label{rem:ass-v}
In the context of Sections~\ref{Pathw:sec:ld-sampl} and \ref{Pathw:sec:proc}, assumption \reff{eq:ext-br} can be checked for $t=0$ by Proposition 3.4 in \cite{Sampl}, and for $t>0$ by Proposition 6.5 in \cite{Sampl}.
\end{rem}
\begin{figure}[ht]
\centering
\includegraphics[scale=0.32]{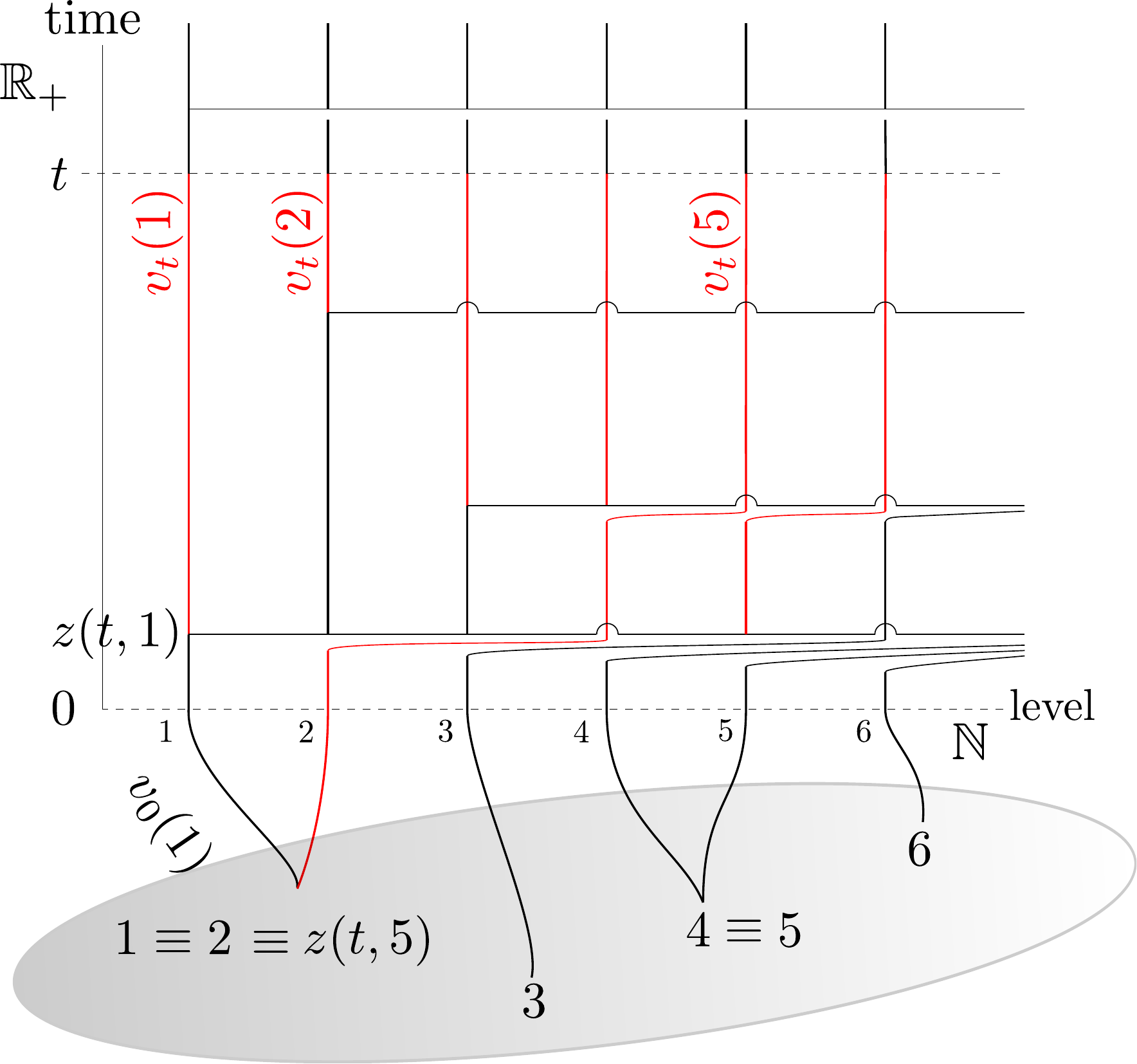}
\caption{\label{Pathw:fig:ld-space-dust}Part of an extended lookdown space. The space $\R_+\times\N$ is represented in the upper part of the figure. Time goes upwards and levels go from the left to the right. In the lower part, the metric space obtained from $\N$, endowed with the semi-metric obtained from $r_0$ is symbolized. For each $i\in\N$, the junction between the individual $(0,i)$ and its ancestor $i$ has length $v_0(i)$. Individuals that are in the same block in a reproduction event have genealogical distance zero and are identified. In the figure, they are connected by horizontal lines. In this example, there are no simultaneous multiple reproduction events. The genealogical distances between the individuals $(t,i)$ and their respective parents $z(t,i)$ equal $v_t(i)$, they are represented by red lines. The genealogical distance between any two individuals is the sum of the lengths of the vertical parts of the path from one individual to the other, plus the distance in the metric space obtained from $r_0$ if this space has to be traversed.}
\end{figure}

\subsection{Extinction of parts of the population}
\label{Pathw:sec:ld-ext}
In the beginning of Section \ref{Pathw:sec:ld-space},
we stated the lookdown model as an infinite particle system. The particles undergo reproduction events that are given by the point measure $\eta$.
In reproduction events encoded by points $(t,\pi)$ of $\eta$ where the partition $\pi$ has finitely many blocks, the particles at time $t$ have only finitely many ancestors among the infinite population at time $t-$, hence particles die at time $t$. A particle can also die due to an accumulation of reproduction events in which its level is pushed to infinity.

For $(s,i)\in\R_+\times\N$ and $t\in[s,\infty)$, let $D_t(s,i)$\label{Pathw:not:Dtsi} be the lowest level that is occupied at time $t$ by a descendant of the individual $(s,i)$, that is,
\[D_t(s,i)=\inf\{j\in\N:A_s(t,j)=i\}\]
with $D_t(s,i)=\infty$ if and only if there exists no $j\in\N$ with $i=A_s(t,j)$. This quantity corresponds to the forward level process in \cite{PW06} and to the fixation line in \cite{H15}. The map $t\mapsto D_t(s,i)$ is non-decreasing. Let $\tau_{s,i}$ be the extinction time of the part of the population
that descends from the individual $(s,i)$, that is,
\[\tau_{s,i}=\inf\{t\in[s,\infty): D_t(s,i)=\infty\}.\label{Pathw:not:tausi}\]
Then the set of times at which such parts of the population become extinct is given by
\begin{equation}
\label{Pathw:eq:Thetaext}
\Thetaext:=\{\tau_{s,i}:s\in\R_+,i\in\N\}.
\end{equation}
\begin{rem}
\label{Pathw:rem:die}
In a reproduction event that is encoded by a point $(t,\pi)$ of $\eta$ with $\#\pi=\infty$, every individual that sits on a level $i\in\N$ at time $t-$ has a descendant
at time $t$.
Hence the two mechanisms mentioned in the beginning of this subsection are the only possibilities for a particle to die.
\end{rem}
\begin{figure}
\centering
\includegraphics[scale=0.32]{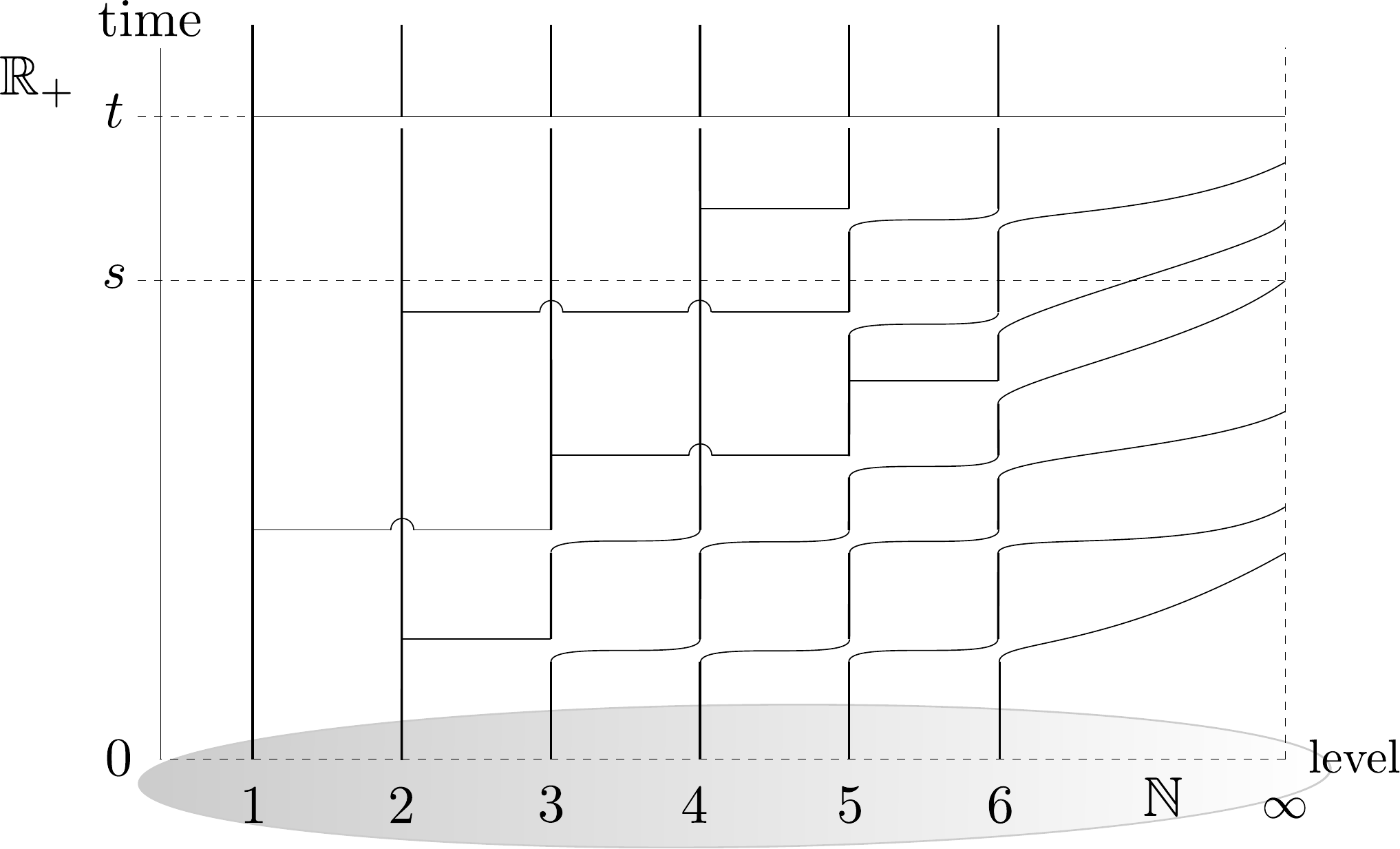}
\caption{\label{Pathw:fig:ld-space-nd}Part of a lookdown space.
Only the reproduction events with offspring on the first $6$ levels are drawn. At time $t$, a reproduction event occurs in which the whole population is replaced by the offspring of the individual on level $1$. The limits $(t-,i)=\lim_{r\uparrow t}(r,i)$, $i\geq 2$ are part of the boundary of the lookdown space. Due to accumulations of jumps, the lines $t'\mapsto D_{t'}(s',i)$ may hit infinity, here symbolized by a dashed line, similarly to illustrations of the lookdown graph in \cite{PW06}. This occurs for instance at time $s$, the limit $\lim_{r\uparrow s}D_r(0,3)$ is part of the boundary. Further elements of the boundary are obtained from Cauchy sequences at fixed times.}
\end{figure}

\subsection{The \texorpdfstring{$\Xi$}{Xi}-lookdown model}
\label{Pathw:sec:ld-Xi}
We recall here the simple point measure $\eta$ that is used to drive the lookdown  model in \cite{Sampl}, cf.\ also the references therein.

Let $\Delta$ be the simplex
\[\Delta=\{x=(x_1,x_2,\ldots): x_1\geq x_2\geq \ldots\geq 0,\left|x\right|_1\leq 1\},\]
where $|x|_p=(\sum_{i\in\N} x_i^p)^{1/p}$.
Let $\kappa$\label{Pathw:not:kappa} be the probability kernel from $\Delta$ to $\p$ associated with Kingman's correspondence, that is, $\kappa(x,\cdot)$ is the distribution of the paintbox partition associated with $x\in\Delta$, see e.\,g.\ Section 2.3.2 in \cite{Bertoin}.

Let $\Xi$ be a finite measure on $\Delta$. We decompose
\[\label{Pathw:not:Xi0}\Xi=\Xi_0+\Xi\{0\}\delta_0.\]
For distinct integers $i,j\in\N$, we denote by $K_{i,j}$\label{Pathw:not:Kij} the partition of $\N$ that contains the block $\{i,j\}$ and apart from that only singleton blocks. We define a $\sigma$-finite measure $H_\Xi$\label{Pathw:not:HXi} on $\p$ by
\begin{equation*}
H_\Xi(d\pi)=\int_\Delta\kappa(x,d\pi)\left|x\right|_2^{-2} \Xi_0(dx) + \Xi\{0\}\;\sum_{1\leq i<j}\delta_{K_{i,j}}(d\pi).
\end{equation*}
In the following sections, $\eta$\label{Pathw:not:eta-Poisson} is always a Poisson random measure on $(0,\infty)$ with intensity $dt\;H_\Xi(d\pi)$.
Then $\eta$ satisfies a.\,s.\ condition~\reff{Pathw:eq:ass-eta}, as checked in equation (5.5) of \cite{Sampl}. The lookdown model can therefore be driven by $\eta$.

The measure $\Xi$ is called dust-free
if and only if
\begin{equation*}
\Xi\{0\}>0\quad\text{or}\quad\int\left|x\right|_1\left|x\right|_2^{-2}\Xi_0(dx)=\infty.
\end{equation*}
In this case, each particle reproduces with infinite rate, hence $v_t(i)=0$ for all $t\in(0,\infty)$ and $i\in\N$ a.\,s.
We denote by $\nd$\label{Pathw:not:nd} the subset of finite measures on $\Delta$ that are dust-free, and by $\dust$ its complement in the set of finite measures on $\Delta$.

The extended lookdown space in Figure~\ref{Pathw:fig:ld-space-dust} could be the extended lookdown space associated with a typical realization of an appropriate Poisson random measure $\eta$ with dust, and a decomposed semimetric $(r_0,v_0)$. Figure~\ref{Pathw:fig:ld-space-nd} illustrates the dust-free case.

We speak of a large reproduction event when a particle has offspring on a positive proportion of the levels. We denote by $\Theta_0$ the set of times at which large reproduction events occur:
\begin{equation}
\label{Pathw:eq:Theta0}
\Theta_0=\{t\in(0,\infty):\text{there exist $\pi\in\p$ and $B\in\pi$ with $\eta\{(t,\pi)\}>0$ and $|B|>0$}\}.
\end{equation}
Here $|B|=n^{-1}\lim_{n\to\infty}B\cap\sn$ denotes the asymptotic frequency of the block $B$ if it exists. (Note that a.\,s., all the asymptotic frequencies in the definition of $\Theta_0$ exist by Kingman's correspondence and as $\eta$ has countably many points.)
The measure $\Xi_0$ governs the large reproduction events, as opposed to $\Xi\{0\}$, which gives the rates of the binary reproduction events (i.\,e.\ the reproduction events in which the reproducing particle has only one offspring).
A.\,s.\ by definition of $H_\Xi$ and Kingman's correspondence, each reproduction event is either binary or large.
In case $\Xi\in\dust$, the set $\Theta_0$ equals a.\,s.\ the set $\{t\in(0,\infty):\eta(\{t\}\times\p)>0\}$ of reproduction times, as there are a.\,s.\ no binary reproduction events.

\section{Sampling measures and jump times}
\label{Pathw:sec:ld-sampl}
In this section, we consider mathematical objects that are defined from a realization of the Poisson random measure $\eta$ and a random (decomposed) distance matrix on an event of probability $1$. Stochastic processes are read off from these constructions in Section~\ref{Pathw:sec:proc}. We defer the proofs of the main statements in Section~\ref{Pathw:sec:ld-sampl} to the second part of the article which begins in Section~\ref{Pathw:sec:outline-proofs}.

\subsection{The case without dust}
\label{Pathw:sec:ld-sampl-nd}
We construct a family of probability measures on the lookdown space. We consider regularity of this family in the weak topology. In the case of coming down from infinity, we also consider regularity of a family of subsets of the lookdown space with respect to the Hausdorff distance.

Let the Poisson random measure $\eta$ be defined from the finite measure $\Xi$ on $\Delta$ as in Section~\ref{Pathw:sec:ld-Xi}. Let $(X,r,\mu)$ be a metric measure space and $\rho_0$ be an independent $\R^{\N^2}$-valued
random variable that has the distance matrix distribution of $(X,r,\mu)$.
That is, we can and will assume $\rho_0=(r(x(i),x(j)))_{i,j\in\N}$ for a $\mu$-iid sequence $x(1),x(2),\ldots$ in $(X,r)$.
We view $\rho_0$ as a random semi-metric on $\N$.
Let $(Z,\rho)$ be the lookdown space associated with $\eta$ and $\rho_0$ as defined in Section~\ref{Pathw:sec:ld-space}.

For each $t\in\R_+$ and $n\in\N$, let the probability measure $\mu^n_t$\label{Pathw:not:munt} on $(Z,\rho)$ be the uniform measure on the first $n$ individuals at time $t$, that is,
\[\mu^n_t=\frac{1}{n}\sum_{i=1}^n\delta_{(t,i)}.\]
Let $\mu_t$\label{Pathw:not:mut} be the weak limit of $\mu^n_t$, provided it exists.
For almost all realizations of $\eta$ and $\rho_0$, these weak limits exist simultaneously for all $t\in\R_+$ by Theorem~\ref{Pathw:thm:ld-nd} below. The convergence is uniform for the Prohorov metric $\dP^Z$ on $(Z,\rho)$ for $t$ in compact intervals.
\begin{thm}
\label{Pathw:thm:ld-nd}
Assume $\Xi\in\nd$. Then there exists an event of probability $1$ on which the following assertions hold:
\begin{enumerate}[label=(\roman{*}),ref=(\roman{*})]
\item\label{Pathw:eq:thm:ld-nd:conv} For all $T\in\R_+$,
\begin{equation*}
\lim_{n\to\infty}\sup_{t\in[0,T]}\dP^Z(\mu^n_t,\mu_t)=0.
\end{equation*}
\item\label{Pathw:item:thm:ld-nd:reg} The map $t\mapsto \mu_t$ is càdlàg in the weak topology on the space of probability measures on $(Z,\rho)$. The set $\Theta_0$, defined in~\reff{Pathw:eq:Theta0}, is the set of jump times.
\item\label{Pathw:item:thm:ld-nd:atoms-jumps} For all $t\in\Theta_0$, the measure $\mu_t$ contains atoms, and the left limit $\mu_{t-}$ is non-atomic. 
\end{enumerate}
\end{thm}
The proof is given in Section \ref{sec:constr-ld-nd}. Note that the measure $\mu^n_t$ assigns mass $1/n$ to each offspring in a reproduction event. At the time $t$ of a large reproduction event, the measure $\mu_t$ has an atom on each family of individuals that descend from the same individual at time $t-$.
Indeed, these individuals have genealogical distance zero and are identified in the lookdown space. In the Kingman case, that is, if $\Xi=\delta_0$, there are a.\,s.\ no large reproduction events and $t\mapsto \mu_t$ is a.\,s.\ continuous.

For $t\in\R_+$, let $X_t$\label{Pathw:not:Xt} be the closure of the set of individuals $\{t\}\times\N$ at time $t$, defined as a subspace of the complete space $(Z,\rho)$. Clearly, the probability measures $\mu^n_0$ do not depend on $\eta$. Their weak limit exists a.\,s.\ by the following lemma which is essentially Vershik's proof \cite{Vershik02}*{Theorem 4} of the Gromov reconstruction theorem, see also \cite{Sampl}*{Proposition 10.5}.
We write $\supp\mu'$ for the closed support of a measure $\mu'$.
\begin{lem}
\label{Pathw:lem:mu0}
The weak limit $\mu_0$ of the probability measures $\mu_0^n$ on $Z$ exists a.\,s. The metric measure spaces $(\supp\mu,r,\mu)$ and $(X_0,\rho,\mu_0)$ are a.\,s.\ measure-preserving isometric.
\end{lem}
\begin{proof}
For our $\mu$-iid sequence $(x(i),i\in\N)$, the empirical measures
$\mu^n:=n^{-1}\sum_{i=1}^n\delta_{x(i)}$
converge to $\mu$ a.\,s.\ by the Glivenko-Cantelli theorem. A.\,s., the isometry $\{x(i):i\in\N\}\to Z$, $x(i)\mapsto (0,i)$ can be extended to a surjective isometry $\varphi:(\supp\mu,r)\to(X_0,\rho)$ with $\varphi(\mu^n)=\mu_0^n$, hence also $\varphi(\mu)=\mu_0$.
\end{proof}

In Theorem~\ref{Pathw:thm:ld-CDI} below, we consider measures $\Xi$ that satisfy the ``coming down from infinity''-assumption that there exists an event of probability $1$ on which the number $\#\{A_s(t,j):j\in\N\}$ of ancestors at time $s$ of the individuals at time $t$ is finite for all $t\in(0,\infty)$ and $s\in(0,t)$.
Let $\MCDI$\label{Pathw:not:MCDI} denote the subset of those finite measures $\Xi$ on $\Delta$ that satisfy this assumption.
Then $\MCDI\subset\nd$. Indeed, if $\Xi\in\dust$, then the rate at which a given ancestral lineage merges with any other ancestral lineage is finite (as discussed e.\,g.\ in \cite{Sampl}*{Section 6.2}) which implies $\#\{A_s(t,j):j\in\N\}=\infty$ a.\,s.\ for all $0\leq s<t$. Furthermore, if $\Xi\in\MCDI$, then the set $\Thetaext$ of extinction times, defined in~\reff{Pathw:eq:Thetaext}, is a.\,s.\ dense in $(0,\infty)$. Indeed, if $\Thetaext$ has no points in an interval $[s,t]$, then the individuals on all levels at time $s$ are ancestors of individuals at time $t$.
\begin{rem}
``Coming down from infinity'' usually refers to a property of a $\Xi$-coalescent. A $\Xi$-coalescent describes the genealogy at a fixed time in our lookdown model. The assumption on $\Xi$ that for each fixed $t\in(0,\infty)$, there exists an event of probability $1$ on which $\#\{A_s(t,j):j\in\N\}<\infty$ holds for all $s\in(0,t)$ already implies $\Xi\in\MCDI$. Indeed, from this a priori weaker assumption, it follows that there exists an event of probability $1$ on which $\#\{A_s(t,j):j\in\N\}<\infty$ for all $t\in(0,\infty)\cap\Q$ and $s\in(0,t)$. For a general $t\in(0,\infty)$, there exists a rational $t'\in(s,t)$, and monotonicity of the number of ancestors yields
\[\#\{A_s(t,j):j\in\N\}=\#\{A_s(t',A_{t'}(t,j)):j\in\N\}\leq\#\{A_s(t',j):j\in\N\}<\infty.\]
\end{rem}
\begin{rem}
\label{Pathw:rem:compact-CDI}
If $\Xi\in\MCDI$, then there exists an event of probability one on which all subsets $X_t\subset Z$ with $t\in(0,\infty)$ are compact.
Indeed, $\#\{A_s(t,j):j\in\N\}<\infty$ for all $s\in(0,t)$ implies that the complete subspace $X_t$ is totally bounded by definition of the metric $\rho$ on the lookdown space $Z$.

If $(X,r)$ is compact, then $X_0$ is a.\,s.\ compact by Lemma~\ref{Pathw:lem:mu0}. This also implies assertion~\ref{Pathw:item:thm:ld-CDI-supp} of Theorem~\ref{Pathw:thm:ld-CDI} below for $t=0$.
\end{rem}

Recall that the Hausdorff distance between two subsets $A,B$ of a metric space $(Y,d)$ is defined as the infimum over those $\ep>0$ such that $d(a,B)<\ep$ for all $a\in A$ and $d(A,b)<\ep$ for all $b\in B$. The Hausdorff distance is a metric on the set of closed subspaces of $(Y,d)$, see e.\,g.\ \cite{BBI}.
\begin{thm}
\label{Pathw:thm:ld-CDI}
Assume $\Xi\in\MCDI$ and that $(X,r)$ is compact. Then the following assertions hold on an event of probability $1$:
\begin{enumerate}[label=(\roman{*}),ref=(\roman{*})]
\item\label{Pathw:item:thm:ld-CDI-supp} For each $t\in\R_+$, the compact set $X_t$ is the closed support of $\mu_t$.
\item\label{Pathw:item:thm:ld-CDI-reg} The map $t\mapsto X_t$ is càdlàg for the Hausdorff distance on the set of closed subsets of $(Z,\rho)$. The set $\Thetaext$
is the set of jump times. For each $t\in \Thetaext$, the set $X_t$ and the left limit $X_{t-}$ are not isometric.
\end{enumerate}
\end{thm}
The proof is given in Section \ref{sec:constr-ld-nd}.
\begin{rem}
By Theorems \ref{Pathw:thm:ld-nd} and \ref{Pathw:thm:ld-CDI}, there exists an event of probability $1$ on which $\mu_t=\mu_{t-}$ and $X_t\subsetneq X_{t-}$ for all $t\in\Thetaext\setminus\Theta_0$. In particular, the closed support of $\mu_{t-}$ is strictly smaller than $X_{t-}$ for these $t$ a.\,s. The set $X_{t-}\setminus\supp \mu_{t-}\subset X_{t-}$ is equal to $X_{t-}\setminus X_t$ a.\,s., this is the part of the population at time $t-$ that dies out at time $t$.
\end{rem}

We conclude this subsection with a side observation (Proposition \ref{Pathw:prop:int-jumps}) on the intersection of the sets $\Theta_0$ and $\Thetaext$ of jump times in Theorems \ref{Pathw:thm:ld-nd} and \ref{Pathw:thm:ld-CDI}.
Reproduction events in which the whole population is replaced by finitely many particles and their offspring occur at the times in the set
\[\Theta_{\rm f}:=\{t\in(0,\infty):
\text{ there exists $\pi\in\p$ with $\eta\{(t,\pi)\}>0$ and $\#\pi<\infty$}\}.\]
By the construction in Section~\ref{Pathw:sec:ld-space}, all particles with level larger than $\#\pi$ die in a reproduction event that is encoded by a point $(t,\pi)$ of $\eta$.
\begin{rem}
If $\Xi$ is concentrated on $\{(x,0,0,\ldots):x\in[0,1]\}\subset\Delta$, then a.\,s., no simultaneous multiple reproduction events occur. This case corresponds to the coalescents with multiple collisions ($\Lambda$-coalescents).
In this case, $\Theta_{\rm f}$ is a.\,s.\ the set of times at which the whole population is replaced by a single particle and its offspring.
If $\Xi$ is concentrated on $\{(x,0,0,\ldots):x\in[0,1)\}\subset\Delta$, then $\Theta_{\rm f}=\emptyset$ a.\,s.
More generally, $\Theta_{\rm f}=\emptyset$ a.\,s.\ if and only if $\Xi\{x\in\Delta:x_1+\ldots+x_k=1\text{ for some }k\in\N\}=0$.
\end{rem}
\begin{prop}
\label{Pathw:prop:int-jumps}
A.\,s., $\Thetaext\cap\Theta_0=\Theta_{\rm f}$.
\end{prop}
For the proof below in this subsection, as well as for later use in Sections~\ref{Pathw:sec:part} and \ref{Pathw:sec:unif-ld-proofs}, we now express $\eta$ in terms of a collection of Poisson processes.
Recall the set $\p_n$ of partitions of $\sn=\{1,\ldots,n\}$, and let $\pi'_1, \pi'_2,\ldots$ be an arbitrary enumeration of the set
\[\bigcup_{n\in\N}(\p_n\setminus\big\{\{\{1\},\ldots,\{n\}\}\big\})\]
of finite partitions that consist not only of singleton blocks.
For $k\in\N$, let
\[\p'_k=\{\pi\in\p:\gamma_n(\pi)=\pi'_k \text{ with $n$ such that }\pi'_k \text{ is a partition of }\sn\},\]
and
\[J_{t,k}=\eta((0,t]\times\p'_k)\]
for $t\in\R_+$.
Then the processes $(J_{t,k},t\in\R_+)$ form a collection (indexed by $k\in\N$) of Poisson processes. We endow $\R_+^{\N}$ with the product topology and consider the $\R_+^\N$-valued stochastic process
\begin{equation}
\label{Pathw:eq:def-J}
J=(J_t,t\in\R_+)=((J_{t,k},k\in\N),t\in\R_+)
\end{equation}
Note that $J$ has independent and stationary increments and a.\,s.\ càdlàg paths. Hence, $J$ is a strong Markov process and Feller continuous (i.\,e.\ the elements of its semigroup preserve the set of bounded continuous functions), so that we obtain from e.\,g.\ Theorem (5.1) in Chapter I of Blumenthal \cite{Blu92} that $J$ is quasi-left-continuous.
Let $\F=(\F_t,t\in\R_+)$ be the complete filtration induced by $J$. Then $\F_t$ is the sigma field generated by the random measure $\eta(\cdot\times((0,t]\times\p))$ and all null events. By condition~\reff{Pathw:eq:ass-eta}, $J$ stays finite a.\,s.\ and the set of jump times of $J$ equals a.\,s.\ the set of reproduction times $\{t\in(0,\infty):\eta(\{t\}\times\p)>0\}$.

\begin{proof}[Proof of Proposition \ref{Pathw:prop:int-jumps}]
We use notation also from Section~\ref{Pathw:sec:ld-ext}. First we show $\Theta_{\rm f}\subset\Thetaext\cap\Theta_0$ a.\,s. Let $t\in\Theta_{\rm f}$. Then there exists a point $(t,\pi)$ of $\eta$ with $i-1:=\#\pi<\infty$. By condition~\reff{Pathw:eq:ass-eta}, as the level of each particle is non-decreasing in time, and as offspring always has a higher level than the reproducing particle, there exists a.\,s.\ a time $s\in(0,t)$ such that $D_{s'}(s,i)=i$ for all $s'\in[s,t)$. In particular, all descendants of $(s,i)$ at time $t-$ have level at least $i$. Hence, $(s,i)$ has no descendants at time $t$, that is, $D_t(s,i)=\infty$, and it follows $t\in\Thetaext$.
Clearly, $\#\pi<\infty$ implies that $\pi$ contains blocks of infinite size. The definition of $H_\Xi$ and Kingman's correspondence imply $\Theta_{\rm f}\subset\Theta_0$ a.\,s. It remains to show that $\Thetaext\cap\Theta_0\subset\Theta_{\rm f}$ a.\,s.

On the event of probability $1$ on which condition~\reff{Pathw:eq:ass-eta} holds, particles on any level remain on that level for a positive amount of time. This implies
\[\Thetaext=\{\tau_{s,i}:s\in\Q_+,i\in\N\}\quad\text{a.\,s.}\]

For $s\in\R_+$, $i,n\in\N$, we define the $\F$-stopping time
\[\tau_{s,i,n}=\inf\{t\geq s: D_t(s,i)\geq n\}.\]
Then $\tau_{s,i,n}$ is non-decreasing in $n$, and
$\tau_{s,i}=\lim_{n\to\infty}\tau_{s,i,n}$.
We assume w.\,l.\,o.\,g.\ $\Xi(\Delta)>0$. Then $\tau_{s,i,n}\in[s,\infty)$ for all $s\in\R_+$ and $i,n\in\N$ a.\,s.
Let $E_{s,i}$ be the event that $\tau_{s,i,n}<\tau_{s,i}$ for all $n\in\N$.

A.\,s.\ by \reff{Pathw:eq:ass-eta}, on the event $E_{s,i}^c$ that $\tau_{s,i}=\tau_{s,i,n}$ for some $n\in\N$, a particle on a level below $n$ at time $\tau_{s,i,n}-$ dies at time $\tau_{s,i,n}$ due to a reproduction event that lies in $\Theta_{\rm f}$.

To show the assertion of the proposition, it now suffices to show that $\tau_{s,i}\notin\Theta_0$ a.\,s.\ on $E_{s,i}$.
We define the $\F$-stopping time $\tilde\tau_{s,i}$
by $\tilde\tau_{s,i}=\tau_{s,i}\Ind_{E_{s,i}}+\infty\Ind_{E_{s,i}^c}$. Then the $\F$-stopping times
\[\tilde\tau_{s,i,n}:=\left\{\begin{aligned}
\tau_{s,i,n}\quad&\text{if }\tau_{s,i,n}<\tau_{s,i}\\
\tau_{s,i,n}\vee n\quad&\text{if }\tau_{s,i,n}=\tau_{s,i}
\end{aligned}\right.\]
form an announcing sequence for $\tilde\tau_{s,i}$, that is, $\tilde\tau_{s,i,n}<\tilde\tau_{s,i}$ a.\,s.\ and $\tilde\tau_{s,i}=\lim_{n\to\infty}\tilde\tau_{s,i,n}$ a.\,s. Quasi-left-continuity of $J$ implies $J_{\tilde\tau_{s,i}-}=J_{\tilde\tau_{s,i}}$ a.\,s. Hence, a.\,s.\ on $E_{s,i}$, no reproduction event occurs at time $\tilde\tau_{s,i}=\tau_{s,i}$, and we have $\tau_{s,i}\notin\Theta_0$.
\end{proof}
\begin{rem}
In the Kingman case, the set $\Thetaext$ of extinction times is described by Poisson processes by Dahmer, Knobloch, and Wakolbinger \cite{DKW14}*{Proposition 1}, see also the references therein.
\end{rem}

\subsection{The general case}
\label{Pathw:sec:ld-sampl-dust}
We construct a family of probability measures on the Cartesian product of the extended lookdown space and the mark
space $\R_+$. Let $(r_0,v_0)$ be an independent $\R^{\N^2}\times\R^\N$-valued random variable that has the marked distance matrix distribution of an $(\R_+)$-marked metric measure space $(X,r,m)$.
That is, $(X,r)$ is a complete and separable metric space, $m$ is a probability measure on the Borel sigma algebra on the product space $X\times\R_+$, and we may assume that $(x(i),v(i))_{i\in\N}$ is an $m$-iid sequence in $X\times\R_+$ and set $r_0(i,j)=r(x(i),x(j))$ for $i,j\in\N$. Then we can view $r_0$ as a random semi-metric on $\N$.
Let $(\hat Z,\rho)$ be the extended lookdown space associated with $\eta$ and $(r_0,v_0)$, as defined in Section \ref{Pathw:sec:ld-dec}.
We endow $\hat Z\times\R_+$ with the product metric
$d^{\hat Z\times\R_+}((z,v),(z',v')) =\rho(z,z')\vee\left|v-v'\right|$. We denote the Prohorov metric on the space of probability measures on $\hat Z\times\R_+$ by $\dP^{\hat Z\times \R_+}$.
Recall from Section~\ref{Pathw:sec:ld-dec} also the parent $z(t,i)$ and the genealogical distance $v_t(i)$ between the individual $(t,i)$ and its parent. For each $t\in\R_+$ and $n\in\N$, we define a probability measure $m^n_t$ on $\hat Z\times\R_+$ by
\begin{equation}
\label{eq:def-mt}
m^n_t=\frac{1}{n}\sum_{i=1}^n\delta_{(z(t,i),v_t(i))}
\end{equation}
The first component $m^n_t(\cdot\times\R_+)$ lays mass on the parents of the first $n$ individuals at time $t$.
The second component $m^n_t(\hat Z\times\cdot)$ records the genealogical distances to these parents. Let $m_t$\label{Pathw:not:mt} denote the weak limit of $m^n_t$ provided it exists. This existence is addressed in Theorem~\ref{Pathw:thm:ld-dust} below in the case with dust, and in Lemma~\ref{Pathw:lem:m0}, Remark~\ref{Pathw:rem:Z-consistence-random}, and Corollary~\ref{Pathw:cor:nd-mdm} below in the case without dust.
\begin{thm}
\label{Pathw:thm:ld-dust}
Assume $\Xi\in\dust$. Then the following assertions hold on an event of probability $1$:
\begin{enumerate}[label=(\roman{*}),ref=(\roman{*})]
\item\label{Pathw:eq:thm:ld-dust:conv} For all $T\in\R_+$,
\begin{equation*}
\lim_{n\to\infty}\sup_{t\in[0,T]}\dP^{\hat Z\times\R_+}(m^n_t,m_t)=0.
\end{equation*}
\item\label{Pathw:item:thm:ld-dust:reg} The map $t\mapsto m_t$ is càdlàg in the weak topology on the space of probability measures on $\hat Z\times\R_+$. The set $\Theta_0$, defined in~\reff{Pathw:eq:Theta0}, is the set of jump times.
\item\label{Pathw:item:thm:ld-dust:zero} For each $t\in(0,\infty)$, the left limit $m_{t-}$ satisfies $m_{t-}(\hat Z\times\{0\})=0$. For each $t\in\Theta_0$, it holds $m_t(\hat Z\times\{0\})>0$.
\item\label{Pathw:item:thm:ld-dust:atomic} If $m$ is purely atomic, then $m_t$ and $m_{t-}$ are purely atomic for all $t\in(0,\infty)$.
\end{enumerate}
\end{thm}
The proof is given in Section \ref{sec:constr-ld-dust}.
In Proposition~\ref{Pathw:prop:repr-mt}, the measures $m_t$ are stated explicitly. At the times of large reproduction events, $v_t(i)=0$ for all individuals $i$ with levels in a non-singleton block. This yields the positive mass of $m_t(\hat Z\times\cdot)$ in zero asserted in Theorem~\ref{Pathw:thm:ld-dust}\ref{Pathw:item:thm:ld-dust:zero}.

\begin{lem}
\label{Pathw:lem:m0}
The weak limit $m_0$ of the probability measures $m_0^n$ on $\hat Z\times\R_+$ exists a.\,s.
\end{lem}
\begin{proof}
This follows analogously to Lemma~\ref{Pathw:lem:mu0}.
\end{proof}

\begin{rem}
\label{Pathw:rem:v-nd}
We recall that in case $\Xi\in\nd$, the rate at which each particle reproduces is infinite. In this case, there exists an event of probability $1$ on which $v_t(i)=0$ and $(t,i)=z(t,i)$ for all $t\in(0,\infty)$ and $i\in\N$, cf.\ \cite{Sampl}*{Section 6.2}.
\end{rem}
As in Subsection~\ref{Pathw:sec:ld-sampl-nd}, we define the probability measures
\begin{equation}
\label{eq:def-mut-ext}
\mu^n_t=\frac{1}{n}\sum_{i=1}^n\delta_{(t,i)}
\end{equation}
on $\hat Z$, and we denote their weak limits by $\mu_t$.
\begin{rem}
\label{Pathw:rem:Z-consistence-random}
If $m(X\times\{0\})=1$, then $m=\mu\otimes\delta_0$ for a probability measure $\mu$ on $X$, hence $v_0=0$ a.\,s.\ and $r_0$ has the distance matrix distribution of the metric measure space $(X,r,\mu)$.
By Remark~\ref{Pathw:rem:Z-consistence}, the extended lookdown space $(\hat Z,\rho)$ coincides in this case a.\,s.\ with the lookdown space associated with $\eta$ and $r_0$. With this identification, the assertions of Theorem~\ref{Pathw:thm:ld-nd} also hold in the context of the present subsection if $\Xi\in\nd$ and $m(\hat Z\times\{0\})=1$.
\end{rem}

Moreover, in the case without dust, the following corollary to Theorem~\ref{Pathw:thm:ld-nd} also holds for the extended lookdown space and the more general initial configuration $(r_0,v_0)$ in the present subsection.
\begin{cor}
\label{Pathw:cor:nd-mdm}
Assume $\Xi\in\nd$. Then a.\,s., the probability measures $\mu_t$ exist for all $t\in(0,\infty)$. The map $t\mapsto\mu_t$ is a.\,s.\ càdlàg on $(0,\infty)$ in the weak topology on the space of probability measures on $(\hat Z,\rho)$, and $\Theta_0$ is a.\,s.\ the set of jump times. Moreover, the family of probability measures $(\mu_t,t\in(0,\infty))$ satisfies a.\,s.\ assertion~\ref{Pathw:item:thm:ld-nd:atoms-jumps} of Theorem~\ref{Pathw:thm:ld-nd}. A.\,s., also the probability measures $m_t$ exist for all $t\in(0,\infty)$ and satisfy $m_t=\mu_t\otimes\delta_0$.
\end{cor}
\begin{proof}
Let $(Z',\rho')$ be the lookdown space associated with $\eta$ and $(0)_{i,j\in\N}$. By Theorem~\ref{Pathw:thm:ld-nd}, the probability measures $\mu'_t$ defined on $Z'$ analogously to $\mu_t$ satisfy the assertion.

Let $\ep>0$, let $Z_\ep$ be the closure of $[\ep,\infty)\times\N$ in $(\hat Z,\rho)$, and $Z'_\ep$ the closure of $[\ep,\infty)\times\N$ in $(Z',\rho')$.
The construction in the beginning of Section~\ref{Pathw:sec:ld-space} yields
$\rho'((t,i),(u,j))\wedge\ep=\rho((t,i),(u,j))\wedge\ep$
for all $(t,i)$, $(u,j)\in[\ep,\infty)\times\N$. Hence, the map from $[\ep,\infty)\times\N\subset Z'_\ep$ to $Z_\ep$, given by $(t,i)\mapsto(t,i)$, can be extended to a homeomorphism $h:Z'_\ep\to Z_\ep$.
Hence a.\,s., the weak limits
\[\mu_t=\text{w-}\lim_{n\to\infty}\mu^n_t=\text{w-}\lim_{n\to\infty}h(\mu'^n_t)=h(\mu'_t)\]
exist for all $t\in[\ep,\infty)$ and the assertion on $(\mu_t,t\in(0,\infty))$ follows.
The assertion on $(m_t,t\in(0,\infty))$ now follows from Remark~\ref{Pathw:rem:v-nd}, the definitions \reff{eq:def-mt} and \reff{eq:def-mut-ext} of $m^n_t$ and $\mu^n_t$, and the definitions of $m_t$ and $\mu_t$ as weak limits.
\end{proof}

\section{Stochastic processes}
\label{Pathw:sec:proc}
\subsection{The case without dust}
\label{Pathw:sec:proc-nd}
From the construction on the lookdown space in Section~\ref{Pathw:sec:ld-sampl-nd}, we now read off stochastic processes with values in the space $\M$ of isomorphy classes of metric measure spaces and in the space $\bM$ of strong isomorphy classes of compact metric measure spaces. First we recall these state spaces from the literature \cites{GPW09,Miermont09,EW06}.

As stated in the introduction, we call two metric measure spaces $(X',r',\mu')$, $(X'',r'',\mu'')$ isomorphic if there exists an isometry $\varphi$ from the closed support $\supp\mu'\subset X'$ to $\supp\mu''\subset X''$ with $\mu''=\varphi(\mu')$. We denote the isomorphy class by $\ew{X',r',\mu'}$.
We endow the space $\M$ of isomorphy classes of metric measure spaces with the Gromov-Prohorov metric $\dGP$ which is complete and separable and induces the Gromov-weak topology, as shown in \cite{GPW09}.

Moreover, we call two metric measure spaces $(X',r',\mu')$, $(X'',r'',\mu'')$ strongly isomorphic if they are measure-preserving isometric, that is, if there exists a surjective isometry $\varphi:X\to X'$ with $\mu''=\varphi(\mu')$. We denote the strong isomorphy class by $\es{X',r',\mu'}$. We endow the space $\bM$ of strong isomorphy classes of compact metric measure spaces with the Gromov-Hausdorff-Prohorov metric $\dGHP$, given by
\[\dGHP((X',r',\mu'),(X'',r'',\mu''))=\inf_{Y,\varphi',\varphi''}\{\dP^Y(\varphi'(\mu'),\varphi''(\mu''))\vee\dH^Y(\varphi'(X'),\varphi''(X''))\}\]
where the infimum is over all isometric embeddings $\varphi':X'\to Y$, $\varphi'':X''\to Y$ into complete and separable metric spaces $Y$. Here we denote by $\dP^Y$ and by $\dH^Y$ the Prohorov and the Hausdorff distance, respectively, over a metric space $Y$.
Then $(\bM,\dGHP)$ is a complete and separable metric space, see \cites{Miermont09,EW06},
and $\dGHP$ induces the Gromov-Hausdorff-Prohorov topology on $\bM$.
The Hausdorff distance in the definition of $\dGHP$ compares the metric spaces also where the probability measures charges them with negligible mass.

We work with the lookdown space $(Z,\rho)$, the families of sampling measures $(\mu_t,t\in\R_+)$, and the subspaces $X_t$ of the lookdown space from Section~\ref{Pathw:sec:ld-sampl-nd}.
Recall that the randomness comes from a Poisson random measure $\eta$ that is characterized by $\Xi\in\nd$,
and from an independent random variable $\rho_0$ with the distance matrix distribution of a metric measure space $(X,r,\mu)$.

We say that a Markov process $(Y_t,t\in\R_+)$ with values in a (not necessarily locally compact)  metric space is Feller continuous if for each $t\in\R_+$, the law of $Y_t$ depends continuously with respect to the weak topology on the initial state.
\begin{prop}
\label{Pathw:thm:proc-nd-MvFV}
Assume $\Xi\in\nd$. Then a Feller-continuous strong Markov process with values in $\M$ is given a.\,s.\ by $(\ew{Z,\rho,\mu_t},t\in\R_+)$. 
\end{prop}
\begin{prop}
\label{Pathw:thm:proc-nd-Mvec}
Assume $\Xi\in\MCDI$ and that $(X,r)$ is compact. Then a Feller-continuous strong Markov process with values in $\bM$ is given a.\,s.\ by $(\es{X_t,\rho,\mu_t},t\in\R_+)$.
\end{prop}
By Remark \ref{Pathw:rem:compact-CDI}, the assumption $\Xi\in\MCDI$ ensures that a.\,s., the spaces $X_t$ are compact for all $t\in\R_+$. The proofs of Propositions \ref{Pathw:thm:proc-nd-MvFV} and \ref{Pathw:thm:proc-nd-Mvec} are given further below in this subsection.

We call the process in Proposition~\ref{Pathw:thm:proc-nd-MvFV} an $\M$-valued $\Xi$-Fleming-Viot process with initial state $\ew{X,r,\mu}$.
By Remark \ref{Pathw:rem:psi} below, this process is the $\U$-valued $\Xi$-Fleming-Viot process from \cite{Sampl}*{Section 7.1} if $\rho_0$ is a semi-ultrametric.
We call the process in Proposition~\ref{Pathw:thm:proc-nd-Mvec} an $\bM$-valued evolving $\Xi$-coalescent starting from $\es{X,r,\mu}$.
Note that $\ew{Z,\rho,\mu_t}$ in Proposition~\ref{Pathw:thm:proc-nd-MvFV} above can be replaced by $\ew{X_t,\rho,\mu_t}$ as $\mu_t$ is supported by $X_t$ for all $t\in\R_+$ a.\,s.
\begin{rem}[Isomorphy classes and strong isomorphy classes]
\label{Pathw:rem:iso-strong-iso}
In this remark, we assume $\Xi\in\MCDI$ and that $(X,r)$ is compact. By construction, the $\bM$-valued $\Xi$-Fleming-Viot process in Proposition~\ref{Pathw:thm:proc-nd-Mvec} depends on $(X,r,\mu)$ only through the isomorphy class $\ew{X,r,\mu}$.
In particular, it does not depend on $X\setminus\supp\mu$. As $\supp\mu_0=X_0$ a.\,s.\ by Lemma~\ref{Pathw:lem:mu0}, the strong isomorphy class $\es{X,r,\mu}$ is not necessarily the initial state.

Moreover, let $\Mc$ be the space of isomorphy classes of compact metric measure spaces, and let
$f:\Mc\to\bM,\quad \ew{X',r',\mu'}\mapsto\es{\supp \mu',r',\mu'}$
be the function that maps an isomorphy class to the strong isomorphy class of a representative where the measure has full support. Using Theorem~\ref{Pathw:thm:ld-CDI}\ref{Pathw:item:thm:ld-CDI-supp}, we then obtain that $(f(\chi_t),t\in\R_+)$ is an $\bM$-valued evolving $\Xi$-coalescent if $(\chi_t,t\in\R_+)$ is an $\M$-valued $\Xi$-Fleming-Viot process with initial state $\chi\in\Mc$. Conversely, let $g:\bM\to\M$, $\es{X',r',\mu'}\mapsto\ew{X',r',\mu'}$. Then $(g(\X_t),t\in\R_+)$ is an $\M$-valued $\Xi$-Fleming-Viot process if $(\X_t,t\in\R_+)$ is an $\bM$-valued evolving $\Xi$-coalescent.
Also note that $g$ is continuous.
By Remark~\ref{Pathw:rem:meas-f} below, the function $f$ is measurable.
\end{rem}

We denote by $\Dd$\label{Pathw:not:Dd} the space of semi-metrics on $\N$. We do not distinguish between a semi-metric $\rho'\in\Dd$ and the distance matrix $(\rho'(i,j))_{i,j\in\N}$ and we consider $\Dd$ as a subspace of $\R^{\N^2}$ which we endow with the product topology. Let the $\Dd$-valued Markov process $(\rho_t,t\in\R_+)$ be defined from $\eta$ and $\rho_0$ as in equation~\reff{Pathw:eq:def-rhot}. For each $t\in\R_+$, the random variable $\rho_t$ is exchangeable by \cite{Sampl}*{Proposition 5.8}, that is, $\rho_t$ and $p(\rho_t)$ are equal in distribution for all bijections $p:\N\to\N$. Here the action of the group of bijections $\N\to\N$ on $\Dd$ is defined by $p(\rho')=(\rho'(p(i),p(j)))_{i,j\in\N}$ for $\rho'\in\Dd$ and a bijection $p:\N\to\N$. 
\begin{rem}
\label{Pathw:rem:psi}
Recall the measurable map $\psi:\Dd\to\M$ from \cite{Sampl}*{Section 3.3}.
By construction and Theorem~\ref{Pathw:thm:ld-nd}\ref{Pathw:eq:thm:ld-nd:conv}, we have
$\ew{Z,\rho,\mu_t}=\psi(\rho_t)$ for all $t\in\R_+$ a.\,s.
\end{rem}
\begin{proof}[Proof of Proposition~\ref{Pathw:thm:proc-nd-MvFV} (beginning)]
By Remark \ref{Pathw:rem:psi}, $\ew{Z,\rho,\mu_t}$ is a random variable. The Markov property follows, precisely as in \cite{Sampl}*{Theorem 4.1}, from an application of Theorem 2 of Rogers and Pitman \cite{RP81} to the Markov process $(\rho_t,t\in\R_+)$, the measurable map $\psi:\Dd\to\M$, and the probability kernel from $\M$ to $\Dd$ given by $(\chi,B)\mapsto\nu^\chi(B)$. Here we use the exchangeability of $\rho_t$.

Feller continuity can be shown as in Corollary 8.2 of \cite{Sampl}.
\end{proof}

The proof of Proposition~\ref{Pathw:thm:proc-nd-Mvec} is analogous. Instead of the map $\psi$, we need another map $\upsilon$ which we now define.
Let $\Dd_c\subset\Dd$ be the space of totally bounded semi-metrics on $\N$,
\[\Dd_c=\{\rho'\in\Dd:\lim_{n\to\infty}\sup_{j>n}\inf_{i\leq n}\rho'(i,j)=0\}.\]
Let $\upsilon:\Dd_c\to\bM$ be the function that maps $\rho'\in\Dd_c$ to the strong isomorphy class $\es{X',\rho',\mu'}$ of the compact metric measure space $(X',\rho',\mu')$ defined as follows:
$(X',\rho')$ is the completion of the metric space obtained by identifying the elements with $\rho'$-distance zero in $(\N,\rho')$. The probability measure $\mu'$ on $(X',\rho')$ is the weak limit
$\text{w-}\lim_{n\to\infty} n^{-1}\sum_{i=1}^n\delta_i$
if it exists, else we set $\mu'=\delta_1$.
The following lemma is analogous to \cite{Sampl}*{Proposition 3.7}. 
\begin{lem}
\label{Pathw:lem:upsilon}
The function $\upsilon:\Dd_c\to\bM$ is measurable.
\end{lem}
Again, we refer by an element of a semi-metric space also to the corresponding element of the completion of the metric space that is obtained by identifying points with distance zero.
\begin{proof}
For $n\in\N$, let $\Dd_n\subset\R^{n^2}$ be the space of semi-metrics on $\sn$, again we do not distinguish between semi-metrics and distance matrices. Let $\upsilon_n:\Dd_n\to\bM$ be the function that maps $\rho'\in\Dd_n$ to the strong isomorphy class of the metric measure space $(X',\rho',n^{-1}\sum_{i=1}^n\delta_n)$, where $(X',\rho')$ is the metric space obtained by identifying the elements of $\sn$ with $\rho'$-distance zero. Clearly, the map $\upsilon_n$ is continuous. To show this formally, we define analogously a metric measure space $(X'',\rho'',n^{-1}\sum_{i=1}^n\delta_i)$ from another $\rho''\in\Dd_n$. From \cite{Miermont09}*{Proposition 6}, it follows
\[\dGHP((X',\rho',n^{-1}\sum_{i=1}^n\delta_i),(X'',\rho'',n^{-1}\sum_{i=1}^n\delta_i))
\leq\tfrac{1}{2}\max_{i,j\leq n}\left|\rho'(i,j)-\rho''(i,j)\right|,\]
we use the coupling $\nu=n^{-1}\sum_{i=1}^n\delta_{(i,i)}$ on $X'\times X''$ and the correspondence $\mathfrak{R}=\{(i,i):i\in\sn\}\subset X'\times X''$. 

Now let $\rho'\in\Dd_c$, and let $(X',\rho',\mu')$ be defined as in the definition of $\upsilon(\rho')$ above. Using the definition of the Gromov-Hausdorff-Prohorov metric, we obtain that
\[\dGHP(\upsilon(\rho'),\upsilon_n((\rho'(i,j))_{i,j\leq n}))
\leq \dP^{X'}(\mu',n^{-1}\sum_{i=1}^n\delta_i)\vee\dH^{X'}(X',\sn)\to 0\quad (n\to\infty)\]
if the weak limit $\mu'$ of the measures $n^{-1}\sum_{i=1}^n\delta_i$ on $X'$ exists. This yields the assertion.
\end{proof}
\begin{proof}[Proof of Proposition~\ref{Pathw:thm:proc-nd-Mvec} (beginning)]
By Remark \ref{Pathw:rem:compact-CDI}, there exists an event of probability $1$ on which $\rho_t\in\Dd_c$ for all $t\in\R_+$.
By construction and Theorems~\ref{Pathw:thm:ld-nd}\ref{Pathw:eq:thm:ld-nd:conv} and \ref{Pathw:thm:ld-CDI}\ref{Pathw:item:thm:ld-CDI-supp}, we have
$\es{X_t,\rho,\mu_t}=\upsilon(\rho_t)$ for all $t\in\R_+$ a.\,s.
Hence, Lemma~\ref{Pathw:lem:upsilon} yields that $\es{X_t,\rho,\mu_t}$ is a random variable. The Markov property follows as in \cite{Sampl}*{Theorem 4.1} from an application of \cite{RP81}*{Theorem 2} to the Markov process $(\rho_t,t\in\R_+)$, the measurable map $\upsilon:\Dd_c\to\bM$, and the probability kernel from $\bM$ to $\Dd_c$ given by $(\chi,B)\mapsto\nu^\chi(B)$. Here we use exchangeability of $\rho_t$.
To check that Condition (a) in \cite{RP81}*{Theorem 2} is satisfied, we note that $\upsilon(\rho')=\upsilon(\rho_t)$ a.\,s.\ for $t\in\R_+$ and a random variable $\rho'$ with conditional distribution $\nu^{\upsilon(\rho_t)}$ given $\upsilon(\rho_t)$. This a.\,s.\ equality follows as in the proof of \cite{Sampl}*{Proposition 10.5}, we also use Theorem~\ref{Pathw:thm:ld-CDI}\ref{Pathw:item:thm:ld-CDI-supp}.
\end{proof}
\begin{rem}
\label{Pathw:rem:meas-f}
As a by-product of Lemma~\ref{Pathw:lem:upsilon}, let us deduce measurability of the canonical map $f:\M_c\to\bM$ from Remark~\ref{Pathw:rem:iso-strong-iso} (this answers a question posed to the author by H.\,Sulzbach).
We consider $\chi\in\M_c$ and a random variable $\rho$ with the distance matrix distribution $\nu^\chi$. As in the proof of \cite{Sampl}*{Proposition 10.5}, it follows that $f(\chi)=\upsilon(\rho)$ a.\,s.
Hence, for a Borel subset $A\subset\bM$, we obtain the equivalence
\[\I{f(\chi)\in A}=1 \quad\Leftrightarrow\quad
\int\nu^{\chi}(d\rho')\I{\upsilon(\rho')\in A}=1
\quad\Leftrightarrow\quad
\nu^\chi(\upsilon^{-1}(A))=1.
\]
The function that maps a metric measure space $\chi\in\M_c$ to its distance matrix distribution $\nu^\chi$ is continuous by definition of the Gromov-Prohorov topology (see \cite{GPW09}). Lemma~\ref{Pathw:lem:upsilon} now implies that $\{\chi\in\M_c:f(\chi)\in A\}$ is a measurable subset of $\M_c$.
We remark that measurability of $f$ can also be obtained as a consequence of \cite{ALW16}*{Corollary 5.6} and e.\,g.\ \cite{Kechris}*{Theorem 15.1}.
\end{rem}

Path regularity of the $\M$-valued $\Xi$-Fleming-Viot process and the $\bM$-valued evolving $\Xi$-coalescent follows from Section~\ref{Pathw:sec:ld-sampl-nd}:
\begin{prop}
\label{Pathw:cor:proc-nd-MvFV}
Assume $\Xi\in\nd$. Then a.\,s., the process $(\ew{Z,\rho,\mu_t},t\in\R_+)$ has càdlàg paths in the Gromov-weak topology and $\Theta_0$ is the set of jump times.
\end{prop}
\begin{prop}
\label{Pathw:cor:proc-nd-Mvec}
Assume $\Xi\in\MCDI$. Then a.\,s., the process $(\es{X_t,\rho,\mu_t},t\in\R_+)$ has càdlàg paths in the Gromov-Hausdorff-Prohorov topology and $\Theta_0\cup\Thetaext$ is the set of jump times.
\end{prop}
\begin{proof}[Proof of Proposition~\ref{Pathw:cor:proc-nd-MvFV}]
By Theorem~\ref{Pathw:thm:ld-nd}\ref{Pathw:item:thm:ld-nd:reg} and the definition of the Gromov-Prohorov metric,
it follows that a.\,s., the map $t\mapsto\ew{Z,\rho,\mu_t}$ is càdlàg and the set of jump times is not larger than $\Theta_0$. By Theorem~\ref{Pathw:thm:ld-nd}\ref{Pathw:item:thm:ld-nd:atoms-jumps} and as the atomicity properties only depend on the isomorphy classes, it follows that a.\,s., the set of jump times is not smaller than $\Theta_0$.
\end{proof}
\begin{proof}[Proof of Proposition~\ref{Pathw:cor:proc-nd-Mvec}]
This is analogous to Proposition \ref{Pathw:cor:proc-nd-MvFV}. We use Theorems~\ref{Pathw:thm:ld-nd}\ref{Pathw:item:thm:ld-nd:reg} and \ref{Pathw:thm:ld-CDI}\ref{Pathw:item:thm:ld-CDI-reg}, the definition of the Gromov-Hausdorff-Prohorov metric, and the atomicity properties from Theorem~\ref{Pathw:thm:ld-nd}\ref{Pathw:item:thm:ld-nd:atoms-jumps} which are determined by the strong isomorphy classes. We also use that $X_{t-}$ and $X_t$ are isometric if $\es{X_{t-},\rho,\mu_{t-}}=\es{X_t,\rho,\mu_t}$.
\end{proof}
\begin{proof}[Proof of Proposition \ref{Pathw:thm:proc-nd-MvFV} (end)]
The strong Markov property can be deduced by standard arguments (cf.\ e.\,g.\ the proof of Theorem 4.2.7 in \cite{EK86}) from Feller continuity and a.\,s.\ right continuity of the sample paths (Proposition \ref{Pathw:cor:proc-nd-MvFV}).
\end{proof}

Now we study Feller continuity of the $\bM$-valued evolving $\Xi$-coalescent.
\begin{lem}
\label{Pathw:lem:Feller-Mvec}
Let $(X^n,r^n,\mu^n)$ be a sequence of compact metric measure spaces such that $\ew{X^n,r^n,\mu^n}$ converges to $\ew{X,r,\mu}$ in the Gromov-weak topology. Assume $\Xi\in\MCDI$ and let $(\X^n_t,t\in\R_+)$ be an $\bM$-valued evolving $\Xi$-coalescent starting from $\es{X^n,r^n,\mu^n}$. Then for each $t\in(0,\infty)$, the random variable $\X^n_t$ converges in distribution to $\es{X_t,\rho,\mu_t}$ in $\bM$, endowed with the Gromov-Hausdorff-Prohorov topology.
\end{lem}
\begin{proof}
Let $t,\ep>0$ and $n\in\N$. Let $\rho^n_0$ be a random variable with distribution $\nu^{\ew{X^n,r^n,\mu^n}}$ that is independent of $\eta$. Recall the definition of $\rho_0$ from Section~\ref{Pathw:sec:ld-sampl-nd}. Let $(Z',\rho')$ be the lookdown space associated with $\eta$ and $\rho^n_0$.
Let $X'_t\subset Z'$ be the closure of the set $\{t\}\times\N$ of individuals at time $t$ therein, and define a probability measure $\mu'_t$ on $Z'$ analogously to $\mu_t$. Then a.\,s., the map $X_t\to X'_t$, $(t,i)\mapsto (t,i)$ can be extended to a measure-preserving homeomorphism $h$ between $(X_t,\rho,\mu_t)$ and $(X'_t,\rho',\mu'_t)$. The correspondence
$\mathfrak{R}=\{(x,h(x)):x\in X_t\}\subset X_t\times X'_t$ has distortion $\max\{|\rho_0^n(i,j)-\rho_0(i,j)|:i,j\in A_0(t,\N)\}$, where we write
$A_0(t,\N)=\{A_0(t,\ell):\ell\in\N\}$.
With the coupling $\nu(dx\,dx')=\mu_t(dx)\delta_{h(x)}(dx')$ of $\mu_t$ and $\mu'_t$, Proposition 6 in \cite{Miermont09} implies
\begin{align*}
&\P(\dGHP(\es{X'_t,\rho',\mu'_t},\es{X_t,\rho,\mu_t})\geq\ep)\\
\leq\; &\P(\max\{|\rho_0^n(i,j)-\rho_0(i,j)|:i,j\in\sk\}\geq 2\ep)
+\P(\#A_0(t,\N)>k)
\end{align*}
for all $k\in\N$.
W.\,l.\,o.\,g., we may assume $\X^n_t=\es{X'_t,\rho',\mu'_t}$ for all $t\in\R_+$ a.\,s.,
and that the distance matrices $\rho_0^n$ converge in probability. We let $n$ and then $k$ tend to infinity.
\end{proof}
\begin{proof}[Proof of Proposition \ref{Pathw:thm:proc-nd-Mvec} (end)]
As the map $g:\bM\to\M$ in Remark \ref{Pathw:rem:iso-strong-iso} is continuous, we can make in Lemma~\ref{Pathw:lem:Feller-Mvec} also the stronger assumption that $\es{X^n,r^n,\mu^n}$ converges to $\es{X,r,\mu}$ in the Gromov-Hausdorff-Prohorov topology. This yields Feller continuity of $([X_t,\rho,\mu_t],t\in\R_+)$. The strong Markov property can now be deduced using a.\,s.\ right continuity of the sample paths (Proposition \ref{Pathw:cor:proc-nd-Mvec}).
\end{proof}
The tree-valued $\Xi$-Fleming-Viot process converges to equilibrium, as shown in \cite{Sampl}*{Proposition 9.1}. Now we show a similar result for the $\bM$-valued evolving $\Xi$-coalescent.

Assume $\Xi\in\MCDI$ and (analogously to Section \ref{Pathw:sec:ld-Xi}), let $\bar\eta$ be a Poisson random measure on $\R\times\p$ with intensity $dt\;H_\Xi(d\pi)$.
From $\bar\eta$, we define a lookdown space $(\bar Z,\bar \rho)$ in two-sided time as the completion of the space of individuals $\R\times\N$ with respect to the semi-metric $\bar\rho$, given by
\begin{equation}
\label{Pathw:eq:barrho}
\bar\rho((t,i),(u,j))=t+u-2\sup\{s\in(-\infty,t\wedge u]:\bar A_s(t,i)=\bar A_s(u,j)\},
\end{equation}
where $\bar A_s(t,i)$ denotes the level of the ancestor of the individual $(t,i)$ when particles and reproduction events are defined precisely as in Section~\ref{Pathw:sec:ld-space}.

Analogously to Theorem~\ref{Pathw:thm:ld-nd}, on an event of probability $1$, the probability measures
$\bar\mu^n_t=n^{-1}\sum_{i=1}^n\delta_{(t,i)}$
on $(\bar Z,\bar\rho)$ weakly converge as $n\to\infty$ for all $t\in\R$, we denote the limits by $\bar\mu_t$. For $t\in\R$, we denote by $\bar X_t$ the closure of $\{t\}\times\N$ in $(\bar Z,\bar\rho)$. A stationary $\bM$-valued evolving $\Xi$-coalescent is given by
$(\es{\bar X_t,\bar\rho,\bar\mu_t},t\in\R)$. We call a random variable that is distributed as $\es{\bar X_0,\bar\rho,\bar\mu_0}$ an $\bM$-valued $\Xi$-coalescent measure tree, in analogy to the coalescent measure trees in \cites{Sampl,GPW13}.
As $\es{\bar X_0,\bar\rho,\bar\mu_0}$ is a.\,s.\ an ultrametric measure space, this random variable can be seen as a random tree.

In the next proposition, we show that the $\bM$-valued $\Xi$-Fleming-Viot process $(\es{X_t,\rho,\mu_t},t\in\R_+)$ that is defined from $\eta$ and $\es{X,r,\mu}$ in this section converges to equilibrium.
\begin{prop}
\label{prop:Mvec-equil}
The $\bM$-valued random variable $\es{X_t,\rho,\mu_t}$, $t\in\R_+$ converges in distribution in $(\bM,\dGHP)$ to an $\bM$-valued $\Xi$-coalescent measure tree as $t\to\infty$.
\end{prop}
As in \cites{Sampl,DK96}, we use a coupling argument in the proof. In the present context, the topology is stronger than in \cite{Sampl}, but as we restrict to $\Xi\in\MCDI$, there exists a coupling of the tree-valued evolving $\Xi$-coalescents with arbitrary initial state and of the stationary process such that these processes coincide after an a.\,s.\ finite random time.
\begin{proof}[Proof of Proposition \ref{prop:Mvec-equil}]
Assume that the Poisson random measure $\eta$ is the restriction of $\bar\eta$ to $(0,\infty)\times\p$. Then $\es{\bar X_t,\bar\rho,\bar\mu_t}=\es{X_t,\rho,\mu_t}$ on the event $\{\diam X_t<2t\}$.
By the properties of the Poisson random measure $\eta$ and as the event $\{\diam X_1<2\}$ is independent of $\rho_0$, the events $\{\diam X_t<2\}$, $t\in\N$ are independent and have the same positive probability.
Hence, the random time $\tau=\inf\{t\in\R_+:\diam X_t<2t\}$ is geometrically bounded. The assertion follows as $\diam X_t< 2t$ for all $t> \tau$, and as $\es{\bar X_t,\bar\rho,\bar\mu_t}$ is an $\bM$-valued $\Xi$-coalescent measure tree.
\end{proof}
\begin{rem}[Convergence of $\M$-valued $\Xi$-Cannings processes]
\label{Pathw:rem:tvC}
As an immediate consequence of the uniform convergence in Theorem~\ref{Pathw:thm:ld-nd}\ref{Pathw:eq:thm:ld-nd:conv}, we obtain
\begin{equation}
\label{Pathw:eq:conv-dGP}
\lim_{n\to\infty}\sup_{t\in[0,T]}\dGP(\ew{Z,\rho,\mu^n_t},\ew{Z,\rho,\mu_t})=0\quad\text{a.\,s.}
\end{equation}
for each $T\in\R_+$. The process $(\ew{Z,\rho,\mu^n_t},t\in\R_+)$ may be called an $\M$-valued $\Xi$-Cannings process. In the case without simultaneous multiple reproduction events, it coincides with the tree-valued $\Lambda$-Cannings process discussed in \cite{KL14}*{Section 4.2},
and in the case without multiple reproduction events with the tree-valued Moran process from \cite{GPW13}*{Definition 2.19}. This can be seen by an application of \cite{RP81}*{Theorem 2} similarly to the proof of Lemma~\ref{Pathw:lem:exch-J} below, see also Section 2 of \cite{DKW14}. Then the convergence~\reff{Pathw:eq:conv-dGP} implies the assertion of Theorem 2 in \cite{GPW13} for a special choice of the approximating sequence of the initial state.

From the uniform convergence in Theorem~\ref{Pathw:thm:ld-dust}\ref{Pathw:eq:thm:ld-dust:conv}, similar statements can be deduced for processes from the next subsection. In \cite{Conv}, convergence of tree-valued Cannings chains is studied by different methods.
\end{rem}

\subsection{The general case}
\label{Pathw:sec:proc-gen}
We include the case with dust by using marked metric measure spaces. Let us first recall some facts from \cite{DGP11}. An $(\R_+)$-marked metric measure space $(X,r,m)$ is a triple that consists of a complete and separable metric space $(X,r)$ and a probability measure $m$ on the Borel sigma algebra on the product space $X\times\R_+$. Two marked metric measure spaces $(X,r,m)$, $(X',r',m')$ are called isomorphic if there exists an isometry $\varphi$ between the closed supports $\supp m(\cdot\times\R_+)\subset X$ and $\supp m'(\cdot\times\R_+)\subset X'$ such that the measurable map $\hat\varphi:\supp m\to \supp m'$, given by $\hat\varphi(x,v)=(\varphi(x),v)$, satisfies $m'=\hat\varphi(m)$. We endow the space $\hat\M$ of isomorphy classes of marked metric measure spaces with the marked Gromov-Prohorov metric, which is defined by
\begin{equation*}
\dMGP((X,r,m),(X',r',m'))=\inf_{Y,\hat\varphi,\hat\varphi'}\dP^Y(\hat\varphi(m),\hat\varphi'(m'))
\end{equation*}
where the infimum is over all isometric embeddings $\varphi:X\to Y$, $\varphi':X'\to Y$ into complete and separable metric spaces $(Y,d^Y)$.
The maps $\hat\varphi:X\times\R_+\to Y\times\R_+$ and $\hat\varphi':X'\times\R_+\to Y\times\R_+$ are defined by $\hat\varphi(x,v)=(\varphi(x),v)$ and $\hat\varphi'(x',v)=(\varphi'(x'),v)$.
The space $Y\times\R_+$ is endowed with the product metric $d^{Y\times\R_+}((y,v),(y',v'))=d^Y(y,y')\vee|v-v'|$.
Then $(\hat\M,\dMGP)$ is a complete and separable metric space. 
The marked distance matrix distribution $\nu^{(X,r,m)}$ of (an isomorphy class of) a marked metric measure space $(X,r,m)$ is defined as the distribution of the random variable $(r(x(i),x(j))_{i,j\in\N},(v(i))_{i\in\N})$ where $(x(i),v(i))_{i\in\N}$ is an $m$-iid sequence in $X\times\R_+$. The metric $\dMGP$ induces the marked Gromov-weak topology in which a sequence of marked metric measure spaces converges if and only if their marked distance matrix distributions converge weakly.

From the construction on the extended lookdown space in Section~\ref{Pathw:sec:ld-sampl-dust}, we now read off a stochastic processes with values in $(\hat\M,\dMGP)$.

We work with the extended lookdown space $(\hat Z,\rho)$ as defined in Section~\ref{Pathw:sec:ld-sampl-dust}. This random metric space is constructed from the Poisson random measure $\eta$ which is characterized by a finite measure $\Xi$ on the simplex $\Delta$, and from an independent random variable $(r_0,v_0)$ with the marked distance matrix distribution of a marked metric measure space $(X,r,\mu)$. The following proposition, which is proved below in this subsection, is the analog of Proposition \ref{Pathw:thm:proc-nd-MvFV}.

\begin{prop}
\label{Pathw:thm:proc}
Assume that $\Xi$ is a finite measure on $\Delta$. Then a Feller-continuous strong Markov process with values in $\hat\M$ is given a.\,s.\ by $(\ew{\hat Z,\rho,m_t},t\in\R_+)$.
\end{prop}
We call $(\ew{\hat Z,\rho,m_t},t\in\R_+)$ an $\hat\M$-valued $\Xi$-Fleming-Viot process with initial state $\ew{X,r,m}$.
By Remark \ref{Pathw:rem:hatpsi} below, this process is the tree-valued $\Xi$-Fleming-Viot process from \cite{Sampl}*{Section 7.2} if the restriction of $\rho$ to $\{0\}\times\N$ is ultrametric.

Recall the space $\Dd$ of semi-metrics on $\N$. We define the space $\hat\Dd=\Dd\times\R_+^\N\subset\R^{\N^2}\times\R^\N$\label{Pathw:not:hatDd}, where $\R^{\N^2}\times\R^\N$ is endowed with the product topology. The elements of $\hat\Dd$ are called marked distance matrices or decomposed semi-metrics on $\N$.
We define a $\hat\Dd$-valued stochastic process $((r_t,v_t),t\in\R_+)$ from the Poisson random measure $\eta$ and $(r_0,v_0)$ as in Remark~\ref{Pathw:rem:rho-r-v}. By Proposition 6.3 of \cite{Sampl}, this process is Markovian. In Lemma~\ref{Pathw:lem:exch-Rt}, we will show that for each $t\in\R_+$, the $\hat\Dd$-valued random variable $(r_t,v_t)$ is exchangeable. That is, its distribution is invariant under the action of the group of bijections $\N\to\N$, defined by
$p(r,v)=((r(p(i),p(j)))_{i,j\in\N},(v(p(i)))_{i\in\N})$ for $(r,v)\in\hat\Dd$ and any bijection $p:\N\to\N$.
\begin{rem}
\label{Pathw:rem:hatpsi}
Recall the measurable map $\hat\psi:\hat\Dd\to\hat\M$ from \cite{Sampl}*{Section 3.3}. Theorem \ref{Pathw:thm:ld-dust} yields $\ew{\hat Z,\rho,m_t}=\hat\psi(r_t,v_t)$ for all $t\in\R_+$ a.\,s.
\end{rem}
\begin{proof}[Proof of Proposition~\ref{Pathw:thm:proc} (beginning)]
The Markov property follows as in \cite{Sampl}*{Theorem 4.1}, we apply \cite{RP81}*{Theorem 2} to the process $((r_t,v_t),t\in\R_+)$, the measurable map $\hat\psi:\hat\Dd\to\hat\M$, and the probability kernel from $\hat\M$ to $\hat\Dd$, given by $(\chi,B)\mapsto\nu^\chi(B)$. Here we use the exchangeability of $(r_t,v_t)$.

Feller continuity can be shown precisely as in Corollary 8.2 of \cite{Sampl}.
\end{proof}
Let us deduce path regularity:
\begin{prop}
\label{Pathw:cor:proc}
Assume that one of the following conditions hold: (i) $\Xi\in\dust$, or (ii) $\Xi\in\nd$ and $m(X\times\{0\})=1$. Then the process $(\ew{\hat Z,\rho,m_t},t\in\R_+)$ has a.\,s.\ càdlàg paths in the marked Gromov-weak topology and $\Theta_0$ is the set of jump times.
\end{prop}
\begin{proof}
We argue as in the proof of Proposition~\ref{Pathw:cor:proc-nd-MvFV}. In case $\Xi\in\dust$, we use Theorem~\ref{Pathw:thm:ld-dust}. If $\Xi\in\nd$ and $m(\hat Z\times\{0\})=1$, we use Remark~\ref{Pathw:rem:Z-consistence-random} and Theorem~\ref{Pathw:thm:ld-nd}.
\end{proof}
\begin{prop}
\label{Pathw:cor:proc-0}
Assume $\Xi\in\nd$ and $m(X\times\{0\})<1$. Then the process $(\ew{\hat Z,\rho,m_t}, t\in(0,\infty))$ has a.\,s.\ càdlàg paths in the marked Gromov-weak topology and $\Theta_0$ is the set of jump times. The process $(\ew{\hat Z,\rho,m_t},t\in\R_+)$ is a.\,s.\ not right-continuous at time $0$.
\end{prop}
\begin{proof}
The first assertion follows from Corollary~\ref{Pathw:cor:nd-mdm}.
The definitions of $m_t$ and of the marked distance matrix distribution, and Remark \ref{Pathw:rem:v-nd} yield that
$\nu^{\ew{\hat Z,\rho,m_t}}(\Dd\times\{0\})=m_t(\hat Z\times\{0\})=1$
for all $t\in(0,\infty)$ a.\,s.
As the marked metric measure spaces $(Z,\rho,m_0)$ and $(X,r,m)$ are a.\,s.\ isomorphic (see e.\,g.\ \cite{Sampl}*{Proposition 10.5}), the assumptions also yield
$\nu^{\ew{\hat Z,\rho,m_0}}(\Dd\times\{0\})=m_0(\hat Z\times\{0\})=m(X\times\{0\})<1$ a.\,s.
Hence, the probability measures $\nu^{\ew{\hat Z,\rho,m_t}}$ do a.\,s.\ not weakly converge to $\nu^{\ew{\hat Z,\rho,m_0}}$. By definition of the marked Gromov-weak topology, it follows that $t\mapsto\ew{\hat Z,\rho,m_t}$ is a.\,s.\ not right-continuous at $0$.
\end{proof}
\begin{proof}[Proof of Proposition \ref{Pathw:thm:proc} (end)]
The strong Markov property can be deduced from Feller continuity and the a.\,s.\ right continuity of the sample paths that is given by Propositions \ref{Pathw:cor:proc} and \ref{Pathw:cor:proc-0}.
\end{proof}

From the construction in Section \ref{Pathw:sec:ld-sampl-dust}, we also read off the $\Xi$-Fleming-Viot process with values in a space of matrix distributions from \cite{Sampl}*{Section 7.3}. As in \cite{Sampl}*{Section 2}, let
\begin{equation}
\label{Pathw:eq:def-alpha}
\alpha:\hat\Dd\to\Dd,\quad (r,v)\to ((r(i,j)+v(i)+v(j))\I{i\neq j})_{i,j\in\N}
\end{equation}
be the continuous function that maps a decomposed semi-metric to the corresponding semi-metric.
The distance matrix distribution of a marked metric measure space $\chi$ is defined as $\alpha(\nu^\chi)$, where $\nu^\chi$ denotes the marked distance matrix distribution of $\chi$. It depends only on the isomorphy class of the marked metric measure space.
Let $\rho_0=\alpha(r_0,v_0)$, then $\rho_0$ is the restriction of $\rho$ to $\{0\}\times\N$. As we will apply \cite{Sampl}*{Proposition 3.4}, we assume in the remainder of this subsection that $\rho_0$ is a semi-ultrametric (in which case the restriction of $\rho$ to $\{0\}\times\N$ is ultrametric).
Then the stochastic process of the distance matrix distributions of the $\hat\U$-valued $\Xi$-Fleming-Viot process, given by
$(\xi_t,t\in\R_+)=(\alpha(\nu^{\ew{\hat Z,\rho,m_t}}),t\in\R_+)$
a.\,s., is the $\UUerg$-valued $\Xi$-Fleming-Viot process from \cite{Sampl}*{Section 7.3} where $\UUerg$ refers to the space
\[\UUerg=\{\alpha(\nu^\chi):\chi\in\hat\M,\alpha(\nu^\chi)\text{-a.\,a.\ $\rho\in\Dd$ are semi-ultrametrics}\}\]
of ultrametric distance matrix distributions of marked metric measure spaces.
(This notation is used in \cite{Sampl} as $\UUerg$ is shown there to be the set of distributions of semi-ultrametrics that are invariant and ergodic under the action of the group of finite permutations.)
We endow $\UUerg$ with the Prohorov metric.

In spite of Proposition~\ref{Pathw:cor:proc-0}, the $\UUerg$-valued $\Xi$-Fleming-Viot process is always a.\,s.\ right-continuous at time $0$ by the following proposition which is applied in \cite{Conv}.
\begin{prop}
Assume that $\Xi$ is a finite measure on $\Delta$. Then
$(\xi_t,t\in\R_+)$ has a.\,s.\ càdlàg paths and $\Theta_0$ is a.\,s.\ the set of jump times.
\end{prop}
\begin{proof}
Propositions \ref{Pathw:cor:proc} and \ref{Pathw:cor:proc-0}, the definition of the marked Gromov-weak topology, and continuity of $\alpha$ yield that on an event of probability $1$, the path $(0,\infty)\to\UUerg$, $t\mapsto\xi_t$ is a.\,s.\ càdlàg in $\dP$ with a.\,s.\ no jumps outside $\Theta_0$, and $\ew{\hat Z,\rho,m_t}\neq\ew{\hat Z,\rho,m_{t-}}$
for all $t\in\Theta_0$. By Propositions 3.4 and 10.5 in \cite{Sampl},
it follows that $\xi_t\neq\xi_{t-}$ for all $t\in\Theta_0$ a.\,s.

Right continuity at time $0$ follows from Proposition \ref{Pathw:cor:proc} under the assumptions therein.
For the general case, we use that $(\xi_t,t\in\R_+)$ solves the martingale problem $(C,\Cc)$ defined in \cite{Sampl}*{Section 7.3}. We briefly recall the definition of the domain $\Cc$. For $n\in\N$, let $\gamma_n$ be the restriction from $\R^{\N^2}$ to $\R^{n^2}$, $\gamma_n(\rho')=(\rho'(i,j))_{i,j\in\sn}$. Let $\C_n$ be the set of functions $\phi\circ\gamma_n:\R^{\N^2}\to\R$, where $\phi$ is a bounded differentiable function $\R^{n^2}\to\R$ with bounded uniformly continuous derivative. Then we set $\C=\bigcup_{n\in\N}\C_n$ and $\Cc=\{\UUerg\to\R,\xi\mapsto\xi\phi:\phi\in\C\}$, where we use the notation $\xi\phi=\int\xi(d\rho')\phi(\rho')$.
By definition of $C$ in \cite{Sampl}*{Section 7.3}, $C\Psi$ is bounded for each $\Psi\in\Cc$.

The set $\Cc$ is convergence determining in $\UUerg$, see \cite{Sampl}*{Remark 4.6}. There also exists a countable subset of $\Cc$ that generates the weak topology on $\UUerg$. Indeed, by smoothing indicator functions of rational intervals, one finds a countable subset $\C'\subset\C$ that generates the product topology on $\R^{\N^2}$. Let $\C''$ be the set of finite products of functions in $\C'$. Then the algebra $\C''$ is convergence determining in $\R^{\N^2}$ by e.\,g.\ \cite{Lohr13}*{Theorem 2.7}. That is, $\Cc':=\{\xi\mapsto\xi\phi:\phi\in\C''\}\subset\Cc$ generates the weak topology on $\UUerg$.

For $\Psi\in\Cc$, we consider the process $(M_t,t\in\R_+)$ that is given by
\[M_t=\Psi(\xi_t)-\int_0^t C\Psi(\xi_s)ds\]
and which is a bounded martingale.
From
\begin{equation*}
\E[\Psi(\xi_t)-\Psi(\xi_0)-\int_0^t C\Psi(\xi_s)ds]=0,
\end{equation*}
continuity of $\Psi$, and as $C\Psi$ is bounded, we obtain that $\xi_t$ converges in distribution to $\xi_0$ as $t\to 0$.
As in the proof of Proposition \ref{Pathw:cor:proc-0}, the marked metric measure spaces $(X,r,m)$ and $(\hat Z,\rho,m_0)$ are a.\,s.\ isomorphic. Hence, $\xi_0=\alpha(\nu^{(X,r,m)})$ a.\,s.
As $(X,r,m)$ is deterministic, it follows that $\xi_t$ converges to $\xi_0$ also in probability.
By martingale convergence
and as $(M_t,t\in(0,\infty))$ has a.\,s.\ càdlàg paths, it follows that the limit $\lim_{t\downarrow 0}M_t$ exists a.\,s., see e.\,g.\ \cite{RW1}*{Theorem II.69.4}.
By the convergence in probability we already know, the definition of $M_t$, and as $C\Psi$ is bounded, the limit must be $\Psi(\xi_0)$ a.\,s., hence $\Psi(\xi_t)$ converges to $\Psi(\xi_0)$ a.\,s. As $\Cc$ contains a countable subset that generates the weak topology on $\UUerg$, it follows that $\xi_t$ converges to $\xi_0$ a.\,s.
\end{proof}

\section{Outline and some definitions for the proof of the central results}
\label{Pathw:sec:outline-proofs}
The aim of the remaining sections is to prove in Section \ref{Pathw:sec:constr-ld} Theorems~\ref{Pathw:thm:ld-nd}, \ref{Pathw:thm:ld-CDI}, and \ref{Pathw:thm:ld-dust}.
In Section~\ref{Pathw:sec:exch}, we prove exchangeability properties for the (decomposed) genealogical distances between the individuals in the lookdown model at various stopping times.
In Section~\ref{Pathw:sec:unif-ld}, we state convergence results for processes of certain asymptotic frequencies that depend on these (decomposed) genealogical distances.
We apply these convergence results in particular to families of partitions in Section~\ref{Pathw:sec:part}:
In the case without dust, we consider the flow of partitions. For the case with dust, we introduce a family of partitions that fits to the decomposition of genealogical distances.
Using these families of partitions, we construct the probability measures on the (extended) lookdown space in Section~\ref{Pathw:sec:constr-ld}.
We prove the convergence results from Section~\ref{Pathw:sec:unif-ld} in Section~\ref{Pathw:sec:unif-ld-proofs} using the exchangeability properties from Section~\ref{Pathw:sec:exch}.

Now we collect some definitions that we will use in the remaining sections.

\subsection{Some notation}
For $n\in\N$, we continue using the notation $\sn=\{1,\ldots,n\}$. We also write $\N_0=\N\cup\{0\}$ and $[0]=\emptyset$.
Recall the space $\Dd\subset\R^{\N^2}$ of distance matrices, and that we do not distinguish distance matrices from semi-metrics on $\N$. Recall also the space $\hat\Dd=\Dd\times\R_+^\N\subset\R^{\N^2}\times\R^\N$ of marked distance matrices or decomposed semi-metrics on $\N$. Here $\R^{\N^2}$ and $\R^{\N^2}\times\R^\N$ are endowed with the product topology.
We denote by $\Dd_n\subset\R^{n^2}$\label{Pathw:not:Ddn} the space of semi-metrics on $\sn$ which we do not distinguish from distance matrices. We denote by $\hat\Dd_n=\Dd_n\times\R_+^n\subset\R^{n^2}\times\R^n$\label{Pathw:not:hatDdn} the space of decomposed semi-metrics on $\sn$ or marked distance matrices.
Recall also the space $\p$ of partitions of $\N$ which is endowed with the topology induced by the restriction maps. We denote by $\p_n$ of set of partitions of $\sn$. We denote by $\gamma_n$ the restriction maps\label{Pathw:not:gamma-outl} $\gamma_n:\Dd\to\Dd_n$, $\rho\mapsto(\rho(i,j))_{i,j\in\sn}$, $\gamma_n:\hat\Dd\to\hat\Dd_n$, $(r,v)\mapsto((r(i,j))_{i,j\in\sn},(v(i))_{i\in\sn})$, and $\gamma_n:\p\to\p_n$.
For a partition $\pi\in\p$ and $i\in\N$, we denote the block of $\pi$ that contains $i$ by $B(\pi,i)$\label{Pathw:not:Bpii}.

We denote the set of the minimal elements of the blocks of $\pi$ by \label{Pathw:not:Mpi}$M(\pi)=\{\min B:B\in\pi\}$. For a subset $B\subset\N$, we denote the relative frequency by $|B|_n=n^{-1}\#(B\cap\sn)$ for $n\in\N$, and the asymptotic frequency by $|B|=\lim_{n\to\infty}|B|_n$ (if it exists).

\subsection{Two-step construction of the point measure of reproduction events}
\label{Pathw:sec:two-step-eta}
We will use the following definitions from Section~\ref{Pathw:sec:exch-ld} onwards.
As in Section~\ref{Pathw:sec:ld-Xi}, let $\Xi$ be a finite measure on the simplex $\Delta$. We decompose $\Xi=\Xi_0+\Xi\{0\}\delta_0$.
Let $\eta$ be a Poisson random measure on $(0,\infty)\times\p$ with intensity $dt\,H_\Xi(d\pi)$ as in Section~\ref{Pathw:sec:ld-Xi}. We assume w.\,l.\,o.\,g.\ that $\eta$\label{Pathw:not:eta-twostep} is constructed as the sum $\eta=\eta_{\rm K}+\eta_0$\label{Pathw:not:eta0etaK} of two independent Poisson random measures, defined as follows. We define $\eta_{\rm K}$ as a Poisson random measure on $(0,\infty)\times\p$ with intensity
\[dt\; \Xi\{0\}\sum_{1\leq i<j}\delta_{K_{i,j}}(d\pi),\]
where $K_{i,j}$ denotes the partition in $\p$ that contains the block $\{i,j\}$ and apart from that only singletons.
The point measure $\eta_{\rm K}$ encodes the Kingman part, that is, the binary reproduction events. If $\Xi_0(\Delta)=0$, we set $\eta_0=0$.
Let $\zeta_0$\label{Pathw:not:xi0} be a Poisson random measure on $(0,\infty)\times\Delta$ with intensity $dt\,\abs{x}_2^{-2}\Xi_0(dx)$. In case $\Xi_0(\Delta)>0$, let $((t^k,y^k),k\in\N)$ be a collection of $\zeta_0$-measurable random variables with values in $(0,\infty)\times\Delta$ such that
\[\zeta_0=\sum_{k\in\N}\delta_{(t^k,y^k)}\quad\text{a.\,s.}\]
Let $(\pi^k,k\in\N)$ be a collection of $\p$-valued random variables that are conditionally independent given $((t^k,y^k),i\in\N)$ such that $\pi^k$ has conditional distribution $\kappa(y^k,\cdot)$ given $((t^k,y^k),k\in\N)$.
Thereby, $\kappa$ is the probability kernel from $\Delta$ to $\p$ associated with Kingman's correspondence as in Section~\ref{Pathw:sec:ld-Xi}.
We define the Poisson random measure $\eta_0$ on $(0,\infty)\times\p$ by
\begin{equation}
\label{eq:def-eta0}
\eta_0=\sum_{k\in\N}\delta_{(t^k,\pi^k)}.
\end{equation}
For each point $(t,y)$ of $\zeta_0$ and the associated point $(t,\pi)$ of $\eta_0$, the vector $y$ gives the asymptotic frequencies (in decreasing order) of the blocks of the partition $\pi$, that is, the relative family sizes in the large reproduction event at time $t$.

Recall that $\Xi\in\dust$ implies $\Xi\{0\}=0$. Also note that the set $\Theta_0$ of times of large reproduction events, defined in~\reff{Pathw:eq:Theta0}, satisfies
\[\Theta_0=\{t\in(0,\infty):\zeta_0(\{t\}\times\Delta)>0\}\quad\text{a.\,s.}\]

\subsection{The general setting}
\label{Pathw:sec:gen-setting}
Let $(r_0,v_0)$ be a $\hat\Dd$-valued random variable that is independent of $\eta$. Using the map $\alpha:\hat\Dd\to\Dd$ from~\reff{Pathw:eq:def-alpha}, we define a $\Dd$-valued random variable by $\rho_0=\alpha(r_0,v_0)$.
Then $\rho_0$ can also be considered as an arbitrary $\Dd$-valued random variable that is independent of $\eta$. In this way, we unify the settings of Sections \ref{Pathw:sec:ld-sampl-nd} and \ref{Pathw:sec:ld-sampl-dust}.

Let $(\hat Z,\rho)$ be the extended lookdown space associated with $\eta$ and $(r_0,v_0)$ as in Section \ref{Pathw:sec:ld-dec}. Recall that $(\hat Z,\rho)$ contains the lookdown space $(Z,\rho)$ associated with $\eta$ and $\rho_0$ as a subspace. We endow $\hat Z\times\R_+$ with the product metric $d^{\hat Z\times\R_+}$ defined in the beginning of Section \ref{Pathw:sec:ld-sampl-dust}. We define the $\Dd$-valued process $(\rho_t,t\in\R_+)$ of the genealogical distances between the individuals at fixed times as in equation~\reff{Pathw:eq:def-rhot}.
Then for each $n\in\N$, the process $(\gamma_n(\rho_t),t\in\R_+)$ of the restrictions to the first $n$ levels jumps only at the reproduction events that are encoded by a partition in $\p^n$ as in the beginning of Section \ref{Pathw:sec:ld-space}. On the event of probability $1$ on which condition \reff{Pathw:eq:ass-eta} holds, $(\gamma_n(\rho_t),t\in\R_+)$ thus jumps only finitely often in bounded time intervals. Between these jumps, the genealogical distances $\rho_t(i,j)$ with $i\neq j$, $i,j\in[n]$ grow linearly with slope $2$. Using also the definition of the metric $\rho$ from the beginning of Section \ref{Pathw:sec:ld-space} (there in particular that the map $t\mapsto A_s(t,i)$ is càdlàg), we deduce that the process $(\rho_t,t\in\R_+)$ has a.\,s.\ càdlàg paths. We denote the left limits by $\rho_{t-}$. The process $(\rho_t,t\in\R_+)$ is also Markovian by \cite{Sampl}*{Proposition 5.4}.

We define the $\hat\Dd$-valued process $((r_t,v_t),t\in\R_+)$ of the decomposed genealogical distances between the individuals at fixed times as in Section~\ref{Pathw:sec:ld-dec} (including Remark~\ref{Pathw:rem:rho-r-v}). In case $\Xi\in\dust$, the following condition is a.\,s.\ satisfied (as in \cite{Sampl}, where this is condition (6.2) which is checked in Section 6.2 therein):
\begin{equation}
\label{Pathw:eq:ass-eta-hat}
\eta((0,t]\times\hat\p^n)<\infty\quad\text{ for all }t\in(0,\infty)\text{ and }n\in\N.
\end{equation}
Here \label{Pathw:not:hatpn}
$\hat\p^n=\{\pi\in\p:\{\{1\},\ldots,\{n\}\}\not\subset\pi\}$
is the set of those partitions of $\N$ in which the first $n$ integers are not all in singleton blocks.
By definition of the population model in the beginning of Section \ref{Pathw:sec:ld-space}, condition \reff{Pathw:eq:ass-eta-hat} implies that the particles on the first $n\in\N$ levels reproduce at only finitely many times in bounded time intervals.
Then for $i\in[n]$, the maps $t\mapsto(z(t,i),v_t(i))$ are càdlàg with jumps only at such reproduction times. This follows from the definition of $v_t(i)$ and $z(t,i)$ in Section \ref{Pathw:sec:ld-dec}. Between such jumps, the parent $z(t,i)$ remains constant and the quantity $v_t(i)$ grows linearly with slope $1$. As a consequence, the process $((r_t,v_t),t\in\R_+)$ is a.\,s.\ càdlàg if $\Xi\in\dust$. Also recall from Remark~\ref{Pathw:rem:v-nd} that $\Xi\in\nd$ implies $v_t=0$ and $\rho_t=r_t$ for all $t\in(0,\infty)$ a.\,s. Hence, $((r_t,v_t),t\in(0,\infty))$ is a.\,s.\ càdlàg for all finite measures $\Xi$ on $\Delta$. We denote the left limits by $(r_{t-},v_{t-})$. Note that the left limits $\rho_{t-}$ and $(r_{t-},v_{t-})$ can be defined like $\rho_t$ and $(r_t,v_t)$, respectively, by ignoring a possible reproduction event at time $t$. By \cite{Sampl}*{Proposition 6.3}, the process $((r_t,v_t),t\in\R_+)$ is also Markovian.

\section{Preservation of exchangeability}
\label{Pathw:sec:exch}
In this section, we extend the exchangeability results on the genealogical distances in the lookdown model from Sections 5.3 and 11.1 of \cite{Sampl}: We consider invariance under permutations that leave the first $b\in\N$ levels unchanged, and we show exchangeability properties of the decomposed genealogical distances at various stopping times. In Subsection~\ref{Pathw:sec:exch-repr}, we show that exchangeability properties are preserved in single reproduction events. In Subsection~\ref{Pathw:sec:exch-ld}, we concatenate these reproduction events to show exchangeability properties in the lookdown model. Exchangeability in the lookdown model is also studied in \cites{DK96,DK99,BBMST09,Lab12,Kurtz98}.

Let us first repeat from \cite{Sampl}*{Sections 5.1 and 6.1} the effect of a reproduction event on the (decomposed) genealogical distances.
For $n\in\N$, $\pi\in\p_n$ and $i\in\sn$, let $\pi(i)$\label{Pathw:not:pii} be the integer $k$ such that $i$ is in the $k$-th block of $\pi$ when blocks are ordered increasingly according to their respective smallest element.
With each element $\pi$ of $\p_n$, we associate a transformation $\Dd_n\to\Dd_n$, which we also denote by $\pi$, by
\[\label{Pathw:not:pntransf}\pi(\rho)=(\rho(\pi(i),\pi(j)))_{i,j\in\N}.\]
By comparison with the construction in the beginning of Section~\ref{Pathw:sec:ld-space}, we see that a reproduction event that is encoded by a point $(t,\pi)$ of $\eta$ results in a jump of the genealogical distances that is given by
\begin{equation}
\label{eq:jump-rho}
\gamma_n(\rho_t)=\gamma_n(\pi)(\gamma_n(\rho_{t-})),
\end{equation}
which is equation (5.3) in \cite{Sampl}. Recall the set $\p^n$ from equation~\reff{Pathw:eq:def-pn}. Clearly, $\gamma_n(\pi)(\rho)=\rho$ for each $\pi\in\p\setminus\p^n$ and $\rho\in\Dd_n$.

To account for the decomposed genealogical distances, we use the set $\S_n$\label{Pathw:not:Sntransf} of semi-partitions of $\sn$. A semi-partition $\sigma$ of $\sn$ is a system of nonempty disjoint subsets of $\sn$, which we call blocks. The union $\cup\sigma$ of the blocks needs not comprise all elements of $\sn$. For each semi-partition $\sigma\in\S_n$, there exists a unique partition $\pi\in\p_n$ that has the same non-singleton blocks as $\sigma$, that is, 
$\{B\in\pi:\#B\geq 2\}=\{B\in\sigma:\#B\geq 2\}$; we define $\sigma(i)=\pi(i)$ for $i\in\sn$.
With each element $\sigma$ of $\S_n$, we associate a transformation $\hat\Dd_n\to\hat\Dd_n$, which we also denote by $\sigma$, by
$\sigma(r,v)=(r',v')$,
where
\[v'(i)=v(\sigma(i))\I{i\notin\cup\sigma}\]
and
\[r'(i,j)=\left(v(\sigma(i))\I{i\in\cup\sigma}
+r(\sigma(i),\sigma(j))+v(\sigma(j))\I{j\in\cup\sigma}\right)\I{i\neq j}\]
for $i,j\in\sn$. Furthermore, we define the map
\begin{equation}
\label{Pathw:eq:varsigma}
\varsigma_n:\p\to\S_n, \quad \pi\mapsto
\{B\cap\sn: B\in\pi, \#B\geq 2\}\setminus\{\emptyset\}
\end{equation}
which removes all singleton blocks from $\pi$ and restricts the obtained semi-partition to $\sn$.
In particular, $\hat\p^n=\{\pi\in\p:\varsigma_n(\pi)\neq\emptyset\}$ for all $n\in\N$.
By construction, each point $(t,\pi)$ of $\eta$ results in a jump of the decomposed genealogical distances at time $t$ given by
\begin{equation}
\label{eq:jump-rv}
\gamma_n(r_t,v_t)=\varsigma_n(\pi)(\gamma_n(r_{t-},v_{t-})),
\end{equation}
this is equation (6.1) in \cite{Sampl}.
For each point $\pi\in\p\setminus\hat\p^n$ and $(r,v)\in\hat\Dd_n$, we have $\varsigma_n(\pi)(r,v)=(r,v)$.

For $n\in\N$ and $b\in\sn\cup\{0\}$, we denote by $\S_{n,b}$ the set of semi-partitions $\sigma\in\S_n$ that satisfy $\sigma(i)=i$ for all $i\in\sb$. These are the $\sigma\in\S_n$ such that no element of $\sb$ is in a non-singleton block, that is, for each $i\in\sb$, either $i\notin\cup\sigma$ or $\{i\}\in\sigma$.
Therefore,
\begin{equation}
\label{Pathw:eq:p-snb}
\S_{n,b}=\{\varsigma_n(\pi):\pi\in\p\setminus\p^b\}.
\end{equation}
Similarly, we define $\p_{n,b}$ as the set of partitions of $\sn$ such that none of the first $b$ integers are in non-singleton blocks. This is the set of partitions $\pi\in\p_n$ with $\pi(i)=i$ for all $i\in\sb$. We have
\begin{equation}
\label{Pathw:eq:p-pnb}
\p_{n,b}=\{\gamma_n(\pi):\pi\in\p\setminus\p^b\}.
\end{equation}

Let $S_n$ be the group of permutations of $\sn$. We define the action of $S_n$ on $\Dd_n$ by
$p(\rho)=(\rho(p(i),p(j)))_{i,j\in\sn}$
for $p\in S_n$ and $\rho\in\Dd_n$. Analogously, we set
$p(r,v)=((r(p(i),p(j)))_{i,j\in\sn},(v(p(i)))_{i\in\sn})$ for $(r,v)\in\hat\Dd_n$, and
$p(\sigma)=\{p(B):B\in\sigma\}$ for $\sigma\in\S_n\supset\p_n$.
For $n\in\N$ and $b\in\sn\cup\{0\}$, we define the group
\begin{equation}
\label{eq:Snb}
S_{n,b}=\{p\in S_n:\text{ $p(i)=i$ for all $i\in\sb$}\}
\end{equation}
of permutations of $\sn$ that leave the first $b$ integers unchanged.
We say that a random (marked) distance matrix in $\Dd_n$ or $\hat\Dd_n$, or a random (semi-)partition of $\sn$ is $(n,b)$-exchangeable if its distribution is invariant under the action of $S_{n,b}$. For $n\in\N$, the usual exchangeability is recovered as $(n,0)$-exchangeability. For $b\in\N$, no restriction is meant by $(b,b)$-exchangeability.

\subsection{Single reproduction events}
\label{Pathw:sec:exch-repr}
In Lemma~\ref{Pathw:lem:exch-pi} below, we show that $(n,b)$-exchangeability of distance matrices is preserved under the transformations associated with independent $(n,b)$-exchangeable partitions which will later encode reproduction events in the lookdown model. This extends Lemma 11.1 of \cite{Sampl} to $(n,b)$-exchangeability. In the subsequent Lemma~\ref{Pathw:lem:exch-sigma}, we consider marked distance matrices and semi-partitions.
\begin{lem}
\label{Pathw:lem:exch-pi}
Let $n\in\N$ and $b\in\sn\cup\{0\}$. Let $\tilde\rho$ be an $(n,b)$-exchangeable random variable with values in $\Dd_n$ and let $\pi$ be an independent $(n,b)$-exchangeable random variable with values in $\p_{n,b}$. Then $\pi(\tilde\rho)$ is $(n,b)$-exchangeable.
\end{lem}
\begin{proof}
Let $p\in S_{n,b}$. As in the proof of Lemma 11.1 in \cite{Sampl}, see equations (6.2) and (6.3) therein, there exists a map $f:\p_n\to S_n$ that satisfies
\begin{equation}
\label{Pathw:eq:exch-pi-i}
\pi'(i)=f(\pi')(p(\pi')(p(i)))
\end{equation}
for all $\pi'\in\p_n$ and $i\in\sn$, and
\begin{equation}
\label{Pathw:eq:exch-pi-rho}
\pi'(\rho')=p(p(\pi')(f(\pi')(\rho')))
\end{equation}
for all $\pi'\in\p_n$ and $\rho'\in\Dd_n$.

For each $\pi'\in\p_{n,b}$, the definition of $p$ implies $p(\pi')\in\p_{n,b}$ and
$\pi'(i)=i=p(\pi')(p(i))$ for all $i\in\sb$. Hence, $f(\pi')\in S_{n,b}$ for each $\pi'\in\p_{n,b}$.

This allows to conclude analogously to the proof of Lemma 11.1 in \cite{Sampl}.
\end{proof}
\begin{lem}
\label{Pathw:lem:exch-sigma}
Let $n\in\N$ and $b\in\sn\cup\{0\}$. Let $(\tilde r,\tilde v)$ be a $(n,b)$-exchangeable random variable with values in $\hat\Dd_n$ and let $\sigma$ be an independent $(n,b)$-exchangeable random variable with values in $\S_{n,b}$. Then $\sigma(\tilde r,\tilde v)$ is $(n,b)$-exchangeable.
\end{lem}
\begin{proof}
Recall from~\reff{Pathw:eq:def-alpha} the map $\alpha:\hat\Dd_n\to\Dd_n$.
Let $(r',v')\in\hat\Dd_n$, $\rho'=\alpha(r',v')$, $\sigma'\in\S_n$, and let $\pi'$ be the partition in $\p_n$ with the same non-singleton blocks as $\sigma'$.
Then,
\begin{equation}
\label{eq:rho-rv}
\alpha(\sigma'(r',v'))=\pi'(\rho')
\end{equation}
by definition of the transformations on $\Dd_n$ and $\hat\Dd_n$ associated with each element of $\p_n$ and $\S_n$, respectively.
Writing $\sigma'(r',v')=(r'',v'')$, we obtain from equation \reff{eq:rho-rv} and the definition of the map $\alpha$ that
\begin{equation}
\label{Pathw:eq:r-exch-sigma}
r''=(((\pi'(\rho'))_{i,j}-v''(i)-v''(j))\I{i\neq j})_{i,j\in\sn}.
\end{equation}

Let $p\in S_{n,b}$ and the map $f:\p_n\to S_n$ be defined as in the proof of Lemma~\ref{Pathw:lem:exch-pi}. For $i\in\sn$, it holds
$i\in\cup\sigma'$ if and only if $p(i)\in\cup p(\sigma')$.
From equation~\reff{Pathw:eq:exch-pi-i}, we obtain
\[v''(i)=v'(\pi'(i))\I{i\notin\cup\sigma'}
=v'(f(\pi')(p(\pi')(p(i))))\I{p(i)\notin\cup p(\sigma')}.\]
Using also equations~\reff{Pathw:eq:exch-pi-rho} and \reff{Pathw:eq:r-exch-sigma}, we deduce
\begin{equation*}
\label{Pathw:eq:exch-sigma-R}
\sigma'(r',v')=p\left(p(\sigma')\left(f(\pi')(r',v')\right)\right).
\end{equation*}
We conclude analogously to the proof of Lemma 11.1 in \cite{Sampl}.
\end{proof}
Recall the partitions of the form $K_{i,j}$ which contain only $\{i,j\}$ as a non-singleton block. These partitions encode binary reproduction events.
In the next lemma, we consider the exchangeability after a transformation associated with a partition of the form $K_{i,j}$ is applied to an $(n,b)$-exchangeable distance matrix.
\begin{lem}
\label{Pathw:lem:exch-binary}
Let $n\in\N$, $b\in[n-1]\cup\{0\}$, $i,j\in[ b+1]$ with $i<j$, and $\pi=K_{i,j}$. Let $\tilde\rho$ be an $(n,b)$-exchangeable random variable with values in $\Dd_n$. Then $\pi(\tilde\rho)$ is $(n,b+1)$-exchangeable.
\end{lem}
\begin{proof}
For all $k\in\sn$,
\[\pi(k)=\left\{
\begin{aligned}
&k\quad\text{if }k<j\\
&i\quad\text{if }k=j\\
&k-1\quad\text{if }k>j.
\end{aligned}\right.\]
Let $p\in S_{n,b+1}$. Then,
\[\pi(p^{-1}(k))=\left\{
\begin{aligned}
&\pi(k)\quad\text{if }k\leq b+1\\
&p^{-1}(k)-1\quad\text{if }k>b+1.
\end{aligned}\right.\]
Let $p'$ be the permutation in $S_{n,b}$ such that $p'(k-1)=p^{-1}(k)-1$ for all $k\in\sn$ with $k>b+1$. The permutation $p'$ indeed exists and is unique as $p'\in S_{n,b}$ implies $p'(k-1)=k-1$ for all $k\in[b+1]$ with $k\geq 2$, as $p\in S_{n,b+1}$ implies $b<p^{-1}(k)-1\leq n$ for all $k\in\sn$ with $k>b+1$, and and only one possibility remains for $p'(n)$ as $p^{-1}$ is injective. 

It follows
\[\pi(p^{-1}(k))=p'(\pi(k))\]
for all $k\in\sn$. To see this, we use that $\pi(k)\leq b$ for all $k\in[b+1]$, and that $\pi(k)=k-1$ for all $k\in\sn$ with $k>b+1$. From
$\pi(k)=p'(\pi(p(k)))$
for all $k\in\sn$, and by definition of the transformations on $\Dd_n$ associated with the elements of $\p_n$,
it follows $\pi(\tilde\rho)=p(\pi(p'(\tilde\rho)))$. The assertion follows as $p'(\tilde\rho)$ and $\tilde\rho$ are equal in distribution.
\end{proof}

\subsection{In the lookdown model}
\label{Pathw:sec:exch-ld}
We formulate most results in this subsection for the process of marked distance matrices $((r_t,v_t),t\in\R_+)$. To apply these results to the distance matrices $(\rho_t,t\in\R_+)$, note that the construction in Section~\ref{Pathw:sec:ld-dec} implies that $\rho_t=\alpha(r_t,v_t)$ for all $t\in\R_+$, with $\alpha$ defined in~\reff{Pathw:eq:def-alpha}. Also recall from Remark~\ref{Pathw:rem:v-nd} that $\Xi\in\nd$ implies that $v_t=0$ and $\rho_t=r_t$ for all $t\in(0,\infty)$ a.\,s.

For $n\in\N$ and $b\in\sn\cup\{0\}$, we say a random variable is conditionally $(n,b)$-exchangeable (given a sigma-algebra or a random variable) if its conditional distribution is a.\,s.\ invariant under the action of the group $S_{n,b}$.

The following lemma generalizes Proposition 5.8 in \cite{Sampl} and shows that conditioned on the event that until time $t$ no reproduction events affect the genealogical distances between the first $b$ individuals, the marked distance matrix $\gamma_n(r_t,v_t)$ is $(n,b)$-exchangeable if this holds for $\gamma_n(r_0,v_0)$. This assertion also holds conditionally given the point measure $\zeta_0$.
\begin{lem}
\label{Pathw:lem:exch-Rt}
Let $n\in\N$, $b\in[n]\cup\{0\}$, and $t\in\R_+$. Assume that $\gamma_n(r_0,v_0)$ is $(n,b)$-exchangeable. Then conditionally given $\zeta_0$, the marked distance matrix
\[\I{\eta((0,t]\times\p^b)=0}\gamma_n(r_t,v_t)\]
is $(n,b)$-exchangeable. 
\end{lem}
The proof is analogous to \cite{Sampl}*{Proposition 5.8}. It relies on representations (as in equation~\reff{Pathw:eq:exch-Rt-repr-nd}) of the (decomposed) genealogical distances in terms of reproduction events and the growth of the (decomposed) genealogical distances. We use that these constituents
preserve exchangeability and that reproduction events that affect the (decomposed) genealogical distances between the first $n$ levels do a.\,s.\ not accumulate. The latter property is ensured in case $\Xi\in\nd$ by condition~\reff{Pathw:eq:ass-eta}, and in case $\Xi\in\dust$ by condition~\reff{Pathw:eq:ass-eta-hat}.
To account for the growth of the (decomposed) genealogical distances according to the description in Section \ref{Pathw:sec:gen-setting},
we define for $t\in\R_+$ and $n\in\N$ the maps
\[\lambda_t:\Dd_n\to\Dd_n,\quad \rho\mapsto\rho+\underline{\underline{2}}_n t\]
and
\[\hat\lambda_t:\hat\Dd_n\to\hat\Dd_n,\quad (r,v)\mapsto(r,v+\underline{1}_n t),\]
where we write $\underline{1}_n=(1)_{i\in\sn}$ and $\underline{\underline{2}}_n=2(\I{i\neq j})_{i,j\in\sn}$.
The jumps of the (decomposed) genealogical distances are described in equations \reff{eq:jump-rho} and \reff{eq:jump-rv}.
Recall also the restriction $\gamma_n:\p\to\p_n$ and the map $\varsigma_n$ defined in~\reff{Pathw:eq:varsigma}.
\begin{proof}[Proof of Lemma \ref{Pathw:lem:exch-Rt}]
In this proof, we always condition on the event $\{\eta((0,t]\times\p^b)=0\}$. This does not affect the distribution of the Poisson random measure $\eta(\cdot\cap((0,\infty)\times(\p\setminus\p^b)))$.

On an event of probability $1$, let $(t_1,\pi_1)$, $(t_2,\pi_2),\ldots$ with $0<t_1<t_2<\ldots$ be the points of $\eta$ in $(0,\infty)\times(\p^n\setminus\p^b)$. Let $L=\eta((0,t]\times(\p^n\setminus\p^b))$. Conditionally given $\zeta_0$ and $(t_1,t_2,\ldots)$, the random partitions $\pi_1,\pi_2,\ldots$ are independent and $\gamma_n(\pi_k)$ is $(n,b)$-exchangeable for each $k\in\N$. Moreover, equation~\reff{Pathw:eq:p-pnb} yields $\gamma_n(\pi_k)\in\p_{n,b}$ for all $k\in[ L]$ a.\,s.
By accounting for the growth of the genealogical distances and their jumps, and as we condition on $\{\eta((0,t]\times\p^b)=0\}$, we obtain
\begin{equation}
\label{Pathw:eq:exch-Rt-repr-nd}
\gamma_n(\rho_t)=\lambda_{t-t_L}\circ\gamma_n(\pi_L)\circ\lambda_{t_L-t_{L-1}}\circ\ldots\circ\gamma_n(\pi_1)\circ\lambda_{t_1}(\gamma_n(\rho_0))\quad\text{a.\,s.}
\end{equation}
on $\{L\geq 1\}$, and $\gamma_n(\rho_t)=\lambda_t(\gamma_n(\rho_0))$ a.\,s.\ on $\{L=0\}$. Lemma~\ref{Pathw:lem:exch-pi} implies that the distance matrix $\gamma_n(\rho_t)$ is $(n,b)$-exchangeable conditionally given $\zeta_0$. In case $\Xi\in\nd$, this also holds for the marked distance matrix $(r_t,v_t)$ as $v_t=0$ and $r_t=\rho_t$ a.\,s.\ by Remark~\ref{Pathw:rem:v-nd}.

An analogous argument applies in case $\Xi\in\dust$. As \reff{Pathw:eq:ass-eta-hat} holds a.\,s.\ in this case, we can now define $(t_1,\pi_1),(t_2,\pi_2),\ldots$ with $0<t_1<t_2<\ldots$ to be the points of $\eta$ in $(0,\infty)\times(\hat\p^n\setminus\p^b)$ a.\,s. Now let $L=\eta((0,t]\times(\hat\p^n\setminus\p^b)$. Conditionally given $\zeta_0$ and $(t_1,t_2,\ldots)$, the random partitions $\pi_1,\pi_2,\ldots$ are independent and the random semi-partition $\varsigma_n(\pi_k)$ is $(n,b)$-exchangeable for all $k\in\N$. Moreover, equation~\reff{Pathw:eq:p-snb} yields $\varsigma_n(\pi_k)\in\S_{n,b}$ for all $k\in[ L]$ a.\,s.
By accounting for the growth of the decomposed genealogical distances and their jumps, we obtain
\[\gamma_n(r_t,v_t)=\hat\lambda_{t-t_L}\circ\varsigma_n(\pi_L)\circ\hat\lambda_{t_L-t_{L-1}}\circ\ldots\circ\varsigma_n(\pi_1)\circ\hat\lambda_{t_1}(\gamma_n(r_0,v_0))\quad\text{a.\,s.}\]
on $\{L\geq 1\}$, and $\gamma_n(r_t,v_t)=\hat\lambda_t(\gamma_n(r_0,v_0))$ a.\,s.\ on $\{L=0\}$. The assertion follows from Lemma~\ref{Pathw:lem:exch-sigma}.
\end{proof}
In the next corollary, we set $(r_{0-},v_{0-})=(r_0,v_0)$.
\begin{cor}
\label{Pathw:cor:exch-R}
Assume that $(r_0,v_0)$ is exchangeable. Let $\tau$ be a $\zeta_0$-measurable and a.\,s.\ finite random time. Then $(r_\tau,v_\tau)$ and $(r_{\tau-},v_{\tau-})$ are exchangeable.
\end{cor}
\begin{proof}
For $k\in\N$, we define a $\zeta_0$-measurable random time $\tau^k$ that assumes countably many values by $\tau^k=(j+1)/k$ on the event $\{\tau\in[j/k,(j+1)/k)\}$ for $j\in\N_0$. For $n\in\N$, $p\in S_n$, and bounded continuous $\phi$, Lemma~\ref{Pathw:lem:exch-Rt} with $b=0$ yields
\begin{align*}
\E[\phi(\gamma_n(r_{\tau^k},v_{\tau^k}))]
& =\sum_{j\in\N_0}\E[\phi(\gamma_n(r_{j/k},v_{j/k}));\tau^k=j/k]\\
& =\sum_{j\in\N_0}\E[\phi(p(\gamma_n(r_{j/k},v_{j/k})));\tau^k=j/k]
=\E[\phi(p(\gamma_n(r_{\tau^k},v_{\tau^k})))].
\end{align*}
We let $k$ tend to infinity. The assertion follows as $t\mapsto \gamma_n(r_t,v_t)$ is càdlàg a.\,s. To prove the assertion for $(r_{\tau-},v_{\tau-})$, we replace $\tau^k$ with $\tilde\tau^k=\lfloor \tau k\rfloor/k$.
\end{proof}

In the next lemma, we see that at the time of the first reproduction event that affects the genealogical distances between the first $b$ individuals, conditioned on this reproduction event being binary,
the matrix of the genealogical distances between the first $n$ individuals is $(n,b+1)$-exchangeable if it is $(n,b)$-exchangeable at time zero.
\begin{lem}
\label{Pathw:lem:exch-b}
Let $n\in\N$ and $b\in[n-1]$ with $n,b\geq 2$. Assume $\Xi\{0\}>0$ and that $\gamma_n(\rho_0)$ is $(n,b)$-exchangeable. Let
\[\tau=\inf\{t>0:\eta((0,t]\times\p^b)>0\}.\]
Then, $\tau<\infty$ a.\,s.\ and the distance matrix
\[\I{\zeta_0(\{\tau\}\times\Delta)=0}\gamma_n(\rho_\tau)\]
is $(n,b+1)$-exchangeable.
\end{lem}
\begin{proof}
The assumptions $b\geq 2$ and $\Xi\{0\}>0$ imply $\tau<\infty$ a.\,s.\ and $\P(\zeta_0(\{\tau\}\times\Delta)=0)>0$. In this proof, we always condition on the event $\{\zeta_0(\{\tau\}\times\Delta)=0\}$.

On an event of probability $1$, let $(t_1,\pi_1),(t_2,\pi_2)\ldots$ with $0<t_1<t_2<\ldots$ be the points of $\eta$ in $(0,\infty)\times\p^n$, and let $L=\eta((0,\tau]\times\p^n)$. Then, $\tau=t_L$ a.\,s. The partitions $\gamma_n(\pi_1),\gamma_n(\pi_2),\ldots$ are conditionally independent given $\zeta_0$, $L$, and $(t_1,\ldots,t_L)$. Conditionally given $\zeta_0$, $L$, and $(t_1,\ldots,t_L)$, the partitions $\gamma_n(\pi_k)$ for $k\in[ L-1]$ are also $(n,b)$-exchangeable and by equation~\reff{Pathw:eq:p-pnb} in $\p_{n,b}$. On the event $\{\zeta_0(\{\tau\}\times\Delta)=0\}$, conditionally given $\zeta_0$, $L$ and $(t_1,\ldots,t_L)$, the partition $\gamma_n(\pi_L)$ a.\,s.\ contains one block with two elements in $\sb$ and apart from that only singleton blocks.
This follows from the construction of $\eta$ in Section \ref{Pathw:sec:two-step-eta}. By accounting for the growth and the jumps of the genealogical distances, we obtain
\[\gamma_n(\rho_\tau)=\gamma_n(\pi_{L})\circ\lambda_{t_{L}-t_{L-1}}\circ\ldots\circ\gamma_n(\pi_1)\circ\lambda_{t_1}(\gamma_n(\rho))\quad\text{a.\,s.}\]
The assertion follows from Lemmas~\ref{Pathw:lem:exch-pi} and \ref{Pathw:lem:exch-binary}.
\end{proof}

We also consider exchangeability properties at certain stopping times.
For $n\in\N$, we denote by $S^n_\infty$ the group of bijections $p:\N\to\N$ with $p(i)=i$ for all $i>n$.
We denote the space of probability measures on $\hat\Dd$ by $\Mm_1(\hat\Dd)$, and we define the measurable function
\[\beta_n:\hat\Dd\to\Mm_1(\hat\Dd),\quad (r,v)\mapsto\frac{1}{n!}\sum_{p\in S^n_\infty}\delta_{p(r,v)}.\]
Intuitively, the image of $(r,v)$ under $\beta_n$ corresponds to the orbit of $(r,v)$ under the action of the group $S^n_\infty$.

For $t\in\R_+$, we denote by $\F^n_t$ the sigma-algebra generated by $(\beta_n(r_s,v_s),s\in[0,t])$.
A filtration is defined by $\F^n=(\F^n_t,t\in\R_+)$.
We will use the following lemma in Section~\ref{Pathw:sec:part} to study asymptotic frequencies that are invariant under permutation of the first $n$ levels.
\begin{lem}
\label{Pathw:lem:exch-J}
Let $n\in\N$ and assume that $(r_0,v_0)$ is exchangeable. Let $\tau$ be a finite $\F^n$-stopping time. Then $\gamma_n(r_\tau,v_\tau)$
is exchangeable.
\end{lem}
\begin{proof}
We show that for each $t\in\R_+$, the marked distance matrix $\gamma_n(r_t,v_t)$ is exchangeable conditionally given $\F^n_t$.
The assertion then follows for stopping times that assume countably many values, and by an approximation argument as in the proof of Corollary~\ref{Pathw:cor:exch-R} also for all finite stopping times.

Let $K$ be the probability kernel from $\Mm_1(\hat\Dd)$ to $\hat\Dd$ given by $K(\xi,\cdot)=\xi$ for $\xi\in\Mm_1(\hat\Dd)$. By Lemma~\ref{Pathw:lem:exch-Rt}, the marked distance matrix $(r_t,v_t)$ is exchangeable, hence
\begin{align*}
&\P((r_t,v_t)\in B,\beta_n(r_t,v_t)\in B')\\
&=\P(p(r_t,v_t)\in B,\beta_n( p(r_t,v_t)) \in B')
=\P(p(r_t,v_t)\in B,\beta_n (r_t,v_t)\in B')
\end{align*}
for all measurable $B, B'$ and all $p\in S^n_{\infty}$. This implies that $K$ is a regular conditional distribution of $(r_0,v_0)$ given $\beta_n(r_0,v_0)$, and of $(r_t,v_t)$ given $\beta_n(r_t,v_t)$. Now we apply Theorem 2 of \cite{RP81} to the Markov process $((r_t,v_t),t\in\R_+)$, the measurable function $\beta_n$, and the probability kernel $K$ to obtain that for each $t\in\R_+$, the marked distance matrix $(r_t,v_t)$ has the same conditional distribution given $\beta_n(r_t,v_t)$ as given $\F^n_t$.
This follows from equation (1) in \cite{RP81} and implies the assertion.
\end{proof}

Furthermore, recall the groups of permutations $S_{n,b}$ for $n\in\N$ and $b\in\sn\cup\{0\}$ from \reff{eq:Snb}.
We write $\Mm$ for the space of probability measures on $\hat\Dd_n$, and we define the function
\[\beta_{n,b}:\hat\Dd_n\to\Mm,\quad
(r,v)\mapsto\frac{1}{(n-b)!}\sum_{p\in S_{n,b}}\delta_{p(r,v)}.\]
For $t\in\R_+$, we denote by $\F^{n,b}_t$ the sigma-algebra generated by $(\beta_{n,b}(r_s,v_s),s\in[0,t])$.
A filtration is defined by $\F^{n,b}=(\F^{n,b}_t,t\in\R_+)$.
We will use the following lemma in Section~\ref{Pathw:sec:unif-ld-proofs} where we consider some relative frequencies in the lookdown model that do not change under permutations that leave the first $b$ levels fixed.
\begin{lem}
\label{Pathw:lem:exch-H}
Let $n\in\N$ and $b\in\sn\cup\{0\}$. Assume that $\gamma_n(r_0,v_0)$ is $(n,b)$-exchangeable. Let $\tau$ be a finite $\F^{n,b}$-stopping time. Then the marked distance matrix
\[\I{\eta((0,\tau]\times\p^b)=0}\gamma_n(r_\tau,v_\tau)\]
is $(n,b)$-exchangeable.
\end{lem}
\begin{proof}
We show that for each $t\in\R_+$, the marked distance matrix 
\[\I{\eta((0,t]\times\p^b)=0}\gamma_n(r_t,v_t)\]
is $(n,b)$-exchangeable conditionally given $\F^{n,b}_t$.

We enlarge the spaces $\hat\Dd_n$ and $\Mm$ by a coffin state $\partial$. Let $K$ be the probability kernel from $\Mm$ to $\hat\Dd_n$ such that $K(\partial,\{\partial\})=1$, and such that $K(\xi,\cdot)=\xi$ for all $\xi\in\Mm\setminus\{\partial\}$.

Let $\tau'=\inf\{t>0:\eta((0,t]\times\p^b)>0\}$ and set $R_t=\gamma_n(r_t,v_t)$ for $t<\tau'$, and $R_t=\partial$ for $t\geq\tau'$. Analogously to the proof of Lemma~\ref{Pathw:lem:exch-J}, $K$ is a regular conditional distribution of $\gamma_n(r_0,v_0)$ given $(\beta_{n,b}(\gamma_n(r_0,v_0))$ as a consequence of the assumed exchangeability. For all $t\in\R_+$, Lemma~\ref{Pathw:lem:exch-Rt} implies that $K$ is a regular conditional distribution of $R_t$ given $\beta_{n,b}(R_t)$, where we set $\beta_{n,b}(\partial)=\partial$. We apply Theorem 2 of Rogers and Pitman \cite{RP81} to the Markov process $(R_t, t\in\R_+)$, the measurable function $\beta_{n,b}$, and the probability kernel $K$ to obtain that for each $t\in\R_+$, the random variable $R_t$ has the same conditional distribution given $\F^{n,b}_t$ as given $\beta_{n,b}(R_t)$.
This implies the assertion as in the proof of Lemma~\ref{Pathw:lem:exch-J}.
\end{proof}

\section{Uniform convergence in the lookdown model}
\label{Pathw:sec:unif-ld}
Donnelly and Kurtz \cite{DK99} prove that the measure-valued processes whose states are the uniform measures
of the types on the first $n$ levels in the lookdown model converge a.\,s.\ as $n$ tends to infinity. Lemmas~\ref{Pathw:lem:BC-bound} and \ref{Pathw:lem:BC-bound-dust} below give the existence of and uniform convergence to asymptotic frequencies of subsets of individuals that are characterized by the (marked) genealogical distances.  

For $\ell\in\N$, let
\[b_\ell:\p\to\N_0,\quad b_\ell(\pi)=\ell-\#(\gamma_\ell(\pi))\]
and
\[\hat b_\ell:\p\to\N_0,\quad \hat b_\ell(\pi)=\#\{i\in\sell:\{i\}\notin\pi\}.\]
Moreover, let
\[N^\ell(I)=\int_{I\times\p}b_\ell(\pi)\,\eta(ds\;d\pi)\]
and
\[\hat N^\ell(I)=\int_{I\times\p}\hat b_\ell(\pi)\,\eta(ds\;d\pi)\]
for each interval $I\subset(0,\infty)$. The random variable $N^\ell(I)$ is the number of newborn particles on the first $\ell$ levels in the time interval $I$. The random variable $\hat N^\ell(I)$ counts in the reproducing particles in each reproduction event, in particular, it also takes account of reproduction events in which only the reproducing particles occupy levels in $\sell$.

\begin{lem}
\label{Pathw:lem:BC-bound}
Assume that $\rho_0$ is exchangeable. Let $b\in\N$ and let $f$ be a measurable function from $\R_+\times\R_+^b$ to $\{0,1\}$. For $t\in\R_+$ and $i,n\in\N$, define
\[Y^b_i(t)=f(t,(\rho_t(j,b+1+i))_{j\in\sb})\]
and
\[X_n(t)=\frac{1}{n}\sum_{i=1}^n Y^b_i(t).\]
Impose the assumption on $f$ that
\begin{equation}
\label{eq:cond-newb}
\left|X_n(t)-X_n(s)\right|\I{\eta((s,t]\times\p^b)=0}\leq \frac{1}{n}N^{b+1+n}(s,t]
\end{equation}
for all $n\in\N$ and $0\leq s< t$. Then there exists a process $(X(t),t\in\R_+)$ with $\lim_{n\to\infty}\sup_{t\in[0,T]}\left|X_n(t)-X(t)\right|= 0$ a.\,s.\ for all $T\in\R_+$.
\end{lem}
\begin{lem}
\label{Pathw:lem:BC-bound-dust}
Assume $\Xi\in\dust$ and that $(r_0,v_0)$ is exchangeable. Let $b\in\N$ and let $f$ be a measurable function from $\R_+\times\R_+^{2b}\times\R_+$ to $\{0,1\}$. For $t\in\R_+$ and $i,n\in\N$, define
\[Y^b_i(t)=f(t,(r_t(j,b+i),v_t(j))_{j\in\sb},v_t(b+i))\]
and
\[X_n(t)=\frac{1}{n}\sum_{i=1}^n Y^b_i(t).\]
Impose the assumption on $f$ that
\begin{equation}
\label{eq:cond-repr}
\left|X_n(t)-X_n(s)\right|\I{\eta((s,t]\times\hat\p^b)=0}\leq\frac{1}{n}\hat N^{b+n}(s,t]
\end{equation}
for all $n\in\N$ and $0\leq s< t$. Then there exists a process $(X(t),t\in\R_+)$ with $\lim_{n\to\infty}\sup_{t\in[0,T]}\left|X_n(t)-X(t)\right|= 0$ a.\,s.\ for all $T\in\R_+$.
\end{lem}
We defer the proofs of these lemmas to Section~\ref{Pathw:sec:unif-ld-proofs}.

\section{Two families of partitions}
\label{Pathw:sec:part}
From the process $(\rho_t,t\in\R_+)$, we now read off the flow of partitions $(\Pi_{s,t},0\leq s\leq t)$. This process corresponds to the dual flow of partitions in Foucart \cite{Foucart12} and to the flow of partitions in Labbé \cite{Lab12}.
We define the random partition $\Pi_{s,t}$\label{Pathw:not:Pist} of $\N$ by
\[\text{$i$ and $j$ are in the same block of $\Pi_{s,t}$}\quad\Leftrightarrow\quad \rho_t(i,j)<2(t-s)\]
for all $i,j\in\N$ with $i\neq j$. That is, $i$ and $j$ are in the same block of $\Pi_{s,t}$ if and only if $A_s(t,j)=A_s(t,i)$
which means that the individuals $(t,i)$ and $(t,j)$ have a common ancestor at time $s$.
For each $s\in\R_+$ and $n\in\N$, the process $t\mapsto\gamma_n(\Pi_{s,s+t})$ jumps only at the times of reproduction events that are encoded by a partition in $\p^n$. These times do not accumulate on the event of probability $1$ on which condition~\reff{Pathw:eq:ass-eta} is satisfied.
For all $0\leq s\leq t$, the random partition $\Pi_{s,t}$ is exchangeable. This follows from Lemma~\ref{Pathw:lem:exch-Rt} (where we may assume w.\,l.\,o.\,g., as $\Pi_{s,t}$ is $\eta$-measurable, that $\rho_0$ is exchangeable). We will apply the flow of partitions in the dust-free case.

For application in the case with dust, we define for each $a\in\N$ and $\ep>0$ a family $(\Pi^{a,\ep}_t,t\in\R_+)$ of partitions of $\N$.
Similarly to the partition induced by $\sim^\ep$ in \cite{Sampl}*{Section 10.4}, our intention is that individuals at time $t$ whose levels are in the same block of $\Pi^{a,\ep}_t$ should have parents that are close to each other in the extended lookdown space (the definition of a parent is given Section~\ref{Pathw:sec:ld-dec}). We will define $\Pi^{a,\ep}_t$ accordingly except for at most one block.

First, we define for each $t\in\R_+$ and $I\subset\sa$ the subset of $\N$\label{Pathw:not:CaepIt}
\begin{align*}
C^{a,\ep,I}_t&=\{i\in\N:v_t(i)\geq t\}\cap\bigcap_{k\in I:\,v_t(k)\geq t}\{i\in\N:r_t(i,k)\vee\left|v_t(i)-v_t(k)\right|<\ep\}\\
&\quad\cap\bigcap_{\ell\in\sa\setminus I:\,v_t(\ell)\geq t}\{r_t(i,\ell)\vee\left|v_t(i)-v_t(\ell)\right|\geq\ep\}
\end{align*}
\begin{rem}
\label{Pathw:rem:AI}
Clearly, $(C^{a,\ep,I}_t,I\subset\sa)$ is a family of subsets whose union is
$\{i\in\N:v_t(i)\geq t\}$. Any two such subsets are either disjoint or equal. For $I\subset[a]$ and $i,j\in C^{a,\ep,I}_t$, the construction in Section~\ref{Pathw:sec:ld-dec} implies that the parents of the individuals $(t,i)$ and $(t,j)$ are also the parents of the individuals $(0,A_0(t,i))$ and $(0,A_0(t,j))$, respectively. If moreover $I\neq\emptyset$, then the genealogical distance $r_t(i,j)$ between these parents is less than $2\ep$, and $|v_t(i)-v_t(j)|<2\ep$.
\end{rem}
Now we let $i,j\in\N$ be in the same block of $\Pi^{a,\ep}_t$\label{Pathw:not:Piaept} if and only if one of the following two conditions is satisfied:
\begin{enumerate}[label=(\roman{*}),ref=(\roman{*})]
\item\label{Pathw:item:Pit-above} $v_t(i)=v_t(j)<t$ and $r_t(i,j)=0$
\item\label{Pathw:item:Pit-below} There exists $I\subset\sa$ such that $i,j\in C^{a,\ep,I}_t$.
\end{enumerate}

Condition~\ref{Pathw:item:Pit-above} means that the individuals $(t,i)$ and $(t,j)$ have the same parent in the (extended) lookdown space, and that this parent lives after time zero. That is, $(t,i)$ and $(t,j)$ have a common ancestor at time zero, and the individuals on each of the ancestral lineages of the individuals $(t,i)$ and $(t,j)$ are in singleton blocks in each reproduction event until these ancestral lineages merge. The individual in which these ancestral lineages merge is the parent of both the individuals $(t,i)$ and $(t,j)$, when we identify individuals with genealogical distance zero. In this sense, the individuals $(t,i)$ and $(t,j)$ may be called siblings.

Using the definitions of $v_t(i)$ and $r_t(i,j)$, it can be seen that for each $n\in\N$, on the event of probability $1$ on which condition~\reff{Pathw:eq:ass-eta-hat} is satisfied, the process $t\mapsto\gamma_n(\Pi^{a,\ep}_t)$ jumps only at the times of reproduction events that are encoded by a partition in $\hat\p^n$, and that these times do not accumulate.

In the next two lemmas, we show that the asymptotic frequencies in $\Pi_{s,t}$ and $\Pi^{a,\ep}_t$, respectively, exist simultaneously for uncountably many $t$ on an event of probability $1$, and that the relative frequencies converge uniformly for $t$ in compact intervals.
On an event of probability $1$, the left limits
\[\Pi_{s,t-}:=\lim_{s'\uparrow t}\Pi_{s,s'}=\{\{j\in\N:\rho_{t-}(i,j)<2(t-s)\}:i\in\N\}\]
and
\[\Pi^{a,\ep}_{t-}:=\lim_{s\uparrow t}\Pi^{a,\ep}_{s}=\{\{j\in\N:v_{t-}(i)=v_{t-}(j)<t,r_{t-}(i,j)=0\}:i\in\N\}\]
exist for all $t\in(0,\infty)$ and $s\in[0,t)$.
The partitions $\Pi_{s,t-}$ and $\Pi^{a,\ep}_{t-}$ are left limits with respect to the topology on $\p$ that is generated by the restriction maps $\gamma_n$, $n\in\N$. They can also be defined like $\Pi_{s,t}$ and $\Pi^{a,\ep}_t$, respectively, except that
a possible reproduction event at time $t$ is ignored.
In the next two lemmas, we also show regularity properties in $t$, and that taking (left) limits in $t$ commutes with taking asymptotic frequencies.
\begin{lem}
\label{Pathw:lem:part-nd}
Let $s,T\in\R_+$ and $b\in\N$. Then,
\begin{equation}
\label{Pathw:eq:part-unif}
\lim_{n\to\infty}\sup_{t\in[s,s+T]}
\big|\left|B(\Pi_{s,t},b)\right|_n-\left|B(\Pi_{s,t},b)\right|\big|=0\quad\text{a.\,s.}
\end{equation}
The paths $[s,\infty)\to[0,1]$, $t\mapsto\left|B(\Pi_{s,t},b)\right|$ are càdlàg a.\,s. Furthermore, $\lim_{\ep\downarrow 0}\left|B(\Pi_{s,t-\ep},b)\right|=\left|B(\Pi_{s,t-},b)\right|$ for all $t\in(s,\infty)$ a.\,s.
\end{lem}
\begin{lem}
\label{Pathw:lem:part-dust}
Let $T\in\R_+$, $a,k\in\N$, and $\ep>0$. Assume $\Xi\in\dust$ and that $(r_0,v_0)$ is exchangeable. Then,
\begin{equation}
\label{Pathw:eq:part-dust}
\lim_{n\to\infty}\sup_{t\in[0,T]}\big|\left|B(\Pi^{a,\ep}_t,k)\right|_n-\left|B(\Pi^{a,\ep}_t,k)\right|\big|=0\quad\text{a.\,s.}
\end{equation}
The paths $t\mapsto\left|B(\Pi^{a,\ep}_t,k)\right|$ are càdlàg a.\,s. Furthermore, $\lim_{s\uparrow t}\left|B(\Pi^{a,\ep}_s,k)\right|=\left|B(\Pi^{a,\ep}_{t-},k)\right|$ for all $t\in(0,\infty)$ a.\,s.
\end{lem}
A result similar to Lemma~\ref{Pathw:lem:part-nd} is Proposition 2.13 in Labbé \cite{Lab12} which is applied there to study relations between the lookdown model and flows of bridges.
\begin{proof}[Proof of Lemma~\ref{Pathw:lem:part-nd}]
By $\eta$-measurability of the random variables in the assertion, we can assume w.\,l.\,o.\,g.\ that $\rho_0$ is exchangeable. By time homogeneity, it suffices to consider the case $s=0$. We choose $f$ in Lemma~\ref{Pathw:lem:BC-bound} such that \[f(t,(\rho_t(j,b+1+i))_{j\in\sb})=\I{\rho_t(b,b+1+i)<2(t-s)}\]
for all $t\in\R_+$ and $i\in\N$. On the right-hand side, we have the indicator variable of the event that $b$ and $b+1+i$ are in the same block of $\Pi_{s,t}$, which can be written as $b+1+i\in B(\Pi_{s,t},b)$. Hence,
\begin{equation}
\label{eq:part-nd-X-B}
X_n(t)\leq \frac{b+1+n}{n}\left|B(\Pi_{s,t},b)\right|_{b+1+n}\leq X_n(t)+\frac{b+1}{n}
\end{equation}
for all $t\in\R_+$ and and $n\in\N$. Here $X_n(t)$ is defined as in Lemma~\ref{Pathw:lem:BC-bound} so that $nX_n(t)$ counts the individuals at time $t$ on the levels $b+1+1,\ldots,b+1+n$ that have a common ancestor with the individual $(t,b)$ after time $s$.

By the construction in the beginning of Section \ref{Pathw:sec:ld-space}, the number $nX_n(t)$ can only jump at times $t$ with $\eta(\{t\}\times\p^{b+1+n})>0$. If there are no newborn particles on the first $b$ levels at such a jump time (i.\,e.\ if the reproduction event at time $t$ is encoded by a partition in $\p^{b+1+n}\setminus\p^b$), then $|nX_n(t)-nX_n(t-)|$ is bounded from above by the number of newborn particles on the first $b+1+n$ levels. Indeed, we then have $nX_n(t)-nX_n(t-)=a-a'+a''$, where $a$ is the number of newborn particles on levels $b+1+1,\ldots,b+1+n$ that descend from an individual at time $t-$ on a level in $B(\Pi_{s,t-},b)$.
By $a'$, we denote here the number of particles at time $t-$ that are moved from a level in $B(\Pi_{s,t-},b)\cap\{b+1+1,\ldots,b+1+n\}$ to a level higher than $b+1+n$ at time $t$.
We set $a''=1$ if the particle on level $b+1$ at time $t-$ belongs to the same block of $\Pi_{s,t-}$ as $b$ and is pushed at time $t$ to a level in $b+1+1,\ldots,b+1+n$ to make way for a newborn on level $b+1$. Else we set $a''=0$.
Note that neither $a+a''$ nor $a'$ exceed the number of newborns on the first $b+1+n$ levels.

Hence, if there are no newborn particles on the first $b$ levels between time $t'$ and $t$ for some $s\leq t'\leq t$, then the numbers of levels $b+1+1,\ldots,b+1+n$ that are in the same block as $b$ of the partitions $\Pi_{s,t'}$ and $\Pi_{s,t}$ can differ by at most the number of newborn particles on the first $b+1+n$ levels between time $t'$ and $t$. This is condition \reff{eq:cond-newb}, hence by Lemma~\ref{Pathw:lem:BC-bound}, $X_n(t)$ converges a.\,s.\ uniformly in compact intervals to $X(t):=|B(\Pi_{0,t},b)|$. The estimates~\reff{eq:part-nd-X-B} yield the convergence~\reff{Pathw:eq:part-unif}.

On the event of probability $1$ on which condition~\reff{Pathw:eq:ass-eta} holds, the paths
$t\mapsto\left|B(\Pi_{0,t},b)\right|_n$
are càdlàg and $\lim_{s'\uparrow t}\left|B(\Pi_{0,s'},b)\right|_n=\left|B(\Pi_{0,t-},b)\right|_n$
for all $t\in(0,\infty)$ and $n\in\N$. This implies that the paths $t\mapsto\abs{B(\Pi_{0,t},b)}$ are càdlàg a.\,s., and that $\lim_{s'\uparrow t}\left|B(\Pi_{0,s'},b)\right|=c_t$ for some $c_t\in[0,1]$ for each $t\in(0,\infty)$ a.\,s. To show the assertion on the left limits, let $\ep>0$, and choose on an event of probability $1$ a sufficiently large integer $n_0$ such that $\left|\left|B(\Pi_{0,t},b)\right|_n-\left|B(\Pi_{0,t},b)\right|\right|<\ep$ for all $t\in[0,T]$ and $n\geq n_0$. Then, $\limsup_{s'\uparrow t}\left|B(\Pi_{0,s'},b)\right|_n\leq c_t+\ep$ and $\liminf_{s'\uparrow t}\left|B(\Pi_{0,s'},b)\right|_n\geq c_t-\ep$ for all $n\geq n_0$ and $t\in[0,T]$. It follows $c_t=\left|B(\Pi_{0,t-},b)\right|$ for all $t\in(0,\infty)$ a.\,s.
\end{proof}
\begin{proof}[Proof of Lemma~\ref{Pathw:lem:part-dust}]
We choose $b=a\vee k$ and $f$ in Lemma~\ref{Pathw:lem:BC-bound-dust} such that
\begin{align*}
&f(t,(r_t(j,b+i),v_t(j))_{j\in\sb},v_t(b+i))\\
&=\mathbf{1}(\{r_t(k,b+i)=0,v_t(b+i)< t\}\cup\bigcup_{I\subset[a]}\{k,b+i\in C^{a,\ep,I}_t\})
\end{align*}
for all $t\in\R_+$ and $i\in\N$. Here, the event that $k$ and $b+i$ are in the same block of $\Pi^{a,\ep}_t$ stands in the the indicator variable. Hence,
\begin{equation}
\label{eq:part-dust-X-B}
X_n(t)\leq\frac{b+n}{n}\left|B(\Pi^{a,\ep}_t,k)\right|_{b+n}\leq X_n(t)+\frac{b}{n}
\end{equation}
for all $n\in\N$ and $t\in\R_+$ a.\,s., with $X_n(t)$ as defined in Lemma~\ref{Pathw:lem:BC-bound-dust}.

By the construction in the beginning of Section \ref{Pathw:sec:ld-space} and in Section \ref{Pathw:sec:ld-dec}, the number $nX_n(t)$ can only jump at times $t$ with $\eta(\{t\}\times\hat\p^{b+n})>0$. If there are no particles on the first $b$ levels that are newborn or progenitor in the reproduction event at such a time $t$ (i.\,e.\ if the reproduction event at time $t$ is encoded by a partition in $\hat\p^{b+n}\setminus\hat\p^b$), then $nX_n(t)$ cannot increase and can decrease by at most the number of particles on levels $b+1,\ldots,b+n$ that are newborn or progenitor in the reproduction event at time $t$. Indeed, such a newborn or progenitor particle becomes its own parent so that the its new level cannot belong to the block $B(\Pi_{t}^{a,\ep},k)$.
As there are no newborns on the first $b$ levels, each progenitor particle on a level in $b+1,\ldots,b+n$ at time $t$ had at time $t-$ a level in $b+1,\ldots,b+n$ which may have been in $B(\Pi_{t-}^{a,\ep},k)$.
The number of newborn particles on the first $b+n$ levels equals the number of other particles that are pushed from one of these levels to a level above $b+n$.
The particles at time $t-$ on a level in $B(\Pi^{a,\ep}_{t-},k)\cap\{b+1,\ldots,b+n\}$ that are not progenitors and not pushed above level $n+b$ have a level in $B(\Pi^{a,\ep}_t,k)\cap\{b+1,\ldots,b+n\}$ at time $t$.
This yields condition \reff{eq:cond-repr}, hence by Lemma~\ref{Pathw:lem:BC-bound-dust} and the estimates~\reff{eq:part-dust-X-B}, $X_n(t)$ converges a.\,s.\ uniformly in compact intervals to $X(t):=|B(\Pi^{a,\ep}_t,k)|$, which yields the convergence~\reff{Pathw:eq:part-dust}.

On the event of probability $1$ on which condition~\reff{Pathw:eq:ass-eta-hat} holds, the processes $t\mapsto\left|B(\Pi^{a,\ep}_t,k)\right|_n$ are càdlàg for all $n\in\N$ a.\,s., hence the other assertions can be deduced from~\reff{Pathw:eq:part-dust}.
\end{proof}

To construct probability measures in the next section, we will need families of partitions with proper frequencies. A partition $\pi$ is said to have proper frequencies if $\sum_{B\in\pi}|B|=1$, that is, the asymptotic frequencies of its blocks exist and sum up to $1$. In case $\Xi\in\nd$, the partition $\Pi_{s,t}$ has proper frequencies a.\,s.\ for each $0\leq s< t$. This follows from \cite{Schw00}*{Proposition 30}. Indeed, that the $\p$-valued process $(\Pi_{t,(t-s)-},s\in[0,t))$ is a $\Xi$-coalescent up to time $t$ can be seen, for instance, by comparing the construction from the point measure $\eta$ in Section~\ref{Pathw:sec:ld-space} and the Poisson construction of Schweinsberg \cite{Schw00}*{Section 3}. The next two lemmas show that uncountably many partitions have proper frequencies on an event of probability $1$. They give the existence of lower bounds on the number of blocks whose asymptotic frequencies add up to $1-\ep$ that are uniform for $t$ in compact intervals.
\begin{lem}
\label{Pathw:lem:proper-nd}
Assume $\Xi\in\nd$. Let $s\in\R_+$, $T\in(0,\infty)$, and $\ep\in(0,T)$. Then, on an event of probability $1$, there exists an integer $k$ such that
\[\sum_{i\in M(\Pi_{s,t})\cap\sk}\left|B(\Pi_{s,t},i)\right|>1-\ep\]
for all $t\in[s+\ep,s+T]$. In particular, the partition $\Pi_{s,t}$ has proper frequencies for all $t\in(s,\infty)$ a.\,s. 
\end{lem}
\begin{lem}
\label{Pathw:lem:proper-dust}
Assume $\Xi\in\dust$ and that $(r_0,v_0)$ is exchangeable. Let $a\in\N$, $\ep,\tilde\ep>0$, and $T\in\R_+$. Then, on an event of probability $1$, there exists an integer $k$ such that
\[\sum_{i\in M(\Pi^{a,\ep}_{t})\cap\sk}|B(\Pi^{a,\ep}_{t},i)|>1-\tilde\ep\]
for all $t\in[0,T]$. In particular, the partition $\Pi^{a,\ep}_t$ has proper frequencies for all $t\in\R_+$ a.\,s.
\end{lem}
A property like the assertion of Lemma~\ref{Pathw:lem:proper-nd} is also considered in Section 6.1 of Labbé \cite{Lab12}.
\begin{proof}[Proof of Lemma~\ref{Pathw:lem:proper-nd}]
Again we assume w.\,l.\,o.\,g.\ that $\rho_0$ is exchangeable, and it suffices to consider
the case $s=0$.

\emph{Step 1. } There exists an event of probability $1$ on which for all $t\in(\ep,\infty)$ with $\#\Pi_{0,t}=\infty$, the partitions $\Pi_{0,t}$ and $\Pi_{0,t-}$ do not contain any singleton blocks.
Indeed, the partition $\Pi_{0,\ep}$ contains a.\,s.\ no singletons. For $t\in(\ep,\infty)$, an implication of $\#\Pi_{0,t}=\infty$ ($\#\Pi_{0,t-}=\infty$) is that $\#\Pi_{\ep,t}=\infty$ ($\#\Pi_{\ep,t-}=\infty$, respectively). This also implies that all individuals at time $\ep$ have a descendant at time $t$ (at time $t-$) as the trajectories of the particles in the population model do not cross, see \cite{Sampl}*{Remark 5.1}.
Hence, the assertion of Step 1 holds with the event of probability $1$ that $\Pi_{0,\ep}$ contains no singletons.

\emph{Step 2. } In this step, we assume $\Xi_0(\Delta)>0$ and we show that $\Pi_{0,t^k-}$ has proper frequencies for all $k\in\N$ a.\,s. By Corollary~\ref{Pathw:cor:exch-R}, the random partition $\Pi_{0,t^k-}$ is exchangeable. On the event that the partition $\#\Pi_{0,t^k-}$ has finitely many blocks, it has proper frequencies a.\,s. If $\P(\#\Pi_{0,t^k-}=\infty)>0$, then also $\Pi_{0,t^k-}$, conditioned on the event $\{\#\Pi_{0,t^k-}=\infty\}$, is exchangeable. The assertion of step 2 now follows from step 1 and Kingman's correspondence.

\emph{Step 3. } Recall the filtration $\F^\ell$ from Section~\ref{Pathw:sec:exch-ld}. For each $\ell\in\N$, the process
\[\left(\sum_{i\in M(\Pi_{0,t})\cap\sell}|B(\Pi_{s,t},i)|,t\in\R_+\right)\]
is adapted with respect to the usual augmentation of $\F^\ell$ and has a.\,s.\ càdlàg paths by Lemma~\ref{Pathw:lem:part-nd}. Hence,
\[\vartheta_{\ep,\ell}:=\inf\left\{t\geq \ep:
\sum_{i\in M(\Pi_{0,t})\cap\sell}|B(\Pi_{0,t},i)|<1-\ep\right\}\]
is a stopping time with respect to the usual augmentation
of $\F^n$ for all $n\in\sell$. As $\vartheta_{\ep,j}\geq\vartheta_{\ep,\ell}$ for integers $j\geq \ell$, it follows that
\[\vartheta_\ep:=\sup_{\ell\in\N}\vartheta_{\ep,\ell}=\lim_{\ell\to\infty}\vartheta_{\ep,\ell}\]
is a stopping time with respect to the usual augmentation of $\F^n$ for all $n\in\N$.

By Lemma~\ref{Pathw:lem:exch-J}, the distance matrix $\rho_{\vartheta_\ep\wedge T}$ is exchangeable, hence the partition $\Pi_{0,\vartheta_\ep\wedge T}$ is exchangeable. On the event that the partition $\Pi_{0,\vartheta_\ep\wedge T}$ has finitely many blocks, it has proper frequencies a.\,s. If $P(\#\Pi_{0,\vartheta_\ep\wedge T}=\infty)>0$, then $\Pi_{0,\vartheta_\ep\wedge T}$, conditioned on the event
$\{\#\Pi_{0,\vartheta_\ep\wedge T}=\infty\}$, remains exchangeable. It follows from step 1 and Kingman's correspondence that $\Pi_{0,\vartheta_\ep\wedge T}$ has a.\,s.\ proper frequencies.

\emph{Step 4. }
A.\,s.,
$\left|\left|B_i(\Pi_{0,t})\right|_n-\left|B_i(\Pi_{0,t-})\right|_n\right|\leq 1/n$
for all $n,i\in\N$ and all $t\in(0,\infty)\setminus\{t^k:k\in\N\}$ as only binary reproduction events can occur at these times.
Lemma~\ref{Pathw:lem:part-nd} now implies $\left|B_i(\Pi_{0,t})\right|=\left|B_i(\Pi_{0,t-})\right|$ for all $i\in\N$ and $t\in(0,\infty)\setminus\{t^k:k\in\N\}$ a.\,s. Hence, $\left|\Pi_{0,t}\right|_1=\left|\Pi_{0,t-}\right|_1$ for all $t\in(0,\infty)\setminus\{t^k:k\in\N\}$ a.\,s., where $\left|\pi\right|_1$ denotes the sum of the asymptotic frequencies of the blocks of a partition $\pi\in\p$.

It follows that the partitions $\Pi_{0,(\vartheta_\ep\wedge T)-}$ and $\Pi_{0,\vartheta_\ep\wedge T}$ have proper frequencies a.\,s. Hence, there exists a.\,s.\ $\ell\in\N$ such that
\[\sum_{i\in M(\Pi_{0,\vartheta_\ep\wedge T})\cap\sell}|B(\Pi_{0,\vartheta_\ep\wedge T},i)|>1-\ep\]
and
\[\sum_{i\in M(\Pi_{0,(\vartheta_\ep\wedge T)-})\cap\sell}|B(\Pi_{0,(\vartheta_\ep\wedge T)-},i)|>1-\ep.\]
By Lemma~\ref{Pathw:lem:part-nd}, there exists a.\,s.\ $\delta>0$ such that
\[\sum_{i\in M(\Pi_{0,t})\cap\sell}|B(\Pi_{0,t},i)|>1-\ep\]
for all $t\in(\vartheta_\ep\wedge T-\delta,\vartheta_\ep\wedge T+\delta)$. This implies $\vartheta_{\ep,j}\notin(\vartheta_\ep-\delta,\vartheta_\ep+\delta)$ for all $j\geq \ell$ a.\,s.\ on the event $\{\vartheta_\ep<T\}$, hence $\{\vartheta_\ep<T\}$ is a null event.

The assertion follows as $T\in(0,\infty)$ and $\ep\in(0,T)$ can be chosen arbitrarily.
\end{proof}
In the proof of Lemma~\ref{Pathw:lem:proper-dust} which is given below, we will use the following lemma which strengthens Lemma 11.2 of~\cite{Sampl}.
\begin{lem}
\label{Pathw:lem:v-dust}
Assume $\Xi\in\dust$ and let $t\in(0,\infty)$, $i\in\N$. Then a.\,s.\ on the event $\{v_t(i)<t\}$, there exists an integer $j\in\N\setminus\{i\}$ with
$v_t(i)=v_t(j)$ and $r_t(i,j)=0$.
\end{lem}
\begin{proof}
We use the points $(t^k,\pi^k)$ of $\eta_0$ from Section \ref{Pathw:sec:two-step-eta}.
W.\,l.\,o.\,g., we assume $\Xi_0(\Delta)>0$ and that the times $t^k\wedge t$ are stopping times with respect to the filtration $(\F_s,s\in\R_+)$ that is defined after equation \reff{Pathw:eq:def-J}.
The latter property can be achieved e.\,g.\ by setting $\Delta_{\ell}=\{x\in\Delta:1/(\ell+1)<x_1\leq 1/\ell\}$ and $t^{\ell,0}=0$ for each $\ell\in\N$,
and $t^{\ell,n}=\inf\{t>t^{\ell,n-1}:\zeta_0((t^{\ell,n-1},t]\times\Delta_{\ell})>0\}\wedge t$ for each $n\in\N$, and associating each $k\in\N$ with a pair $(\ell,n)$.

Lemma \ref{Pathw:lem:exch-Rt} yields for each $k\in\N$ that the sequence $(\I{t-v_t(j)\leq t^k},j\in\N)$ is exchangeable. Here we also use that this sequence is $\eta$-measurable so that we can assume w.\,l.\,o.\,g.\ that $(r_0,v_0)$ is exchangeable. The de Finetti theorem implies that a.\,s.\ either no or infinitely many elements of this sequence equal $1$.
Note that the random subset $A_k:=\{A_{t^k}(t,j):j\in\N,t-v_t(j)\leq t^k\}\subset\N$ is measurable with respect to $t^k$ and the process $(J_{t^k+s}-J_{t^k},s\in\R_+)$, where $J$ is defined in \reff{Pathw:eq:def-J}. Now the strong Markov property of $J$, our assumption on $t^k$, and the definition of $\pi^k$ yield that $A_k$ and $\pi^k$ are independent. As $\Xi\{0\}=0$, by Kingman's correspondence, and as $A_k$ and $\pi^k$ are independent, all non-singleton blocks of $\pi^k$ have an infinite intersection with $A_k$ a.\,s.\ on the event $\{\#A_k=\infty\}$. Furthermore, the definition of $v_{t^k}$ yields that $A_{t^k}(t,j)\neq A_{t^k}(t,j')$ for all distinct $j,j'\in\N$ with $t-v_t(j)\leq t^k$, $t-v_t(j')\leq t^k$. Hence, $\#A_k=\#\{j\in\N:t-v_t(j)\leq t^k\}$.

By definition of $v_t(i)$ and condition \reff{Pathw:eq:ass-eta-hat}, there exists a.\,s.\ on $\{v_t(i)<t\}$ an integer $k\in\N$ such that $t-v_t(i)=t^k$ and $A_{t^k}(t,i)$ is in a non-singleton block of $\pi^k$. The above implies that there exists a.\,s.\ on $\{v_t(i)<t\}$ an integer $j\in\N\setminus\{i\}$ such that $t-v_t(j)\leq t^k$ and $A_{t^k}(t,j)$ is in the same block of $\pi^k$ as $A_{t^k}(t,i)$. The definition of $v_t(j)$ now yields $v_t(j)=v_t(i)$ and $r_t(i,j)$=0.
\end{proof}
We will apply Lemma \ref{Pathw:lem:v-dust} through the following corollary.
\begin{cor}
\label{Pathw:cor:nosingl-dust}
Assume $\Xi\in\dust$ and that $(r_0,v_0)$ is exchangeable. Let $a\in\N$ and $\ep>0$. Then on an event of probability $1$, none of the partitions $\Pi^{a,\ep}_t$, $t\in\R_+$ contains singleton blocks.
\end{cor}
\begin{proof}
Let $q\in\R_+$. A.\,s.\ by Lemma \ref{Pathw:lem:v-dust}, only integers $i\in\N$ with $v_q(i)\geq q$ can form singleton blocks of $\Pi^{a,\ep}_q$. As those integers belong to finitely many blocks by definition of $\Pi^{a,\ep}_q$,
it follows that $\Pi^{a,\ep}_q$ contains a.\,s.\ at most finitely many singleton blocks. As the partition $\Pi^{a,\ep}_q$ is exchangeable by Lemma~\ref{Pathw:lem:exch-Rt}, Kingman's correspondence implies that $\Pi^{a,\ep}_q$ contains a.\,s.\ no singleton blocks.

On the event of probability $1$ on which condition~\reff{Pathw:eq:ass-eta-hat} holds, there exists for each $t\in\R_+$ and $i\in\N$ a time $q(t,i)\in(t,\infty)\cap\Q$ with $\eta((t,q(t,i)]\times\hat\p^i)=0$, as \reff{Pathw:eq:ass-eta-hat} implies that the points of $\eta(\cdot\times\hat\p^i)$ do not accumulate. By construction, if $i$ forms a singleton block in $\Pi^{a,\ep}_t$, then $i$ forms a singleton block also in $\Pi^{a,\ep}_{q(t,i)}$, as the particle remains on level $i$ and does not reproduce. This implies the assertion.
\end{proof}
\begin{proof}[Proof of Lemma~\ref{Pathw:lem:proper-dust}]
We proceed similarly to the proof of Lemma~\ref{Pathw:lem:proper-nd}. We assume w.\,l.\,o.\,g.\ $\Xi_0(\Delta)>0$.

\emph{Step 1. } For all $k\in\N$, the partition $\Pi^{a,\ep}_{t^k-}$ is exchangeable by Corollary~\ref{Pathw:cor:exch-R}. By Corollary~\ref{Pathw:cor:nosingl-dust} and Kingman's correspondence, it follows that $\Pi^{a,\ep}_{t^k-}$ has proper frequencies a.\,s.

\emph{Step 2. } We set for $\ell\in\N$ \[\vartheta_{\tilde\ep,\ell}=\inf\left\{t\geq 0:\sum_{i\in M(\Pi^{a,\ep}_t)\cap\sell}|B(\Pi^{a,\ep}_t,i)|<1-\tilde\ep\right\},\]
We deduce from Lemma~\ref{Pathw:lem:part-dust} that $\vartheta_{\tilde\ep,\ell}$ is a stopping time with respect to the usual augmentation of $\F^n$ for all $n\in\sell$. Then we define $\vartheta_{\tilde\ep}=\sup_{\ell\in\N}\vartheta_{\tilde\ep,\ell}$ which is for all $n\in\N$ a stopping time with respect to the usual augmentation of $\F^n$. Let $T\in\R_+$. The partition $\Pi^{a,\ep}_{\vartheta_{\tilde\ep}\wedge T}$ is exchangeable by Lemma~\ref{Pathw:lem:exch-J}. We deduce as in step 1 that it has proper frequencies a.\,s.

\emph{Step 3. } We conclude as in the proof of Lemma~\ref{Pathw:lem:proper-nd}, using Lemma~\ref{Pathw:lem:part-dust}.
\end{proof}

Recall the set of measures $\MCDI$ from Section~\ref{Pathw:sec:ld-sampl-nd}. The following lemma will be used in the proof of Theorem~\ref{Pathw:thm:ld-CDI}\ref{Pathw:item:thm:ld-CDI-supp}.
\begin{lem}
\label{Pathw:lem:Pi-supp}
Assume $\Xi\in\MCDI$. Let $s\in\R_+$. Then a.\,s., $|B(\Pi_{s,t},i)|>0$ for all $t\in(s,\infty)$ and $i\in\N$.
\end{lem}
\begin{proof}
Again we assume w.\,l.\,o.\,g.\ that $\rho_0$ is exchangeable, and it suffices to consider the case $s=0$. Let $\ep>0$. For each $k\in\N$, the process
\[\left(\min_{i\in\sk}|B(\Pi_{0,t},i)|,t\in\R_+\right)\]
is adapted with respect to the usual augmentation of the filtration $\F^k$ (defined in Section~\ref{Pathw:sec:exch-ld}) and has a.\,s.\ càdlàg paths by Lemma~\ref{Pathw:lem:part-nd}. Consequently,
\[\vartheta_k:=\inf\left\{t\geq \ep:
\min_{i\in\sk}|B(\Pi_{0,t},i)|=0\right\}\]
is a stopping time with respect to the usual augmentation of $\F^n$ for all $n\in\sk$. It follows that
\[\vartheta:=\inf_{k\in\N}\vartheta_k=\lim_{k\to\infty}\vartheta_k\]
is a stopping time with respect to the usual augmentation of $\F^n$ for all $n\in\N$. By Lemma~\ref{Pathw:lem:exch-J}, the distance matrix $\gamma_n(\rho_{\vartheta\wedge T})$ is exchangeable for each $T\in [\ep,\infty)$ and $n\in\N$. Hence, the partition $\Pi_{0,\vartheta\wedge T}$ is exchangeable. The assumption $\Xi\in\MCDI$ implies $\#\Pi_{0,\vartheta\wedge T}<\infty$ a.\,s. Kingman's correspondence now implies that each block of $\Pi_{0,\vartheta\wedge T}$ has a positive asymptotic frequency a.\,s. Hence, by Lemma~\ref{Pathw:lem:part-nd}, there exists an event of probability $1$ on which all blocks of the partitions $\Pi_{0,t}$ with $t$ in a right neighborhood of $\vartheta\wedge T$ have positive asymptotic frequencies. By definition of $\vartheta$, it follows $\P(\vartheta<T)=0$. The assertion follows as $T$ and $\ep$ can be chosen arbitrarily.
\end{proof}

\section{The construction on the lookdown space}
\label{Pathw:sec:constr-ld}
Now we apply the results from the last section to prove the assertions from Section~\ref{Pathw:sec:ld-sampl}. We use the coupling characterization of the Prohorov distance, namely (see e.\,g.\ \cite{EK86}*{Theorem 3.1.2}) that in a separable metric space $(Y,d)$, the Prohorov distance between probability measures $\mu$ and $\mu'$ on the Borel sigma algebra is given by
\begin{equation*}
\label{Pathw:eq:Proh-coupling}
\dP^Y(\mu,\mu')=\inf_\nu\inf\{\ep>0:\nu\{(y,y')\in Y^2:d(y,y')>\ep\}<\ep\},
\end{equation*}
where the first infimum is over all couplings $\nu$ of the probability measures $\mu$ and $\mu'$.

\subsection{The case with dust}
\label{sec:constr-ld-dust}
In this subsection, we always consider the case $\Xi\in\dust$. Let $(X,r,m)$ be a marked metric measure space, and let $X\times\R_+$ be endowed with the product metric $d((x,v),(x',v'))=r(x,x')\vee|v-v'|$. Let $(x(i),v(i))_{i\in\N}$ be an $m$-iid sequence in $X\times\R_+$ that is independent of $\eta$. We set $(r_0,v_0)=((r(x(i),x(j)))_{i,j\in\N},v)$. Then $(r_0,v_0)$ is distributed according to the marked distance matrix distribution of $(X,r,m)$. With the extended lookdown space $(\hat Z,\rho)$ associated with $\eta$ and $(r_0,v_0)$, we are in the setting of Subsection~\ref{Pathw:sec:ld-sampl-dust}.

\begin{proof}[Proof of Theorem~\ref{Pathw:thm:ld-dust} (beginning)]
We begin with the proof of item \ref{Pathw:eq:thm:ld-dust:conv}.
For $\ep>0$ and $a,n\in\N$, we define the probability measures
\[m^{a,\ep,n}_t=\sum_{i\in M(\Pi^{a,\ep}_t)}\left|B(\Pi^{a,\ep}_t,i)\right|_n\delta_{(z(t,i),v_t(i))}.\]
on $\hat Z\times\R_+$. Clearly, $\left|B(\Pi^{a,\ep}_t,i)\right|_n=0$ for all $i\in M(\Pi^{a,\ep}_t)$ with $i>n$.

Let $T\in\R_+$. Using that the Prohorov distance is bounded from above by the total variation distance, we obtain
\begin{align*}
&\lim_{n,\ell\to\infty}\sup_{t\in[0,T]}\dP^{\hat Z\times\R_+}(m^{a,\ep,n}_t,m^{a,\ep,\ell}_t)\\
\leq&\lim_{n,\ell\to\infty}\sup_{t\in[0,T]}\sum_{i\in M(\Pi^{a,\ep}_t)}
\big|\left|B(\Pi^{a,\ep}_t,i)\right|_n-\left|B(\Pi^{a,\ep}_t,i)\right|_\ell\big|\\
\leq&\lim_{k\to\infty}\lim_{n,\ell\to\infty}\sup_{t\in[0,T]}\sum_{i\in M(\Pi^{a,\ep}_t)\cap\sk}
\big|\left|B(\Pi^{a,\ep}_t,i)\right|_n-\left|B(\Pi^{a,\ep}_t,i)\right|_\ell\big|\\
&+2\lim_{k\to\infty}\lim_{n\to\infty}\sup_{t\in[0,T]}\sum_{\substack{i\in M(\Pi^{a,\ep}_t):\\i>k}}\left|B(\Pi^{a,\ep}_t,i)\right|_n
=0\quad\text{a.\,s.}
\end{align*}
The first summand on the right-hand side equals zero a.\,s.\ by Lemma~\ref{Pathw:lem:part-dust}.
A.\,s., the second summand equals zero as for each $\tilde\ep>0$, there exist integers $k$ and $n_0$ such that
\[\inf_{t\in[0,T]}\sum_{i\in M(\Pi^{a,\ep}_t)\cap\sk}
\left|B(\Pi^{a,\ep}_t,i)\right|_n\geq 1-\tilde\ep\]
for all $n\geq n_0$ by Lemmas~\ref{Pathw:lem:proper-dust} and \ref{Pathw:lem:part-dust}.

Now we compare the probability measures $m^n_t$ and $m^{a,\ep,n}_t$. A coupling $\nu$ of these probability measures is given by
\[\nu=\sum_{i=1}^n\frac{1}{n}\delta_{((z(t,i)),v_t(i)),(z(t,j)),v_t(j)))},\]
where $j=\min B(\Pi^{a,\ep}_t,i)$ in each summand.
By definition of $\Pi^{a,\ep}_t$ and Remark~\ref{Pathw:rem:AI},
\[d^{\hat Z\times\R_+}((z(t,i),v_t(i)),(z(t,j),v_t(j)))=r_t(i,j)\vee\left|v_t(i)-v_t(j)\right|<2\ep\]
for all $i,j\in\N$ that are in the same block of $\Pi^{a,\ep}_t$ and not in $C^{a,\ep,\emptyset}_t$. The coupling characterization of the Prohorov metric implies
\begin{align*}
&\dP^{\hat Z\times\R_+}(m^n_t,m^{a,\ep,n}_t)\\
&\leq\nu\{(y,y')\in (\hat Z\times\R_+)^2:d^{\hat Z\times\R_+}(y,y')\geq 2\ep\}+2\ep\\
&\leq |C^{a,\ep,\emptyset}_t|_n+2\ep.
\end{align*}

By construction, $|C^{a,\ep,\emptyset}_t|_n\leq|C^{a,\ep,\emptyset}_0|_n$.
This follows from the definition of $C_t^{a,\ep,\emptyset}$ in Section \ref{Pathw:sec:part}, from the definition of $(r_t,v_t)$, as a particle at time $s$ on a level $i$ loses the property that $v_{s}(i)\geq s$ if it reproduces, and as it can only increase its level in a reproduction event.
By exchangeability (or Lemma~\ref{Pathw:lem:part-dust}), we have
$\lim_{n\to\infty}|C^{a,\ep,\emptyset}_0|_n=|C^{a,\ep,\emptyset}_0|$ a.\,s.

The triangle inequality yields
\begin{equation}
\label{Pathw:eq:thm:dust-triangle}
\lim_{n,\ell\to\infty}\sup_{t\in[0,T]}\dP^{\hat Z\times\R_+}(m^n_t,m^\ell_t)
\leq 2|C^{a,\ep,\emptyset}_0|+4\ep\quad\text{a.\,s.}
\end{equation}
Letting first $a\to\infty$ and then $\ep\downarrow 0$, we obtain from Lemma~\ref{Pathw:lem:A0} below that the left-hand side of \reff{Pathw:eq:thm:dust-triangle} equals zero a.\,s. By definition of $m_t$ and as $\hat Z\times\R_+$ is complete, this implies assertion~\ref{Pathw:eq:thm:ld-dust:conv}.

As recalled in Section~\ref{Pathw:sec:gen-setting}, the map $t\mapsto (z(t,i),v_t(i))$ is a.\,s.\ càdlàg with respect to $d^{\hat Z\times\R_+}$, hence the map $t\mapsto m^n_t$ is a.\,s.\ càdlàg in the weak topology on $\hat Z\times\R_+$. Jump times can only lie in the set $\Theta_0$, which equals a.\,s.\ the set of reproduction times.
The uniformity of the convergence in assertion~\ref{Pathw:eq:thm:ld-dust:conv} implies that also $t\mapsto m_t$ is a.\,s.\ càdlàg in the weak topology on $\hat Z\times\R_+$ with no jump times outside $\Theta_0$.

W.\,l.\,o.\,g., we assume $\Xi(\Delta)>0$. As $\Xi\in\dust$, this implies $\Xi_0(\Delta)>0$ and $\Theta_0=\{t^k:k\in\N\}$ a.\,s.
Now we deduce that $m_{t^k}(\hat Z\times\{0\})>0$ a.\,s.\ all for $k\in\N$, which is part of assertion~\ref{Pathw:item:thm:ld-dust:zero}.
By Corollary \ref{Pathw:cor:exch-R}, the sequence $(v_{t^k}(i),i\in\N)$ is exchangeable. From the definition \reff{Pathw:eq:Theta0} of $\Theta_0$, the definition of $v_t(i)$ in Section \ref{Pathw:sec:ld-dec}, and as condition \reff{Pathw:eq:ass-eta-hat} is a.\,s.\ satisfied, it follows that $|\{i\in\N:v_{t^k}(i)=0\}|>0$ a.\,s.
Hence, the de Finetti theorem and the definition of $m_t$ yield that a.\,s., the empirical measures $n^{-1}\sum_{i=1}^n\delta_{v_{t^k}(i)}$ converge weakly to the directing measure $m_{t^k}(\hat Z\times\cdot)$ on $\R_+$ which satisfies $m_{t^k}(\hat Z\times\{0\})>0$.

On the event of probability $1$ on which condition \reff{Pathw:eq:ass-eta-hat} is satisfied, $v_{t^k-}(i)>0$ for all $i\in\N$ a.\,s. Analogously to the above, this yields $m_{t^k-}(\hat Z\times\{0\})=0$ a.\,s.
This is another part of assertion \ref{Pathw:item:thm:ld-dust:zero}. As a consequence, the set of jump times is a.\,s.\ not smaller than $\Theta_0$, which yields assertion \ref{Pathw:item:thm:ld-dust:reg}. After completing the proof of assertion~\ref{Pathw:eq:thm:ld-dust:conv} by proving Lemma \ref{Pathw:lem:A0} below, we give in Proposition \ref{Pathw:prop:repr-mt} a representation of the probability measures $m_t$ from which assertion \ref{Pathw:item:thm:ld-dust:atomic} will follow by Remark \ref{Pathw:rem:repr-mt-leli}.
A.\,s., $v_t(i)>0$ (by condition~\reff{Pathw:eq:ass-eta-hat}) and $m_t=m_{t-}$ (by item~\ref{Pathw:item:thm:ld-dust:reg}) for all $t\in(0,\infty)\setminus\Theta_0$ and $i\in\N$. Hence, the remainder of assertion \ref{Pathw:item:thm:ld-dust:zero} also follows from Proposition~\ref{Pathw:prop:repr-mt} below.
\end{proof}
\begin{lem}
\label{Pathw:lem:A0}
Let $\ep>0$, and let $C^{a,\ep,\emptyset}_0$ be defined as in this subsection from a random variable $(r_0,v_0)$ that has the marked distance matrix distribution of a marked metric measure space $(X,r,m)$. Then,
$\lim_{a\to\infty}|C^{a,\ep,\emptyset}_0|=0$ a.\,s.
\end{lem}
\begin{proof}
From the definitions of $C^{a,\ep,\emptyset}_0$ and of the extended lookdown space, we obtain
\[C^{a,\ep,\emptyset}_0=\{j\in\N: (x(j),v(j))\in(X\times\R_+)\setminus\bigcup_{i=1}^a\oB^{X\times\R_+}_\ep(x(i),v(i))\}\]
for all $a\in\N$, where $\ep$-balls are defined by
$\oB^{X\times\R_+}_\ep(x',v')=\{(x'',v'')\in X\times\R_+:r(x',x'')\vee|v'-v''|<\ep\}$ for $(x',v')\in X\times\R_+$.

As $(x(a+j),v(a+j))_{j\in\N}$ is an $m$-iid sequence in $X\times\R_+$ that is independent of $(x(i),v(i))_{i\in[a]}$, the law of large numbers yields
\[|C^{a,\ep,\emptyset}_0|=m((X\times\R_+)\setminus\bigcup_{i=1}^a\oB^{X\times\R_+}_\ep(x(i),v(i)))\quad\text{a.\,s.}\]
By separability, $X\times\R_+$ can be covered by countably many balls of diameter $\ep/2$, and each ball with positive mass contains elements of the sequence $(x(i),v(i))_{i\in\N}$ a.\,s. Using also continuity of $m$ from above, this implies 
\[\lim_{a\to\infty} m((X\times\R_+)\setminus\bigcup_{i=1}^a\oB^{X\times\R_+}_\ep(x(i),v(i)))=0\quad\text{a.\,s.}\]
This yields the assertion.
\end{proof}

Now we give an explicit representation of the probability measures $m_t$. For $t\in\R_+$, we denote by $\Pi_t$\label{Pathw:not:Pit} the partition of $\N$ in which integers $i,j$ are in the same block if and only if
$v_t(i)=v_t(j)<t$ and $r_t(i,j)=0$,
which is condition~\ref{Pathw:item:Pit-above} on p.\,\pageref{Pathw:item:Pit-above}. We may call $\Pi_t$ the partition of siblings.

The individuals at time $t$ whose parents are also parents of individuals at time $0$ are on the levels in the set
\[C_t=\{i\in\N:v_t(i)\geq t\}.\]
All elements of $C_t$ form singleton blocks in $\Pi_t$. For $a\in\N,\ep>0$, the non-singleton blocks of $\Pi_t$ are also blocks of $\Pi_t^{a,\ep}$. Hence, Corollary~\ref{Pathw:cor:nosingl-dust} implies that a.\,s., $C_t$ equals the union of the singleton blocks of $\Pi_t$ for each $t\in\R_+$.
We also define the map
\[\theta_t:\hat Z\times\R_+\to \hat Z\times\R_+,\quad (z,s)\mapsto (z,s+t)\]
with is continuous for $d^{\hat Z\times\R_+}$.
\begin{prop}
\label{Pathw:prop:repr-mt}
Assume $\Xi\in\dust$. Then on an event of probability $1$,
\begin{equation}
\label{Pathw:eq:mu_t-dust}
m_t=\sum_{i\in M(\Pi_t)} \left|B(\Pi_t,i)\right|\delta_{(z(t,i),v_t(i))}+\left|C_t\right| \theta_t(m_0)
\end{equation}
for all $t\in\R_+$.
\end{prop}
\begin{rem}
Proposition~\ref{Pathw:prop:repr-mt} allows to describe the probability measure $m_t(\cdot\times\R_+)$ on the extended lookdown space $\hat Z$ as follows. With probability given by the asymptotic frequency of the individuals at time $t$ whose ancestral lineages do not coalesce with other ancestral lineages within the time interval $(0,t]$, we sample according to $m_0(\cdot\times\R_+)$. For each block in a reproduction event at a time $\tau'$ in $(0,t]$, we draw the individual on, say, the lowest level in this block (which is identified with the individuals on all other levels in this block, as they have genealogical distance zero) with probability given by the asymptotic frequency of the individuals at time $t$ that descend from this block and whose ancestral lineages do not coalesce with any other ancestral lineages in the time interval $(\tau',t]$.
\end{rem}
\begin{rem}
\label{Pathw:rem:repr-mt-leli}
The marked metric measure spaces $(X,r,m)$ and $(\hat Z,\rho,m_0)$ are isomorphic (as defined in Section \ref{Pathw:sec:proc-gen}). This follows from the definition of $(\hat Z,\rho)$ and the Gromov reconstruction theorem (cf.\ e.\,g.\ \cite{Sampl}*{Proposition 10.5}). Hence, $m_0$ is purely atomic if and only if $m$ is purely atomic, and the assertion on the atomicity of $m_t$ in Theorem \ref{Pathw:thm:ld-dust}\ref{Pathw:item:thm:ld-dust:atomic} follows from Proposition~\ref{Pathw:prop:repr-mt}.

Proposition~\ref{Pathw:prop:repr-mt} and Lemmas \ref{Pathw:lem:part-dust} and \ref{Pathw:lem:proper-dust} imply that on an event of probability $1$, each weak limit $m_{t-}$ with $t\in(0,\infty)$ is also the sum of countably many atoms and a multiple of $\theta_t(m_0)$. Here we also use that the map $\R_+\times\hat Z\times\R_+\to\hat Z\times\R_+$, $(s,z',v')\mapsto\theta_s(z',v')$ is continuous also in $s$. This yields the assertion on the atomicity of the left limits $m_{t-}$ in Theorem \ref{Pathw:thm:ld-dust}\ref{Pathw:item:thm:ld-dust:atomic}.
\end{rem}
\begin{proof}[Proof of Proposition~\ref{Pathw:prop:repr-mt}]
For all $a\in\N$, $\ep>0$, and $t\in\R_+$, the definition of $C^{a,\ep,I}_t$ implies
\[C_t=\bigcup_{I\subset\sa}C^{a,\ep,I}_t.\]
Lemma~\ref{Pathw:lem:part-dust} implies the existence of the asymptotic frequencies
\[|C_t|=\sum_{I\subset\sa}|C^{a,\ep,I}_t|\]
for all $t\in\R_+$ a.\,s.
Let us denote the right-hand side of equation~\reff{Pathw:eq:mu_t-dust} by $m_t'$.
On an event of probability one, $m_t'$ is a well-defined probability measure for all $t\in\R_+$, and $t\mapsto m_t'$ is càdlàg with respect to $\dP^{\hat Z\times\R_+}$.
This follows from Lemmas~\ref{Pathw:lem:part-dust} and \ref{Pathw:lem:proper-dust}, and as a.\,s., $C_t$ equals the union of the singleton blocks of $\Pi_t$ for each $t$.
As also $t\mapsto m_t$ is a.\,s.\ càdlàg, it suffices to show that \reff{Pathw:eq:mu_t-dust} holds a.\,s.\ for a fixed $t\in\R_+$.

For $i,j\in\N$ that are in the same block of $\Pi_t$, we have $v_t(i)=v_t(j)$ and $z(t,i)=z(t,j)$ in $\hat Z$ by definition. Hence,
\[m^n_t=\sum_{i\in M(\Pi_t)}|B(\Pi_t,i)|_n\delta_{(z(t,i),v_t(i))}.\]

Let $A^n_t=\{A_0(t,i):i\in C_t\cap\sn\}$ be the set of the ancestral levels at time $0$ of the individuals on the levels in $C_t\cap\sn$ at time $t$. For finite sets $A\subset\N$, let
\[m_0^A=\frac{1}{\#A\vee 1}\sum_{i\in A}\delta_{(z(0,i),v_0(i))}.\]
Using the definitions of $z(t,i)$ and $v_t(i)$,
we can write
\begin{equation}
\label{Pathw:eq:repr-mt-Cn}
m^n_t=\sum_{i\in M(\Pi_t)\setminus C_t}
|B(\Pi_t,i)|_n\delta_{(z(t,i),v_t(i))}
+|C_t|_n\theta_t(m^{A^n_t}_0).
\end{equation}

Note that $A^n_t$ is the set of those levels at time zero that are occupied by particles that do not reproduce until time $t$ and that are not above level $n$ at time $t$.
This implies that $C_t\cap[n]$ and $A^n_t$ are bijective. In particular, $\bigcup_{n\in\N}A^n_t$ is infinite on the event $\{|C_t|>0\}$.

Now we work with arguments that we encountered already in Lemmas~\ref{Pathw:lem:mu0} and \ref{Pathw:lem:m0}. Recall that $(r_0,v_0)=((r(x(i),x(j)))_{i,j\in\N},(v(i)_{i\in\N})$ where $(x(i),v(i))_{i\in\N}$ is an $m$-iid sequence in $X\times\R_+$ that is independent of $\eta$. For $A\subset\N$, we define the empirical measure
\[m^{A}=\frac{1}{\#A\vee 1}
\sum_{i\in A}\delta_{(x(i),v(i))}\]
on $X\times\R_+$.
As the sets $A^n_t$ and $C_t$ are $\eta$-measurable, the Glivenko-Cantelli theorem implies that the weak convergence
\[m=\text{w-}\lim_{n\to\infty}m^{A^n_t}\]
holds a.\,s.\ on $\{|C_t|>0\}$.
By construction of the lookdown space $\hat Z$,
the map
$\{(x(i),v(i)):i\in\N\}\to \hat Z\times\R_+$, $(x(i),v(i))\mapsto (i,v(i))$ can be extended a.\,s.\ to an isometry $\varphi$ from the closed support of $m$ in $X\times\R_+$ to $\hat Z\times\R_+$ which satisfies $m_0^{A^n_t}=\varphi(m^{A^n_t})$ and $m_0=\varphi(m)$. Thus the weak convergence
\[m_0=\text{w-}\lim_{n\to\infty}m^{A^n_t}_0\]
holds a.\,s.\ on $\{|C_t|>0\}$.

By exchangeability (or Lemma~\ref{Pathw:lem:part-dust}), also the relative frequencies in expression~\reff{Pathw:eq:repr-mt-Cn} converge a.\,s.\ to the corresponding asymptotic frequencies.
Hence, $m^n_t$ converges weakly to $m'_t$ on an event of probability $1$. This yields the assertion.
\end{proof}

\subsection{The case without dust}
\label{sec:constr-ld-nd}
In this subsection, we always consider the case $\Xi\in\nd$. Let $(x(i),i\in\N)$ be an iid sequence in a metric measure space $(X,r,\mu)$ that is independent of $\eta$. We assume $\rho_0=(r(x(i),x(j)))_{i,j\in\N}$.
With the lookdown space $(Z,\rho)$ associated with $\eta$ and $\rho_0$, we are in the setting of Section~\ref{Pathw:sec:ld-sampl-nd}.

\begin{proof}[Proof of Theorem~\ref{Pathw:thm:ld-nd} (beginning)]
We work on an event of probability $1$ on which in particular the assertions of Lemmas~\ref{Pathw:lem:part-nd} and \ref{Pathw:lem:proper-nd} hold simultaneously for all $s\in\Q_+$, and we mostly omit `a.\,s.' We define for each $t\in(0,\infty)$, $s\in(0,t)\cap\Q$, and $n\in\N$ a probability measure $\mu^{(n,s)}_t$ on $Z$ by
\[\mu^{(n,s)}_t=\sum_{i\in M(\Pi_{s,t})}\left|B(\Pi_{s,t},i)\right|_n
\delta_{(t,i)}.\]
There exists a coupling $\nu$ of the probability measures $\mu^{(n,s)}_t$ and $\mu^n_t$ given by
\[\nu=\sum_{i=1}^n\frac{1}{n}\delta_{((t,\min B(\Pi_{s,t},i)),(t,i))}.\]
As $\rho((t,\min B(\Pi_{s,t},i)),(t,i))\leq 2(t-s)$, the coupling characterization of the Prohorov metric implies
\begin{equation}
\label{Pathw:eq:muns-mun}
\dP^Z(\mu^{(n,s)}_t,\mu^n_t)\leq 2(t-s).
\end{equation}

Let $\ep\in(0,\infty)\cap\Q$, $T\in(\ep,\infty)$, and $\tau\in(0,\ep)\cap\Q$. First we consider $t\in[\tau,T]$.  Let $s_0=0$ and $s_j=\tau+(j-1)\ep$ for $j\in\N$. By Lemma~\ref{Pathw:lem:proper-nd}, there exists for each $j\in\N_0$ an integer $\ell_j$ such that
\[\sum_{i\in M(\Pi_{s_j,t})\cap[ \ell_j]}|B(\Pi_{s_j,t},i)|>1-\ep\]
for all $t\in[s_{j+1},s_{j+2}]$. We set $\ell=\max\{\ell_j:j\in\N_0,s_{j+1}\leq T\}$. By Lemma~\ref{Pathw:lem:part-nd}, there exists an integer $n'$ such that
\begin{equation}
\label{Pathw:eq:n'}
\sum_{i\in M(\Pi_{s_j,t})\cap\sell}|B(\Pi_{s_j,t},i)|_n>1-\ep
\end{equation}
for all $n\geq n'$, $j\in\N_0$ with $s_{j+1}\leq T$, and $t\in[s_{j+1},s_{j+2}]$.
For all $t\in[\tau,T]$, all $k,n\geq n'$, and $j\in\N_0$ such that $t\in[s_{j+1},s_{j+2}]$, the bound \reff{Pathw:eq:n'} yields
\[\dP^Z(\mu^{(n,s_j)}_t,\mu^{(k,s_j)}_t)\leq
\sum_{\substack{i\in M(\Pi_{s_j,t}):\\ i\leq \ell}}
\big||B(\Pi_{s_j,t},i)|_n-|B(\Pi_{s_j,t},i)|_k\big|
+\ep\]
as the Prohorov distance is bounded from above by the total variation distance.
By Lemma~\ref{Pathw:lem:part-nd}, this expression converges to $\ep$ uniformly in $t\in[\tau,T]$ as $n,k\to\infty$.
Furthermore, for all $n,k\in\N$, $t\in[\tau,T]$, and $j$ with $t\in[s_{j+1},s_{j+2}]$, the bound \reff{Pathw:eq:muns-mun} yields
\begin{align*}
\dP^Z(\mu^n_t,\mu^k_t)&\leq\dP^Z(\mu^n_t,\mu^{(n,s_j)}_t)+\dP^Z(\mu^{(n,s_j)}_t,\mu^{(k,s_j)}_t)+\dP^Z(\mu^{(k,s_j)}_t,\mu^k_t)\\
&\leq 8\ep + \dP^Z(\mu^{(n,s_j)}_t,\mu^{(k,s_j)}_t).
\end{align*}
Hence,
$\lim_{n\to\infty}\sup_{t\in[\tau,T]}\dP^Z(\mu^n_t,\mu_t)\leq 8\ep$
for all $\tau\in(0,\ep)\cap\Q$ a.\,s.
\end{proof}
In the next part of the proof, we choose a random $\tau$ and show uniform convergence of $\mu^n_t$ in $[0,\tau]$. To this aim, we define for each $b\in\N$, $\ep>0$, and $t\in\R_+$ a collection $(A^{b,\ep,I}_t,I\subset[b])$ of disjoint subsets of $\N$ (whose union is $\N$) by
\[A^{b,\ep,I}_t=\bigcap_{j\in I}\{i\in\N: \rho_t(j,i)<\ep+2t\}
\cap\bigcap_{k\in\sb\setminus I}\{i\in\N:\rho_t(k,i)\geq \ep+2t\}.\]
This partitions the set of individuals at time $t$ into blocks such that, if $t$ is small, then two individuals are close if they are in a common block $A^{b,\ep,I}_t$ with $I\neq\emptyset$.
\begin{lem}
\label{Pathw:lem:var}
Let $\tilde T,\ep>0$, $b\in\N$, and $I\subset\sb$.
Then
\[\lim_{n\to\infty}\sup_{t\in[0,\tilde T]}\big||A^{b,\ep,I}_t|_n-|A^{b,\ep,I}_t|\big|=0\quad\text{a.\,s.}\]
\end{lem}
\begin{proof}
In Lemma~\ref{Pathw:lem:BC-bound}, choose $f$ such that
\[f(t,(\rho_t(j,b+1+i))_{j\in\sb})
=\prod_{j\in I}\I{\rho_t(j,b+1+i)<\ep+2t}
\prod_{k\in \sb\setminus I}\I{\rho_t(k,b+1+i)\geq \ep+2t}\]
for all $t\in\R_+$ and $i\in\N$. The product on the right-hand side is the indicator variable of $\{b+1+i\in A^{b,\ep,I}_t\}$. By the same argument as in the proof of Lemma \ref{Pathw:lem:part-nd}, condition \reff{eq:cond-newb} is satisfied. The assertion now follows from Lemma~\ref{Pathw:lem:BC-bound} as
\[X_n(t)\leq\frac{b+1+n}{n}|A^{b,\ep,I}_t|_{b+1+n}\leq X_n(t)+\frac{b+1}{n}\]
and $X(t)=|A^{b,\ep,I}_t|$ for all $t\in\R_+$ a.\,s., where $X_n(t)$ and $X(t)$ are defined in Lemma~\ref{Pathw:lem:BC-bound}.
\end{proof}
In a metric space $(Y,d)$, we denote by $\oB^Y_\ep(z)=\{y'\in Y:d(y',y)<\ep\}$ for $y\in Y$ some $\ep$-balls.
\begin{lem}
\label{Pathw:lem:A-lln}
Let $b\in\N$ and $\ep>0$. Then
\[|A^{b,\ep,\emptyset}_0|=1-\mu(\bigcup_{i=1}^b\oB^X_\ep(x(i)))\quad\text{a.\,s.}\]
\end{lem}
\begin{proof}
By the definitions of $A^{b,\ep,\emptyset}_0$ and $Z$,
\[A^{b,\ep,\emptyset}_0=\{j\in\N:x(j)\notin\bigcup_{i=1}^b\oB^X_\ep(x(i))\}.\]
Clearly, $(x(b+i),i\in\N)$ is a $\mu$-iid sequence in $X$ that is independent of $(x(i),i\in[b])$. The assertion follows from the law of large numbers.
\end{proof}
\begin{proof}[Proof of Theorem~\ref{Pathw:thm:ld-nd} (end)]
By separability of $(X,r)$ and continuity of $\mu$ from above (analogously to the proof of Lemma \ref{Pathw:lem:A0}), we can (and do) choose a (random) integer $b$ such that
\[\mu\left(\bigcup_{i=1}^b\oB^X_\ep(x(i))\right)> 1-\ep.\]
By condition \reff{Pathw:eq:ass-eta} and Lemmas \ref{Pathw:lem:A-lln} and \ref{Pathw:lem:var}, we can (and do) choose (random) $\tau\in(0,\ep)\cap\Q$ and $n_0\in\N$ such that $\eta((0,\tau]\times\p^b)=0$ and $|A^{b,\ep,\emptyset}_t|_n<\ep$ for all $t\in[0,\tau]$ and $n\geq n_0$.

By construction, $\rho((t,i),(t,j))<2\ep+4t$ for all $t\in\R_+$, $I\neq\emptyset$, and $i,j\in A^{b,\ep,I}_t$. 
Using a coupling $\nu$ of $\mu^n_t$ and $\mu^k_t$ such that
\[\nu((\{t\}\times A^{b,\ep,I}_t)\times(\{t\}\times A^{b,\ep,I}_t))=\mu^n_t(\{t\}\times A^{b,\ep,I}_t)\wedge \mu^k_t(\{t\}\times A^{b,\ep,I}_t)\]
for all $I\subset\sb$,
we obtain for $n,k\geq n_0$ and $t\in[0,\tau]$ from the coupling characterization of the Prohorov metric that
\begin{align*}
&\dP^Z(\mu^n_t,\mu^k_t)\\
&\leq\nu((z,z')\in Z^2:\rho(z,z')>2\ep+4t)+2\ep+4t\\
&\leq\nu(Z\times(\{t\}\times A^{b,\ep,\emptyset}_t))
+\sum_{\substack{I\subset\sb:\\I\neq\emptyset}}
\nu((Z\setminus(\{t\}\times A^{b,\ep,I}_t))\times(\{t\}\times A^{b,\ep,I}_t))+2\ep+4t\\
&\leq \mu^k_t(\{t\}\times A^{b,\ep,\emptyset}_t)+\sum_{\substack{I\subset\sb:\\I\neq\emptyset}}\left|\mu^n_t(\{t\}\times A^{b,\ep,I}_t)-\mu^k_t(\{t\}\times A^{b,\ep,I}_t)\right|+2\ep+4t\\
&\leq 7\ep+\sum_{\substack{I\subset\sb:\\I\neq\emptyset}}
\big||A^{b,\ep,I}_t|_n-|A^{b,\ep,I}_t|_k\big|.
\end{align*}
Lemma~\ref{Pathw:lem:var} implies
$\lim_{n,k\to\infty}\sup_{t\in[0,\tau]}\dP^Z(\mu^n_t,\mu^k_t)=7\ep$.
Altogether,
\[\lim_{n,k\to\infty}\sup_{t\in[0,T]}\dP^Z(\mu^n_t,\mu^k_t)=8\ep.\]
Assertion~\ref{Pathw:eq:thm:ld-nd:conv} follows as $\ep$ can be chosen arbitrarily small, and as $(Z,\rho)$ is complete.

Now we come to assertion~\ref{Pathw:item:thm:ld-nd:reg}. By condition~\reff{Pathw:eq:ass-eta}, the map $\R_+\to Z$, $t\mapsto (t,i)$ is càdlàg for each $i\in\N$. Hence, the map $t\mapsto\mu^n_t$ is càdlàg for each $n\in\N$, and by item~\ref{Pathw:eq:thm:ld-nd:conv} also the map $t\mapsto \mu_t$ is càdlàg.
For $n\in\N$ and $t\in(0,\infty)$, we define the probability measure
\[\mu^n_{t-}=\frac{1}{n}\sum_{i=1}^n\delta_{(t-,i)}\]
on $Z$, with the left limit $(t-,i)=\lim_{s\uparrow t}(s,i)$. Then $\mu^n_{t-}=\text{w-}\lim_{s\uparrow t}\mu^n_s$. Let $\mu_{t-}=\text{w-}\lim_{s\uparrow t}\mu_s$. From the uniform convergence in item~\ref{Pathw:eq:thm:ld-nd:conv}, it follows that $\mu^n_{t-}$ converges weakly to $\mu_{t-}$ as $n\to\infty$. Note that $\dP^Z(\mu^n_{t-},\mu^n_t)\leq 1/n$ for all $t\in(0,\infty)\setminus \Theta_0$ and all $n\in \N$ a.\,s., as only a binary reproduction event can occur at such a time $t$.
It follows that $\mu_{t-}=\mu_t$ for all $t\in(0,\infty)\setminus \Theta_0$ a.\,s. That the set of jump times is not smaller than $\Theta_0$ follows from assertion~\ref{Pathw:item:thm:ld-nd:atoms-jumps} which we now prove.

W.\,l.\,o.\,g., we assume $\Xi_0(\Delta)>0$, then we have $\Theta_0=\{t^k:k\in\N\}$ a.\,s.
(For $\Xi_0(\Delta)=0$, nothing remains to prove as $\Theta_0=\emptyset$ a.\,s.\ in this case.)
Condition \reff{Pathw:eq:ass-eta} implies that $\rho_{t^k-}(i,j)>0$ for all $k,i,j\in\N$ a.\,s.
For each $k\in\N$, the random variable $\rho_{t^k-}$ is exchangeable by Corollary \ref{Pathw:cor:exch-R}.
On an event of probability $1$, let $\chi$ be the isomorphy class of the metric measure space $(Z,\rho,\mu_{t^k-})$. As in Remark~\ref{Pathw:rem:psi}, we have $\chi=\psi(\rho_{t^k-})$. Let $\rho'$ be a random variable whose conditional distribution given $\chi$ is the distance matrix distribution of $\chi$. Then by \cite{Sampl}*{Remark 10.4}, the random variables $\rho'$ and $\rho_{t^k-}$ are (unconditionally) equal in distribution. Hence, $\rho'(i,j)>0$ a.\,s.\ for all $i,j\in\N$ which implies that $\mu_{t^k-}$ is a.\,s.\ non-atomic.

By the definition of the population model in Section \ref{Pathw:sec:ld-space} and the definition \reff{Pathw:eq:Theta0} of $\Theta_0$, there exists for each $k\in\N$ an $i\in\N$ such that
$|\{j\in\N: \rho_{t^k}(i,j)=0\}|=|B(\pi^k,i)|>0$. It can now be shown as above that $\mu_{t^k}$ contains an atom. More simply, the Portmanteau theorem and item \ref{Pathw:eq:thm:ld-nd:conv} imply
\[\mu_{t^k}\{(t,i)\}\geq\limsup_{n\to\infty}\mu_{t^k}^n\{(t,i)\}=|\{j\in\N: \rho_{t^k}(i,j)=0\}|>0.\]
\end{proof}

\begin{proof}[Proof of Theorem~\ref{Pathw:thm:ld-CDI}]
Recall that $\supp\mu_t$ denotes the closed support of $\mu_t$.
For $z\in Z$ and $\ep>0$, we denote $\ep$-balls in $Z$ by
$\oB^Z_\ep=\{z'\in Z:\rho(z',z)<\ep\}$. 
and
$\cB^Z_\ep=\{z'\in Z:\rho(z',z)\leq\ep\}$. 
Up to null events,
\begin{align*}
& \{\supp \mu_t\neq X_t\text{ for some }t\in(0,\infty)\}\\
& \subset \bigcup_{i\in\N}\{\mu_t(\oB^Z_{2(t-s)}(t,i))=0\text{ for some }t\in(0,\infty),s\in(0,t)\cap\Q\}\\
& = \bigcup_{i\in\N}\{|B(\Pi_{s,t},i)|=0\text{ for some }t\in(0,\infty),s\in(0,t)\cap\Q\}
\end{align*}
and we have a null event in the last line by Lemma~\ref{Pathw:lem:Pi-supp}. By Lemma~\ref{Pathw:lem:mu0}, $X_0$ is the closed support of $\mu_0$ a.\,s. This shows assertion~\ref{Pathw:item:thm:ld-CDI-supp}.

For $t\in\R_+$ and $n\in\N$, we denote the subset $\{t\}\times\sn$ of $Z$ by $X^n_t$, and more generally, for $M\subset\N$, we denote the subset $\{t\}\times M$ of $Z$ by $X^M_t$. We mostly omit `a.\,s.' in the following.

In the proof of assertion~\ref{Pathw:item:thm:ld-CDI-reg}, we begin with right continuity. Let $t\in\R_+$ and $\ep>0$. By construction, $X_u\subset \cB^Z_{u-t}(X_t)$ for all $u\geq t$. As $X_t$ is compact and as $\{t\}\times\N$ is dense in $X_t$, there exists $n\in\N$ such that
\[X_t \subset\cB^Z_\ep(X^n_t).\]
By condition~\reff{Pathw:eq:ass-eta}, we may choose $\delta\in(0,\ep)$ sufficiently small such that $\eta((t,t+\delta)\times\p^n)=0$. Then,
\[\cB^Z_\ep(X^n_t)\subset \cB^Z_{2\ep}(X^n_u)\subset \cB^Z_{2\ep}(X_u)\]
for all $u\in[t,t+\delta)$. Thus, we obtain the bound
\[\dH^Z(X_t,X_u)\leq 2\ep\]
for the Hausdorff distance $\dH^Z$ over $Z$. This proves right continuity of the map $t\mapsto X_t$ in $\dH^Z$.

We now turn to the left limits. We write $A_s(t,\N)=\{A_s(t,i):i\in\N\}$. Recall from Section \ref{Pathw:sec:ld-ext} the level $D_t(s,i)$ of the descendant of an individual $(s,i)$ at time $t\geq s$. Let $t\in(0,\infty)$ and $\ep\in(0,t)$. For $s\in[0,t)$, we denote the closure of $(s,t)\times\N$ in $Z$ by $X_{s,t}$. We define the closed subset
\[X_{t-}=\bigcap_{s\in(0,t)}X_{s,t}\]
of $Z$. By construction, $X_{t-}\supset X_t$. We claim that $\dH^Z(X_s,X_{t-})\to 0$ as $s\uparrow t$. From the definition of $X_{t-}$, we have $X_{t-}\subset \cB^Z_{t-s}(X_s)$ for all $s\in(0,t)$. Let $M=\bigcap_{s\in[t-\ep,t)}A_{t-\ep}(s,\N)$. As $\Xi\in\MCDI$, it holds $\#M<\infty$. That is, at time $t-$, the number of families of individuals that descend from the same ancestor at time $t-\ep$ is finite. As the map $(t-\ep/2,\infty)\to2^\N$, $s\mapsto A_{t-\ep}(s,\N)$ is non-increasing and by condition \reff{Pathw:eq:ass-eta-hat} piecewise constant, there exists $\delta\in(0,\ep)$ such that $A_{t-\ep}(s,\N)=M$ for all $s\in(t-\delta,t)$. For all $i\in M$ and $s,s'\in(t-\ep,t)$, the definitions of $M$ and $D_s(t-\ep,i)$ yield $D_s(t-\ep,i)<\infty$ and $\rho((s,D_s(t-\ep,i)),(s',D_{s'}(t-\ep,i)))=|s-s'|$,
hence the limit $x:=\lim_{s\uparrow t}(s,D_s(t-\ep,i))$ exists in the complete subspace $X_{s',t}$. Also note that $\rho(x,(t-\ep,i))=\ep$ and $x\in X_{t-}$. Thus, $X^{M}_{t-\ep}\subset\cB^Z_{\ep}(X_{t-})$,
\[X_s \subset  \cB^Z_{\ep}(X^{M}_{t-\ep})\subset\cB^Z_{2\ep}(X_{t-}),\]
and $\dH^Z(X_s,X_{t-})\leq 2\ep$ for all $s\in(t-\delta,t)$.

Now we show $X_t=X_{t-}$ for $t\in(0,\infty)\setminus \Thetaext$. Let $x\in X_{t-}$. Then there exists a sequence $((s_k,i_k):k\in\N)$ in $(0,t)\times\N\subset Z$ with $0<s_1<s_2<\ldots$ that converges to $x$. For each $k\in\N$, there exists $\ell\in\N$ such that $\rho((s_n,i_n),(s_\ell,i_\ell))<2(s_\ell-s_k)$ for all $n\geq\ell$. This implies $A_{s_k}(s_n,i_n)=A_{s_k}(s_\ell,i_\ell)$ for all $n\geq \ell$. We set $j_k=A_{s_k}(s_\ell,i_\ell)$. Then $j_k=A_{s_k}(s_n,j_n)$ for all $n\geq \ell\in\N$, hence $D_s(s_k,j_k)<\infty$ for all $s\in[s_k,t)$.
By our assumption on $t$, the definition \reff{Pathw:eq:Thetaext} of $\Thetaext$, and condition \reff{Pathw:eq:ass-eta-hat}, this implies $j'_k:=D_t(s_k,j_k)<\infty$ for all $k\in\N$. The sequence $((s_k,j_k):k\in\N)$ converges to $x$ as
\[\rho((s_k,j_k),x)
=\lim_{n\to\infty}\rho((s_k,j_k),(s_n,i_n))
\leq\lim_{n\to\infty}\rho((s_k,i_k),(s_n,i_n))
=\rho((s_k,i_k),x)\]
for all $k\in\N$. Also the sequence $((t,j'_k):k\in\N)$ converges to $x$ as
\[\rho((t,j'_k),(s_k,j_k))=t-s_k\]
for all $k\in\N$. This implies $x\in X_t$.

Now let $t\in\Thetaext$. Then by \reff{Pathw:eq:Thetaext}, there exists $\ep\in(0,t)$ and $i\in\N$ with $D_t(t-\ep,i)=\infty$ and $D_s(t-\ep,i)<\infty$ for all $s\in[t-\ep,t)$.
That is, the descendants of some ancestor at time $t-\ep$ die out at time $t$. The space $(X_t,\rho\wedge t)$ is ultrametric, cf.\ \cite{Sampl}*{Remark 5.2} and equation~\reff{Pathw:eq:def-rhot}.
Hence, a semi-metric $\rho^{(\ep)}$ on $X_t$ is given by $\rho^{(\ep)}=(\rho\wedge t-\ep)\vee 0$. We denote by $X^{(\ep)}_t$ the space obtained from $X_t$ by identifying elements with $\rho^{(\ep)}$-distance $0$. As $X_t$ is compact, the space $X^{(\ep)}_t$ is finite. Also the space $(X_{t-},\rho\wedge t)$ is ultrametric. Indeed, for each $x,y,z\in X_{t-}$ and each $s\in(0,t)$, there exist $x',y',z'\in X_s$ with $\rho(x,x')\leq t-s$, $\rho(y,y')\leq t-s$, and $\rho(z,z')\leq t-s$. This implies
\begin{align*}
(\rho\wedge t)(x,z)&\leq(\rho\wedge s)(x',z')+2(t-s)\\
&\leq\max\{(\rho\wedge s)(x',y'),(\rho\wedge s)(y',z')\}+2(t-s)\\
&\leq\max\{(\rho\wedge t)(x,y),(\rho\wedge t)(y,z)\}+4(t-s).
\end{align*}
Hence, a semi-metric $\rho^{(\ep)}$ on $X_{t-}$ is given by $\rho^{(\ep)}=(\rho\wedge t-\ep)\vee 0$. We denote by $X^{(\ep)}_{t-}$ the space obtained from $X_{t-}$ by identifying elements with $\rho^{(\ep)}$-distance $0$.
With $i$ as above, the limit $x:=\lim_{s\uparrow t}(s,D_s(t-\ep,i))$ exists in $X_{t-}$. As $A_{t-\ep}(t,j)\neq i$ for all $j\in\N$, it follows that $\rho((s,A_s(t,j)),(s,D_s(t-\ep,i)))\geq 2(s-t+\ep)$ for all $s\in(t-\ep,t)$. Taking the limit $s\uparrow t$, we obtain that $\rho((t,j),x)\geq 2\ep$. As $\{t\}\times\N$ is dense in $X_t$, it also follows that $\rho(x',x)\geq 2\ep$ for all $x'\in X_t$. As $X_t\subset X_{t-}$, this implies that the cardinality of $X^{(\ep)}_{t-}$ is greater than the cardinality of $X^{(\ep)}_t$. Hence, $X^{(\ep)}_{t-}$ and $X^{(\ep)}_t$ are not isometric. It follows that $X_{t-}$ and $X_t$ are not isometric. As a consequence, $\dH^Z(X_{t-},X_t)>0$, that is, the map $t\mapsto X_t$ is discontinuous in $\Thetaext$.
\end{proof}

\section{Proof of Lemmas~\ref{Pathw:lem:BC-bound} and \ref{Pathw:lem:BC-bound-dust}}
\label{Pathw:sec:unif-ld-proofs}
We work in the context of Section~\ref{Pathw:sec:unif-ld}, using also the definitions from Section~\ref{Pathw:sec:outline-proofs}. The proofs in the present section are adaptations of the proofs of Lemmas 3.4 and 3.5 of Donnelly and Kurtz \cite{DK99} and of Lemma 3.2 of Birkner et al.\,\cite{BBMST09}. We also mention Lemma 6.2 of Labbé \cite{Lab12}.

We define stochastic processes $(U(t),t\in\R_+)$ and $(\hat U(t),t\in\R_+)$ by
\begin{equation*}
U(t)=\int_{(0,t]\times \Delta}\left|x\right|_2^2\,\zeta_0(ds\;dx)
\end{equation*} 
and
\begin{equation*}
\hat U(t)=\int_{(0,t]\times \Delta}\left|x\right|_1\,\zeta_0(ds\;dx).
\end{equation*}
For $t\in\R_+$, the random variable $U(t)$ equals the sum of the squared asymptotic frequencies of the blocks that encode the reproduction events up to time $t$. The random variable $\hat U(t)$ equals the sum of the asymptotic frequencies of these blocks.
By the properties of the Poisson random measure $\zeta_0$,
\[\E[U(t)]=t\int_{\Delta}\left|x\right|_2^2\left|x\right|_2^{-2}\Xi_0(dx)<\infty,\]
hence the random variable $U(t)$ is a.\,s.\ finite. In case $\Xi\in\dust$, we also have that
\[\E[\hat U(t)]=t\int_{\Delta}\left|x\right|_1\left|x\right|_2^{-2}\Xi_0(dx)<\infty,\]
and $\hat U(t)$ is a.\,s.\ finite.
\begin{proof}[Proof of Lemma~\ref{Pathw:lem:BC-bound}]
We assume w.\,l.\,o.\,g.\ $\Xi(\Delta)>0$. Let $T,\ep,c>0$. For $\ell\in\N$, we set $\alpha^\ell_0=0$ and inductively for $k\in\N_0$
\[\alpha^\ell_{k+1}=
\inf\{t>\alpha^\ell_k:U(t)>U(\alpha^\ell_k)+\ell^{-4}\}
\wedge(\alpha^\ell_k+\ell^{-4}).\]
We also set $k_\ell=2\lceil(c+T)\ell^4\rceil$, then we have that $\P(\alpha^\ell_{k_\ell}\leq T,U(T)\leq c)=0$.

For an arbitrary integer $n_\ell\geq\ell$, we define for $k\in\N_0$
\[\beta^\ell_k=\inf\{t>\alpha^\ell_k: \eta((\alpha^\ell_k,t]\times\p^b)>0\},\]
\[\tilde\alpha^\ell_k=\inf\{t>\alpha^\ell_k:\left|X_{n_\ell}(t)-X_{n_\ell}(\alpha^\ell_k)\right|\geq 4\ep\}\wedge(\alpha^\ell_k+1),\]
and
\[\tilde\beta^\ell_k=\inf\{t>\beta^\ell_k:\left|X_{n_\ell}(t)-X_{n_\ell}(\beta^\ell_k)\right|\geq 4\ep\}\wedge(\beta^\ell_k+1).\]

Recall from Section \ref{Pathw:sec:two-step-eta} the point measure $\eta_0$. Each point $(t,\pi)$ of $\eta_0$ stands for a large (i.\,e.\ non-Kingman) reproduction event at time $t$ that is governed by the partition $\pi$. If $\Xi_0(\Delta)=0$, we have $\eta_0=0$ a.\,s.

We define the event $E^\ell$ as the intersection of the two events $E^\ell_1$ and $E^\ell_2$ that are defined as follows: We set
\[E^\ell_1=\bigcap_{k\in\N_0}\{\eta((\alpha^\ell_k\wedge (T+1),\alpha^\ell_{k+1}\wedge (T+1))\times\p^b)\leq 1\},\]
where we use the notation $(c,c')=\emptyset$ for $c'\leq c$. We define
\[E^\ell_2=\{\left|\pi\right|_2^2>\ell^{-4}\text{ for all }\pi\in\p^b\text{ with }\eta_0((0,T+1]\times\{\pi\})>0\},\]
where $\left|\pi\right|_2^2=\sum_{B\in\pi}|B|^2$.
Note that a.\,s., the asymptotic frequencies of the blocks $B$ of the partitions $\pi$ that form part of the points of $\eta_0$ exist.
This follows by Kingman's correspondence from the definition \reff{eq:def-eta0} of $\eta_0$ and as $\eta_0$ has a.\,s.\ at most countably many points.

As $\alpha^\ell_{k+1}\leq\alpha^\ell_k+\ell^{-4}$, the event $E^\ell_1$ occurs if $\ell^{-4}$ goes below the minimal distance between points of $\eta(\cdot\times\p^b)$ in $(0,T+1]$. This minimal distance is positive on the event of probability $1$ on which $\eta((0,T+1]\times\p^b)<\infty$ holds by condition \reff{Pathw:eq:ass-eta}.
A.\,s., also the event $E^\ell_2$ occurs for all sufficiently large $\ell$. Indeed, all points $(t,\pi)$ of $\eta_0$ satisfy $|\pi|_2^2>0$ a.\,s.\ by the definition of $\eta_0$ from $\zeta_0$ in Section \ref{Pathw:sec:two-step-eta} and as $\Xi_0\{0\}=0$ implies that for all points $(t,y)$ of $\zeta_0$, the first component $y_1$ of $y$ is positive a.\,s.
Hence a.\,s., $E^\ell$ occurs for all sufficiently large $\ell$.

Recall the strong Markov property of the process $J$ from~\reff{Pathw:eq:def-J}. For each $k\in\N_0$, the sequence $(Y^b_i(\alpha^\ell_k),i\in\N)$ is exchangeable by Corollary~\ref{Pathw:cor:exch-R} as $\alpha^\ell_k$ is $\zeta_0$-measurable.
The distance matrix
\[\I{\eta((\alpha^\ell_k,\tilde\alpha^\ell_k]\times\p^b)=0}\gamma_{b+1+n_\ell}(\rho_{\tilde\alpha^\ell_k})\]
is $(b+1+n_\ell,b)$-exchangeable by Lemma~\ref{Pathw:lem:exch-H} and the strong Markov property of $J$ at $\alpha^\ell_k$.
Hence, the vector
\[\left(\I{\eta((\alpha^\ell_k,\tilde\alpha^\ell_k]\times\p^b)=0}Y^b_i(\tilde\alpha^\ell_k),i\in[ n_\ell ]\right)\]
is exchangeable.

If $b\geq 2$ and $\Xi\{0\}>0$, then the distance matrix
\[\I{\zeta_0(\{\beta^\ell_k\}\times\Delta)=0}\gamma_{b+1+n_\ell}(\rho_{\beta^\ell_k})\]
is $(b+1+n_\ell,b+1)$-exchangeable by Lemma~\ref{Pathw:lem:exch-b} and the strong Markov property of $J$ at $\alpha^\ell_k$. Hence, the vector
\[(\I{\zeta_0(\{\beta^\ell_k\}\times\Delta)=0}Y^b_i(\beta^\ell_k),i\in n_\ell)\]
is exchangeable. The distance matrix
\[\I{\zeta_0(\{\beta^\ell_k\}\times\Delta)=0,\eta((\beta^\ell_k,\tilde\beta^\ell_k]\times\p^b)=0}\gamma_{b+1+n_\ell}(\rho_{\tilde\beta^\ell_k})\]
is $(b+1+n_\ell,b+1)$-exchangeable by Lemma~\ref{Pathw:lem:exch-H} and the strong Markov property of $J$ at $\beta^\ell_k$. Hence, the vector 
\[\left(\I{\zeta_0(\{\beta^\ell_k\}\times\Delta)=0,\eta((\beta^\ell_k,\tilde\beta^\ell_k]\times\p^b)=0}Y^b_i(\tilde\beta^\ell_k),i\in[ n_\ell]\right)\]
is exchangeable. If $b<2$ or $\Xi\{0\}=0$, then it suffices to work with the stopping times $\alpha^\ell_k$ and $\tilde\alpha^\ell_k$.

By Lemma A.2 in \cite{DK99}, there exists a number $\eta_\ep$
that depends only on $\ep$ (not on $n_\ell$) such that
\begin{align*}
&\P(|X_{n_\ell}(\alpha^\ell_k)-X_\ell(\alpha^\ell_k)|\geq \ep)\leq
2\expp{-\eta_\ep\ell},\\
&\P(|X_{n_\ell}(\tilde\alpha^\ell_k)-X_\ell(\tilde\alpha^\ell_k)|\geq \ep,\eta((\alpha^\ell_k,\tilde\alpha^\ell_k]\times\p^b)=0)	\leq
2\expp{-\eta_\ep\ell},\\
&\P(|X_{n_\ell}(\beta^\ell_k)-X_\ell(\beta^\ell_k)|\geq \ep,\zeta_0(\{\beta^\ell_k\}\times\Delta)=0)\leq
2\expp{-\eta_\ep\ell},\\
\text{and}\qquad
&\P(|X_{n_\ell}(\tilde\beta^\ell_k)-X_\ell(\tilde\beta^\ell_k)|\geq \ep,\zeta_0(\{\beta^\ell_k\}\times\Delta)=0,\eta((\beta^\ell_k,\tilde\beta^\ell_k]\times\p^b)=0)\leq
2\expp{-\eta_\ep\ell}.
\end{align*}

Let
\begin{align*}
H_k=
&\big|X_{n_\ell}(\alpha^\ell_k)-X_{\ell}(\alpha^\ell_k)\big|
\vee\big|X_{n_\ell}(\tilde\alpha^\ell_k\wedge\beta^\ell_k\wedge \alpha^\ell_{k+1})-X_{\ell}(\tilde\alpha^\ell_k\wedge\beta^\ell_k\wedge \alpha^\ell_{k+1})\big|\\
&\vee\big|X_{n_\ell}(\beta^\ell_k\wedge \alpha^\ell_{k+1})-X_{\ell}(\beta^\ell_k\wedge \alpha^\ell_{k+1})\big|
\vee\big|X_{n_\ell}(\tilde\beta^\ell_k\wedge \alpha^\ell_{k+1})-X_{\ell}(\tilde\beta^\ell_k\wedge \alpha^\ell_{k+1})\big|.
\end{align*}
As
\[\{\beta^\ell_k<\alpha^\ell_{k+1}\}\cap\{\alpha^\ell_k\leq T\}\cap E^\ell
\subset\{\zeta_0(\{\beta^\ell_k\}\times\Delta)=0\},\]
\[\{\tilde\beta^\ell_k<\alpha^\ell_{k+1}\}\cap\{\alpha^\ell_k\leq T\}\cap E^\ell
\subset\{\eta((\beta^\ell_k,\tilde\beta^\ell_k]\times\p^b)=0\},\]
and
\[\{\tilde\alpha^\ell_k<\beta^\ell_k\}\subset
\{\eta((\alpha^\ell_k,\tilde\alpha^\ell_k]\times\p^b)=0\}\]
up to null events for all $k\in\N_0$, the above implies
\begin{equation*}
\P(\max_{k< k_\ell:\,\alpha^\ell_k\leq T} H_k \geq\ep,E^\ell)\leq
16\lceil (c+T)\ell^4\rceil\expp{-\eta_\ep\ell}.
\end{equation*}

For $k\in\N_0$, we have that $\tilde\alpha^\ell_k\geq\beta^\ell_{k}\wedge\alpha^\ell_{k+1}$ a.\,s.\ on the event 
\[\{H_k<\ep\}\cap\left\{\sup_{t\in[\alpha^\ell_k,\beta^\ell_{k}\wedge\alpha^\ell_{k+1})}\left|X_\ell(t)-X_\ell(\alpha^\ell_k)\right|<\ep\right\}.\]
Indeed, the intersection of this event with $\{\tilde\alpha^\ell_k<\beta^\ell_{k}\wedge\alpha^\ell_{k+1}\}$ is a null event as it holds on this event that
\[|X_{n_\ell}(\tilde\alpha^\ell_k)-X_{n_\ell}(\alpha^\ell_k)|<3\ep,\]
whereas we have
\[|X_{n_\ell}(\tilde\alpha^\ell_k)-X_{n_\ell}(\alpha^\ell_k)|\geq 4\ep\quad\text{a.\,s.}\]
by definition of $\tilde\alpha^\ell_k$ and right continuity.

Similarly, $\tilde\beta^\ell_k\geq\alpha^\ell_{k+1}$ a.\,s.\ on the event 
\[\{H_k<\ep\}\cap \left\{\sup_{t\in[\beta^\ell_k,\alpha^\ell_{k+1})}\left|X_\ell(t)-X_\ell(\beta^\ell_k)\right|<\ep\right\}.\]

By Lemma~\ref{Pathw:lem:bound-finite} below and the Markov inequality,
\begin{align}
&\sum_{\ell\in\N}k_\ell\P(N^{b+1+\ell}(0,\alpha^\ell_1)>\ell \ep)\notag\\
&\leq \sum_{\ell\geq b+1}2\lceil (c+T)\ell^4\rceil\frac{1}{(\ell\ep)^4}\E[(N^{2\ell}(0,\alpha^\ell_1))^4]+\sum_{\ell=1}^b k_\ell
<\infty.\label{eq:bound-BC-1}
\end{align}
Using the strong Markov property of $J$ at $\alpha^\ell_k$ and the assumption \reff{eq:cond-newb} on $f$, we obtain
\begin{align}
&\P(\max_{k<k_\ell}\sup_{t\in[\alpha^\ell_k,\beta^\ell_{k}\wedge\alpha^\ell_{k+1})}\left|X_{\ell}(t)-X_{\ell}(\alpha^\ell_k)\right|>\ep)\notag\\
&\leq\sum_{k< k_\ell}\P(\sup_{t\in[\alpha^\ell_k,\beta^\ell_{k}\wedge\alpha^\ell_{k+1})}\left|X_{\ell}(t)-X_{\ell}(\alpha^\ell_k)\right|>\ep)\notag\\
&= k_\ell\P(\sup_{t\in[0,\beta^\ell_0\wedge\alpha^\ell_1)}\left|X_{\ell}(t)-X_{\ell}(\alpha^\ell_1)\right|>\ep)\notag\\
&\leq k_\ell \P(N^{b+1+\ell}(0,\beta^\ell_0\wedge\alpha^\ell_1)>\ell\ep)\notag\\
&\leq k_\ell \P(N^{b+1+\ell}(0,\alpha^\ell_1)>\ell\ep)\label{eq:bound-BC-2}
\end{align}
for all $\ell\in\N$. By \reff{eq:bound-BC-1} and \reff{eq:bound-BC-2},
\[\sum_{\ell\in\N}\P(\max_{k<k_\ell}\sup_{t\in[\alpha^\ell_k,\beta^\ell_{k}\wedge\alpha^\ell_{k+1})}\left|X_{\ell}(t)-X_{\ell}(\alpha^\ell_k)\right|>\ep)<\infty.\]
After replacing $[\alpha^\ell_k,\beta^\ell_{k}\wedge\alpha^\ell_{k+1})$ with $[\beta^\ell_k\wedge (T+1),\alpha^\ell_{k+1}\wedge (T+1))$ and intersecting with the event $E^\ell$ in \reff{eq:bound-BC-2}, the calculation from \reff{eq:bound-BC-1} and \reff{eq:bound-BC-2} also yields
\begin{align*}
&\P(\max_{k<k_\ell}\sup_{t\in[\beta^\ell_k\wedge (T+1),\alpha^\ell_{k+1}\wedge (T+1))}\left|X_{\ell}(t)-X_{\ell}(\beta^\ell_k)\right|>\ep,E^\ell)\\
&\leq k_\ell\P(N^{b+1+\ell}(0,\alpha^\ell_1)>\ell\ep).
\end{align*}
Indeed, we can use assumption \reff{eq:cond-newb} in the same way as in the third line of \reff{eq:bound-BC-2} because $\eta((\beta^\ell_0\wedge (T+1),\alpha^\ell_1\wedge (T+1))\times\p^b)=0$ on the event $E^\ell$. Again we deduce
\[\sum_{\ell\in\N}\P(\max_{k<k_\ell}\sup_{t\in[\beta^\ell_k\wedge (T+1),\alpha^\ell_{k+1}\wedge (T+1))}\left|X_{\ell}(t)-X_{\ell}(\beta^\ell_k)\right|>\ep,E^\ell)<\infty.\]

Altogether, it follows that there exist $\delta_\ell$ which do not depend on $n_\ell$ such that $\sum_{\ell=1}^\infty\delta_\ell<\infty$ and
\[\P(\sup_{t\in [0,T]}\left|X_{n_\ell}(t)-X_{\ell}(t)\right|>4\ep,U(T)\leq c,E^\ell)<\delta_\ell\]
for all $\ell\in\N$.

By Corollary~\ref{Pathw:cor:exch-R} and the de Finetti Theorem, there exists an event of probability $1$ on which the limits $X(t)=\lim_{n\to\infty} X_n(t)$ exist for all $t\in\Q_+$. Hence,
\begin{align*}
& \P(\sup_{t\in[0,T]\cap\Q}\left|X(t)-X_\ell(t)\right|> 4\ep,U(T)\leq c,E^\ell)\\
& = \P(\sup_{t\in[0,T]\cap\Q}\liminf_{n\to\infty}\left|X_n(t)-X_\ell(t)\right|> 4\ep,U(T)\leq c,E^\ell)\\
& \leq \lim_{n\to\infty}\P(\inf_{j\geq n} \sup_{t\in[0,T]\cap\Q}\left|X_j(t)-X_\ell(t)\right|> 4\ep,U(T)\leq c,E^\ell)\leq \delta_\ell.
\end{align*}
The Borel-Cantelli lemma allows to deduce that a.\,s.\ on the event $\{U(T)\leq c\}$,
\[\sup_{t\in[0,T]\cap\Q}\left|X(t)-X_\ell(t)\right|\leq 4\ep\]
for all sufficiently large $\ell$. Here we used that a.\,s., $E^\ell$ occurs for all sufficiently large $\ell$.
Hence, there exists a random integer $L$ such that 
\[|X_n(t)-X_\ell(t)|\leq|X_n(t)-X(t)|+|X_\ell(t)-X(t)|\leq 8\ep\]
for all $t\in[0,T]$ and $n\geq\ell\geq L$ a.\,s. on the event $\{U(T)\leq c\}$. Here we used right continuity of $X_n$ and $X_\ell$. It follows that
\[\P(\lim_{n,\ell\to\infty}\sup_{t\in[0,T]}|X_n(t)-X_\ell(t)|=0,U(T)\leq c)=1.\]
The assertion follows by letting $c$ tend to infinity.
\end{proof}
\begin{lem}
\label{Pathw:lem:bound-finite}
For $\ell\in\N$, let $\alpha^\ell_1$ be defined as in the proof of Lemma~\ref{Pathw:lem:BC-bound}. Then there exists a constant $C$ such that $\E[(N^{2\ell}(0,\alpha^\ell_1))^4]\leq C\ell^{-2}$ for all $\ell\in\N$.
\end{lem}
The proof extends the argument presented on p.\,44 in \cite{BBMST09} where additional assumptions on $\Xi$ are required to ensure that the process used there instead of $U(t)$ is finite.
\begin{proof}
First, let $x\in\Delta$. For $\ell\in\N$, let $(X_1,X_2,\ldots)$ have infinite multinomial distribution with parameters $(\ell,x_1,x_2,\ldots)$, that is, we may consider iid random variables $U_1,U_2,\ldots$ with uniform distribution on $[0,1]$ and set
\[X_i=\#\{j\in\sell: \sum_{k=1}^{i-1}x_k<U_j<\sum_{k=1}^ix_k\}\]
for $i\in\N$.
The infinite multinomial distribution appears in the context of $\Xi$-coalescents e.\,g.\ in \cite{M10}.
For a random partition $\pi$ in $\p$ with distribution $\kappa(x,\cdot)$, Kingman's correspondence implies that $b_{\ell}(\pi)$ and $\lim_{n\to\infty}\sum_{i=1}^n(X_i-1)_+$ are equal in distribution. Here we write $(x-1)_+=\max\{x-1,0\}$. We use the inequalities $[(x-1)_+]^2\leq x^{(2)}$ and $[(x-1)_+]^4\leq 3x^{(4)}+3x^{(2)}$ for $x\in\N_0$, where $x^{(k)}=x!/(x-k)!$. Inserting also mixed factorial moments of multinomial distributions, we obtain a constant $C'$ such that
\begin{align*}
\E\left[\left(\sum_{i=1}^n(X_i-1)_+\right)^2\right]
\leq\sum_{\substack{i,j=1\\i\neq j}}^n\E\left[X_i^{(2)}X_j^{(2)}\right]
+\sum_{i=1}^n\E\left[X_i^{(2)}\right]
\leq \ell^4|x|_2^4+\ell^2|x|_2^2
\end{align*}
and
\begin{align*}
\E\left[\left(\sum_{i=1}^n(X_i-1)_+\right)^4\right]
\leq C'\,\left( \ell^8|x|_2^8+\ell^6|x|_2^6+\ell^6|x|_4^4\,|x|_2^2+\ell^4|x|_2^4+\ell^4|x|_4^4+\ell^2|x|_2^2\right)
\end{align*}
for all $n\in\N$. Taking the limit $n\to\infty$ on the left-hand side, we obtain upper bounds for $\kappa(x,\cdot)b_\ell^2$ and $\kappa(x,\cdot)b_\ell^4$.

Let
\[N_0^\ell(I)=\int_{I\times\p}b_\ell(\pi)\eta_0(ds\;d\pi)\]
for each interval $I\subset\R_+$.  The random variable $N_0^\ell(I)$ is the number of newborn particles on the first $\ell$ levels in the large reproduction events in the time interval $I$. The random variable $N^{2\ell}(0,\alpha^\ell_1)-N^{2\ell}_0(0,\alpha^\ell_1)=\eta_{\rm K}((0,\alpha^\ell_1)\times\p^{2\ell})$ is stochastically bounded from above by a Poisson random variable with mean $\Xi\{0\}\ell^{-4}\binom{2\ell}{2}$.

Here we only show $\E[(N_0^{2\ell}(0,\alpha^\ell_1))^4]\leq C''\ell^{-2}$ for an appropriate constant $C''$ and all $\ell\in\N$. W.\,l.\,o.\,g.\ we assume $\Xi_0(\Delta)>0$. Recall $((t^i,y^i,\pi^i),i\in\N)$ from Section~\ref{Pathw:sec:outline-proofs}. We have
\begin{align*}
&\E[(N_0^{2\ell}(0,\alpha^\ell_1))^4|\zeta_0]\\
=&\E[\sum_{\substack{i_1,\ldots,i_4\in\N:\\t^{i_1},\ldots,t^{i_4}<\alpha^\ell_1}}b_{2\ell}(\pi^{i_1})\cdots b_{2\ell}(\pi^{i_4})|((t^i,y^i),i\in\N)]\\
\leq &C'''\big[(\sum_{i\in\N:\,t^i<\alpha^\ell_1}\kappa(y^i,\cdot)b_{2\ell})^4\\
&+\sum_{i\in\N:\,t^i<\alpha^\ell_1}\kappa(y^i,\cdot)b_{2\ell}^2
\,(\sum_{j\in\N:\,t^j<\alpha^\ell_1}\kappa(y^j,\cdot)b_{2\ell})^2\\
&+(\sum_{i\in\N:\,t^i<\alpha^\ell_1}\kappa(y^i,\cdot)b_{2\ell}^2)^2\\
&+\sum_{i\in\N:\,t^i<\alpha^\ell_1}\kappa(y^i,\cdot)b_{2\ell}^3
\,\sum_{j\in\N:\,t^j<\alpha^\ell_1}\kappa(y^j,\cdot)b_{2\ell}\\
&+\sum_{i\in\N:\,t^i<\alpha^\ell_1}\kappa(y^i,\cdot)b_{2\ell}^4\big]\quad\text{a.\,s.}
\end{align*}
for a combinatorial constant $C'''$. Now we estimate $b_{2\ell}\leq b_{2\ell}^2$ and $b_{2\ell}^3\leq b_{2\ell}^4$, and we insert the bounds for $\kappa(y^i,\cdot)b_{2\ell}^2$ and $\kappa(y^i,\cdot)b_{2\ell}^4$. From the definitions of $(U(t),t\in\R_+)$ and $\alpha^\ell_1$, we have
\[\sum_{i\in\N:\,t^i<\alpha^\ell_1}\left|y^i\right|_2^2\leq\ell^{-4}\]
which yields the assertion.
\end{proof}
\begin{proof}[Proof of Lemma~\ref{Pathw:lem:BC-bound-dust}]
The proof is analogous to the proof of Lemma~\ref{Pathw:lem:BC-bound}. Again we fix $T,c,\ep>0$ and assume $\Xi(\Delta)>0$. We work only with the stopping times $\alpha^\ell_k$ and $\tilde\alpha^\ell_k$ which we define for $\ell\in\N$, $k\in\N_0$, and some arbitrarily large integers $n_\ell\geq\ell$ as follows:
\[\alpha^\ell_0=0,\]
\[\tilde\alpha^\ell_k=\inf\{t>\alpha^\ell_k:\left|X_{n_\ell}(t)-X_{n_\ell}(\alpha^\ell_k)\right|\geq 4\ep\}\wedge(\alpha^\ell_k+1),\]
and
\[\alpha^\ell_{k+1}=\inf\{t>\alpha^\ell_k:\hat U(t)>\hat U(\alpha^\ell_k)+\ell^{-1}\}.\]
We set $k_\ell=2\lceil c\,\ell\rceil$. Then we have that $\P(\alpha_{k_\ell}\leq T,\hat U(T)\leq c)=0$. Let
\[E^\ell=\{\left|\pi\right|_2^2>\ell^{-1}\text{ for all }\pi\in\hat\p^b\text{ with }\eta_0((0,T+1]\times\{\pi\})>0\}.\]
A.\,s., the event $E^\ell$ occurs for all sufficiently large $\ell$ as $\eta((0,T+1]\times\hat\p^b)<\infty$ by \reff{Pathw:eq:ass-eta-hat} and as $\left|\pi\right|_2^2>0$ for all points $(t,\pi)$ of $\eta_0$. A.\,s.\ on $E^\ell$, no reproduction events that are characterized by a partition in $\hat\p^b$ occur in the open time intervals $(\alpha^\ell_k, \alpha^\ell_{k+1})$ as long as $\alpha^\ell_{k+1}<T+1$, that is, $\eta((\alpha^\ell_k,\alpha^\ell_{k+1})\times\hat\p^b)=0$ for such $k$. Here we use that $\eta=\eta_0$ a.\,s.\ by our assumption that $\Xi\in\dust$.

We replace the distance matrices $(\rho_t,t\in\R_+)$ from the proof of Lemma \ref{Pathw:lem:BC-bound} by the marked distance matrices $((r_t,v_t),t\in\R_+)$ to obtain exchangeability for the sequence
$(Y^b_i(\alpha^\ell_k),i\in\N)$ and the vector
\[(\I{\eta((\alpha^\ell_k,\tilde\alpha^\ell_k]\times\hat\p^b)=0}Y^b_i(\tilde\alpha^\ell_k),i\in[n_\ell]).\]
We set
\[H_k=
\big|X_{n_\ell}(\alpha^\ell_k)-X_{\ell}(\alpha^\ell_k)\big|
\vee\big|X_{n_\ell}(\tilde\alpha^\ell_k\wedge\alpha^\ell_{k+1})-X_{\ell}(\tilde\alpha^\ell_k\wedge \alpha^\ell_{k+1})\big|\]
for $k\in\N_0$. Applying Lemma A.2 of \cite{DK99} as in the proof of Lemma \ref{Pathw:lem:BC-bound}, and using that
\[\{\tilde\alpha^\ell_k<\alpha^\ell_{k+1}\}\cap\{\alpha^\ell_k\leq T\}\cap E^\ell
\subset\{\eta((\alpha^\ell_k,\tilde\alpha^\ell_k]\times\hat\p^b)=0\}\]
for all $k\in\N_0$, we obtain
\begin{equation*}
\P(\max_{k< k_\ell:\,\alpha^\ell_k\leq T} H_k \geq\ep,E^\ell)\leq
8\lceil c\ell\rceil\expp{-\eta_\ep\ell}.
\end{equation*}

We claim that
\begin{equation}
\label{eq:claim-bound-BC-dust}
\sum_{\ell\in\N}k_\ell\P(\hat N^{b+\ell}(0,\alpha^\ell_1)>\ell\ep)<\infty
\end{equation}
which we prove further below.
From the strong Markov property of $J$ at $\alpha^\ell_k$ and by the assumption \reff{eq:cond-repr} on $f$, we obtain
\begin{align*}
&\P(\max_{k<k_\ell}\sup_{t\in[\alpha^\ell_k,\alpha^\ell_{k+1})}\left|X_{\ell}(t)-X_{\ell}(\alpha^\ell_k)\right|>\ep, E^\ell)\\
&\leq k_\ell \P(\hat N^{b+\ell}(0,\alpha^\ell_1)>\ell\ep,E^\ell)
\end{align*}
for all $\ell\in\N$, and the claim \reff{eq:claim-bound-BC-dust} implies
\begin{equation}
\label{eq:bound-BC-dust-1}
\sum_{\ell\in\N}\P(\max_{k<k_\ell}\sup_{t\in[\alpha^\ell_k,\alpha^\ell_{k+1})}\left|X_{\ell}(t)-X_{\ell}(\alpha^\ell_k)\right|>\ep,E^\ell)<\infty.
\end{equation}

As in the proof of Lemma~\ref{Pathw:lem:BC-bound}, on the event
\[\{H_k<\ep\}\cap\left\{\sup_{t\in[\alpha^\ell_k,\alpha^\ell_{k+1})}\left|X_\ell(t)-X_\ell(\alpha^\ell_k)\right|<\ep\right\},\]
we have $\tilde\alpha^\ell_k\geq\alpha^\ell_{k+1}$ a.\,s., hence \reff{eq:bound-BC-dust-1} yields
\[\P(\sup_{t\in [0,T]}\left|X_{n_\ell}(t)-X_{\ell}(t)\right|>4\ep,U(T)\leq c,E^\ell)<\delta_\ell\]
for a summable sequence $(\delta_\ell)$ that does not depend on the $n_\ell$, and the assertion follows as in the proof of Lemma~\ref{Pathw:lem:BC-bound}.

It remains to prove the claim \reff{eq:claim-bound-BC-dust}. By the Markov inequality, we have
\[\P(\hat N^{b+\ell}(0,\alpha^\ell_1)>\ell\ep)
\leq\expp{-\ell\ep}\E[\E[\exp(\hat N^{b+\ell}(0,\alpha^\ell_1))|\zeta_0]].\]
We show (similarly to the proof of Lemma \ref{Pathw:lem:bound-finite}) that
$\sup_{\ell\in\N}\E[\exp(\hat N^{2\ell}(0,\alpha^\ell_1))|\zeta_0]$
is bounded, this will imply the claim.

For $x\in\Delta$ and a random partition $\pi$ with distribution $\kappa(x,\cdot)$, the random variable $\hat b_{2\ell}(\pi)$ is binomially distributed with parameters $2\ell$ and $|x|_1$. Using monotone convergence, conditional independence, and inserting moment generating functions of binomial distributions, we obtain
\begin{align*}
&\E[\exp(\hat N^{2\ell}(0,\alpha^\ell_1))|\zeta_0]\\
=&\E[\exp(\sum_{i\in\N:\,t^i<\alpha^\ell_1}\hat b_{2\ell}(\pi^i))|((t^i,y^i),i\in\N)]\\
=&\lim_{n\to\infty}\prod_{i\in\sn:\,t^i<\alpha^\ell_1}\E[\exp(\hat b_{2\ell}(\pi^i))|((t^i,y^i),i\in\N)]\\
=&\lim_{n\to\infty}\prod_{i\in\sn:\,t^i<\alpha^\ell_1}(1-|y^i|_1+|y^i|_1\e )^{2\ell}\\
=&\exp(\sum_{i\in\N:\,t^i<\alpha^\ell_1}2\ell\log(1+|y^i|_1(\e-1) ))\\
\leq&\exp(2\ell\sum_{i\in\N:\,t^i<\alpha^\ell_1}|y^i|_1(\e -1))\leq \exp(2(\e-1))\quad\text{a.\,s.}
\end{align*}
The last inequality follows from the definitions of $(\hat U(t),t\in\R_+)$ and $\alpha^\ell_1$.
\end{proof}

\section*{List of notation}
\sectionmark{List of notation}
Here we collect notation that is used globally in the article.
\small{
\hparagraph{Miscellaneous}
$\R_+=[0,\infty)$, $\N=\{1,2,3,\ldots\}$, $\N_0=\N\cup\{0\}$, $\sn=\{1,\ldots,n\}$ for $n\in\N$, $[0]=\emptyset$\\
$\dP^X$, $\dH^X$: Prohorov metric and Hausdorff distance over $X$\\
$\supp\mu$: closed support of the measure $\mu$\\
$\varphi(\mu)=\mu\circ\varphi^{-1}$: pushforward measure under a measurable function $\varphi$\\ 
$\gamma_n$: restriction map in various contexts (p.\,\pageref{Pathw:not:gamma-p}/l.\,21, p.\,\pageref{Pathw:not:gamma-outl}/l.\,-6)\\
$\oB^X_\ep(x)=\{y\in X: d(x,y)<\ep\}$, $\cB^X_\ep(x)=\{y\in X: d(x,y)\leq\ep\}$: balls in a metric space $(X,d)$

\hparagraph{(Marked) distance matrices}
$\Dd$: space of semi-metrics on $\N$ (p.\,\pageref{Pathw:not:Dd}/l.\,10)\\
$\Dd_n$: space of semi-metrics on $[n]$ (p.\,\pageref{Pathw:not:Ddn}/l.\,-9)\\
$\hat\Dd$: space of decomposed semi-metrics on $\N$ (p.\,\pageref{Pathw:not:hatDd}/l.\,-2)\\
$\hat\Dd_n$: space of decomposed semi-metrics on $[n]$ (p.\,\pageref{Pathw:not:hatDdn}/l.\,-9)

\hparagraph{Partitions and semi-partitions}
$\p$: Set of partitions of $\N$ (p.\,\pageref{Pathw:not:p}/l.\,19)\\
$K_{i,j}$: partition of $\N$ that contains only $\{i,j\}$ and singleton blocks (p.\,\pageref{Pathw:not:Kij}/l.\,-2)\\
$\p_n$: Set of partitions of $[n]$, associated transformations (p.\,\pageref{Pathw:not:pntransf}/l.\,-5)\\
$\p^n$: Set of partitions of $\N$ in which the first $n$ integers are not all in different blocks, equation \reff{Pathw:eq:def-pn}\\
$\hat\p^n$: Set of partitions of $\N$ in which the first $n$ integers are not all in singleton blocks, (p.\,\pageref{Pathw:not:hatpn}/l.\,13)\\
$B(\pi,i)$: block of the partition $\pi\in\p$ that contains $i\in\N$\\
$\pi(i)=k$ such that $i$ is in the $k$-th block of $\pi$ (p.\,\pageref{Pathw:not:pii}/l.\,-7)\\
$M(\pi)$: set of the minimal elements of the blocks of $\pi\in\p$ (p.\,\pageref{Pathw:not:Mpi}/l.\,-3)\\
$|B|_n=n^{-1}(\#B\cap[n])$, $B=\lim_{n\to\infty}|B|_n$ for $B\subset\N$, $n\in\N$: relative and asymptotic frequency\\
$\#\pi$: number of blocks of a partition $\pi$\\
$\S_n$ set of semi-partitions of $[n]$, associated transformations (p.\,\pageref{Pathw:not:Sntransf}/l.\,5)\\
$\Delta=\{x=(x_1,x_2,\ldots):x_1\geq x_2\geq\ldots\geq 0,|x|_1\leq 1\}$\\
$|x|_p=\left(\sum_{i\in\N}x_i^p\right)^{1/p}$ for $x\in\Delta$\\
$\kappa(x,\cdot)$: paintbox distribution associated with $\Delta$ (p.\,\pageref{Pathw:not:kappa}/l.\,6)

\hparagraph{Genealogy in the lookdown model}
$\eta$: point measure on $(0,\infty)\times\p$ that encodes the reproduction events (p.\,\pageref{Pathw:not:eta-det}/l.\,-14, p.\,\pageref{Pathw:not:eta-Poisson}/l.\,3, p.\,\pageref{Pathw:not:eta-twostep}/l.\,6)\\
$\eta_0$, $\eta_{\rm K}$: restrictions of $\eta$ to large and binary reproduction events, respectively (p.\,\pageref{Pathw:not:eta0etaK}/l.\,6)\\
$\zeta_0$ point measure on $(0,\infty)\times\Delta$ that encodes the family sizes in the large reproduction events (p.\,\pageref{Pathw:not:xi0}/l.\,12)\\
$(t,i)$: individual on level $i$ at time $t$ (p.\,\pageref{Pathw:not:ti}/l.\,4)\\
$A_s(t,i)$: level of the ancestor at time $s$ of $(t,i)$ (p.\,\pageref{Pathw:not:ancestor}/l.\,6)\\
$(Z,\rho)$: lookdown space with genealogical distance (p.\,\pageref{Pathw:not:Z}/l.\,22)\\
$\rho_t(i,j)=\rho((t,i),(t,j))$ (equation \reff{Pathw:eq:def-rhot})\\
$(\hat Z,\rho)$: extended lookdown space with genealogical distance (p.\,\pageref{Pathw:not:hatZ}/l.\,-11)\\
$X_t$ closure of $\{t\}\times\N$ in $Z$ (p.\,\pageref{Pathw:not:Xt}/l.\,8)\\
$X_n(t)$: relative frequencies in Lemmas \ref{Pathw:lem:BC-bound} and \ref{Pathw:lem:BC-bound-dust}\\
$z(t,i)$: parent of $(t,i)$ (p.\,\pageref{Pathw:not:zti}/l.\,24)\\
$v_t(i)=\rho((t,i),z(t,i))$ (p.\,\pageref{Pathw:not:vti}/l.\,3)\\
$r_t(i,j)=\rho(z(t,i),z(t,j))$ (equation \reff{Pathw:eq:rho-r-v})\\
$D_t(s,i)$: lowest level at time $t$ of a descendant of $(s,i)$ (p.\,\pageref{Pathw:not:Dtsi}/l.\,-12)\\
$\tau_{s,i}$: extinction time of $(s,i)$ (p.\,\pageref{Pathw:not:tausi}/l.\,-7)\\
$\Thetaext=\{\tau_{s,i}:s\in\R_+,i\in\N\}$ (equation \reff{Pathw:eq:Thetaext})\\
$\Theta_0$: set of large reproduction times (equation \reff{Pathw:eq:Theta0})\\
$\mu^n_t$: uniform measure on $(t,i)$, $i\in[n]$ (p.\,\pageref{Pathw:not:munt}/l.\,-10)\\
$m^n_t$: uniform measure on $(z(t,i),v_t(i))$, $i\in[n]$ (equation \reff{eq:def-mt})\\
$\mu_t$, $m_t$: weak limits of $\mu^n_t$ and $m_t$, respectively (p.\,\pageref{Pathw:not:mut}/l.\,-9, p.\,\pageref{Pathw:not:mt}/l.\,18)\\
$(\Pi_{s,t},0\leq s\leq t)$: flow of partitions (p.\,\pageref{Pathw:not:Pist}/l.\,-16)\\
$\Pi_t$: partition of siblings: (p.\,\pageref{Pathw:not:Pit}/l.\,14)\\
$C^{a,\ep,I}_t$: subsets of individuals at time $t$ (p.\,\pageref{Pathw:not:CaepIt}/l.\,2)\\
$\Pi^{a,\ep}_t$: like $\Pi_t$ but individuals that are in a common $C^{a,\ep,I}_t$ are also in the same block (p.\,\pageref{Pathw:not:Piaept}/l.\,10)

\hparagraph{Characteristic measures}
$\Xi=\Xi_0+\Xi\{0\}\delta_0$ (p.\,\pageref{Pathw:not:Xi0}/l.\,-3)\\
$H_\Xi$: characteristic measure of $\eta$ (p.\,\pageref{Pathw:not:HXi}/l.\,-1)\\
$\dust$, $\nd$, $\MCDI$: Sets of finite measures on $\Delta$ with and without dust, and with the coming down from infinity property, respectively (p.\,\pageref{Pathw:not:nd}/l.\,9, p.\,\pageref{Pathw:not:MCDI}/l.\,22)
}

\bibliography{C:/Users/sg/Documents/Bib/diss}
\bibliographystyle{plain}

\normalsize
\paragraph{Acknowledgments.}
This work is part of the author's PhD thesis.
The author thanks Götz Kersting and Anton Wakolbinger for their valuable advice.
He also thanks them and the referees for very helpful comments.
Partial support from the DFG Priority Programme 1590 ``Probabilistic Structures in Evolution'' is acknowledged. In 2017/2018, the author is supported by a postdoctoral fellowship of the Minerva Foundation.
\end{document}